\documentclass[reqno]{amsart}

\usepackage{amsthm}
\usepackage{amssymb}
\usepackage{amsmath}
\usepackage{amscd}
\usepackage{graphicx}
\usepackage{latexsym}
\usepackage{enumerate}

\newcommand{\g}{{\mathfrak g}}
\newcommand{\fg}{{\mathfrak g}}
\newcommand{\fh}{{\mathfrak h}}

\newcommand{\CC}{{\mathbb C}}
\newcommand{\RR}{{\mathbb R}}

\newcommand{\ZZ}{{\mathbb Z}}

\newcommand{\LL}{{\mathbb L}}
\newcommand{\PP}{{\mathbb P}}
\newcommand{\K}{{\mathcal K}}

\newcommand{\fk}{{\mathfrak k}}
\newcommand{\DD}{{\mathcal D}}
\newcommand{\MM}{{\mathcal M}}
\newcommand{\U}{{\mathcal U}}
\renewcommand{\P}{{\mathcal P}}
\renewcommand{\SS}{{\mathbb S}}
\theoremstyle{definition}
\newtheorem*{definition*}{Definition}
\newtheorem{theorem}{Theorem}[subsection]
\newtheorem{definition}{Definition}[subsection]
\newtheorem*{remark*}{Remark}
\newtheorem*{remarks*}{Remarks}
\newtheorem*{theorem*}{Theorem}
\newtheorem{example}{Example}[subsection]
\newtheorem{example*}{Example}
\newtheorem{lemma}{Lemma}[section]
\newtheorem{corollary}{Corollary}[subsection]
\newtheorem{proposition}{Proposition}[subsection]

\numberwithin{equation}{section}

\gdef\myletter{}

\let\savetheequation\theequation
\def\theequation{\savetheequation\myletter}

\begin{document}

\title{One-skeleta, Betti numbers and equivariant cohomology}
\author[V. Guillemin]{V. Guillemin\footnotemark {*}}
\thanks{* Supported by NSF grant DMS 890771.}
\address{Department of Mathematics, MIT, Cambridge, MA 02139}
\email{vwg@@math.mit.edu}
\author[C. Zara]{C. Zara}
\address{Department of Mathematics, MIT, Cambridge, MA 02139}
\email{czara@math.mit.edu}

\begin{abstract}
The one-skeleton of a $G$-manifold $M$ is the set of points $p \in M$ 
where $\dim G_p \geq \dim G -1$; and $M$ is a GKM manifold if the 
dimension of this one-skeleton is 2. Goresky, Kottwitz and MacPherson 
show that for such a manifold this one-skeleton has the structure of a 
``labeled'' graph, $(\Gamma, \alpha)$, and that the equivariant cohomology 
ring of $M$ is isomorphic to the ``cohomology ring'' of this graph. Hence, 
if $M$ is symplectic, one can show that this ring is a free module 
over the symmetric algebra $\SS(\fg^*)$, with $b_{2i}(\Gamma)$ generators 
in dimension $2i$, $b_{2i}(\Gamma)$ being the ``combinatorial'' $2i$-th 
Betti number of $\Gamma$. In this article we show that this 
``topological'' result is , in fact, a combinatorial result about graphs.
\end{abstract}

\maketitle

\section*{\bf Introduction}
\label{sec:intro}

Let $G$ be a commutative, compact, connected, $n$-dimensional Lie
group, $\fg$ its Lie algebra, $M$ a compact $2d$-dimensional
manifold and $\tau : G \times M \to M$ a faithful action of $G$
on $M$.  We say that $M$ is a \emph{GKM manifold} 
if it has the following properties:

\begin{enumerate}
\item %%1
  $M^G$ is finite.

\item %%2
  $M$ possesses a $G$-invariant almost-complex structure.

\item %%3
  For every $p \in M^G$, the weights
  \begin{equation}
    \label{eq:1.0}
    \alpha_{i,p} \in \fg^* \, , \, i=1, \ldots ,d \, ,
  \end{equation}
of the isotropy representation of $G$ on $T_pM$ are pairwise
linearly independent.
\end{enumerate}
There is an alternate way of formulating this third condition:
Let $M$ be a $G$-manifold which satisfies the first two
conditions, and define the \emph{one-skeleton} of $M$ to be the
set of points, $p \in M$ with $\dim G_p \geq n-1$.  Then $M$
satisfies the third condition if and only if its one-skeleton
consists of $G$-invariant submanifolds which are fixed point free
and $G$-invariant embedded
$2$-spheres, each of which contains
exactly two fixed points. Thus the combinatorial structure of this
one-skeleton is that of a graph $\Gamma$ having the fixed
points of $G$ as vertices and these $2$-spheres as edges.
As we will see in the next section, $\Gamma$ is a
regular graph:  each vertex is the point of intersection of
exactly $d$ edges. 
Moreover, the action of $G$ on $M$ gives one a labeling of the oriented 
edges of $\Gamma$ by one-dimensional representations of $G$. 
Namely, to each oriented
edge, $e$, one can assign the isotropy representation, $\chi_e$, of $G$ 
on the tangent space at the ``north pole'' of the corresponding $S^2$, 
the ``north pole'' corresponding to the initial vertex of $e$. 
Thus one has a map
$$\alpha : E_{\Gamma} \to \fg^*$$
from the set of oriented edges of $\Gamma$ to $\fg^*$, which assigns 
to each oriented edge, $e$, the weight, $\alpha_e$, 
of the representation $\chi_e$. 
We will refer to the pair $(\Gamma, \alpha)$ as the 
\emph{GKM one-skeleton} associated to $M$.

A beautiful result of Goresky-Kottwitz-MacPherson asserts that 
if $M$ is equivariantly formal, the 
equivariant cohomology ring, $H_G(M)$, can be reconstructed from the 
GKM one-skeleton. More explicitly, let $V_{\Gamma}$ be the vertices of 
$\Gamma$ and $H(\Gamma, \alpha)$ the set of all maps, 
$f: V_{\Gamma} \to \SS(\fg^*)$, which satisfy the compatibility condition
\begin{equation}
\label{eq:n1.1}
f(p) - f(q) = 0 \mod{\alpha_e}
\end{equation}
for every pair of vertices $p$ and $q$, and every edge, $e$, 
joining $p$ and $q$. Then the GKM theorem asserts 
\begin{equation}
\label{eq:n1.2}
H_G(M) \simeq H(\Gamma, \alpha).
\end{equation}
One interesting implication of this theorem is that one can prove, 
by topology, combinatorial results about $H(\Gamma, \alpha)$. For instance, 
suppose that the action, $\tau$, of $G$ on $M$  is Hamiltonian. Then, by a 
theorem of Kirwan, $M$ is equivariantly formal; so the GKM theorem applies 
to $M$. Moreover, from the constant map,, $\gamma : M \to pt$, one gets 
a map, $\gamma^* : H_G(pt) \to H_G(M)$, and since $H_G(pt) = \SS(\fg^*)$, 
this map makes $H_G(M)$ into an $\SS(\fg^*)$-module. If $\tau$ is 
Hamiltonian, Kirwan proves that $H_G(M)$ is a free $\SS(\fg^*)$-module 
with $b_{2i}(M)$ generators in dimension $2i$, $b_{2i}(M)$ being the $2i$-th 
Betti number of $M$. In addition, using Morse theory, one can compute 
these Betti numbers directly from the GKM one-skeleton, $(\Gamma, \alpha)$, 
as follows: 
Fix $\xi \in \fg$ with $\alpha_e(\xi) \neq 0$ for all $e \in E_{\Gamma}$ 
and let $b_{2i}(\Gamma)$ be the number of 
vertices, $p$, for which there are exactly $i$ oriented edges, $e$, 
with \emph{initial} vertex $p$, such that $\alpha_e(\xi) < 0$. Then 
\begin{equation}
\label{eq:n1.3}
b_{2i}(\Gamma) = b_{2i}(M).
\end{equation}
Thus, by topology, one proves
\begin{theorem}
\label{th:n1}
$H(\Gamma, \alpha)$ is a free $\SS(\fg^*)$-module. Moreover
\begin{equation}
\label{eq:n1.4}
H(\Gamma, \alpha) \otimes_{\SS(\fg^*)} \CC
\end{equation}
is a finite dimensional graded ring, its $2i$-th graded component 
being of dimension $b_{2i}(\Gamma)$.
\end{theorem}

There are a number of other theorems about the structure of 
$H(\Gamma, \alpha)$ which can be proved ``by topology''. For instance, 
using equivariant Morse theory, one can write down a canonical set of 
generators of $H(\Gamma, \alpha)$, and if $\tau$ is Hamiltonian, one can, 
by methods of Kirwan (\cite{Ki}) prove a number of interesting facts 
about subrings and quotient rings of $H(\Gamma, \alpha)$. 
(See, for instance, \cite{TW}.)

The question we want to explore in this paper is: \emph{Can one prove 
these ``topological'' results about $H(\Gamma, \alpha)$ purely by 
combinatorial methods ?} In other words, 
\emph{are these theorems combinatorial theorems about graphs in disguise ?} 
Two types of GKM manifolds for which this question has a positive answer 
are toric varieties and flag varieties. For toric varieties $H_G(M)$ is 
the Stanley-Reisner ring of the moment polytope of $M$, and for flag 
varieties, $H_G(M)$ is the ring of ``double Schubert polynomials''; and, 
in these cases, the theorems above follow from combinatorial theorems 
about poset cohomology, root systems, Hecke algebras \emph{et al.} 
(See \cite{Bi}, \cite{BH}, \cite{Fu1}, \cite{Fu2}, \cite{Hu}, \cite{LS},
 \cite{St}) Therefore the question above is part of a more open-ended 
question: \emph{Are there analogues of some of these combinatorial 
theorems for GKM manifolds in general ?}

The interplay between graphs and GKM manifolds may have some 
interesting applications in 
graph theory \emph{per se}. We will describe one example of such 
an application: Let $\Delta$ be a convex polytope in 
$\RR^n$ and let $\Gamma$ be its \emph{one-skeleton}, \emph{i.e.}
the graph consisting of the vertices and edges 
of $\Delta$. Then, just as above, the oriented edges of $\Gamma$ 
have a natural labeling: to each oriented edge, $e$, one can assign the edge 
vector, $\alpha_e = q-p$, $p$ and $q$ being the initial and the terminal 
vertices of $e$. In analogy with the case of GKM manifolds, we will call
the pair $(\Gamma, \alpha)$ the \emph{GKM one-skeleton} of $\Delta$.
A problem of interest to combinatorists (see, for 
instance, \cite{CW}) is how to deform $\Delta$ so that the directions of 
its edges are unchanged. In particular, how many such deformations are 
there ? GKM theory suggests an answer: If $\Delta$ is a simple polytope 
and its edge directions are rational, it is the moment polytope of a 
toric variety, $M$; and the number of ways in which one can deform 
$\Delta$ without changing its edge directions is equal, by Delzant's 
theorem (\cite{Del}), to the number of ways in which one can deform the 
symplectic structure of $M$, \emph{i.e.} is equal to $\dim H^2(M)$, or, 
alternatively, by \eqref{eq:n1.2}, is equal to $b_2(\Gamma)$. We will 
show in section \ref{sec:polytopes} that this result is true not just 
for simple 
polytopes but for all convex polytopes which have the following 
``edge-reflecting''  property: If two vertices, $p$ and $q$, of $\Gamma$ 
are joined by an edge, $e$, then for every edge, $e'$, containing $p$, 
there exists a unique edge, $e''$, containing $q$, such that, $e'$ and $e''$ 
are coplanar. (This is a joint result with Ethan Bolker.)

This article consists of three chapters. In chapter one we review the 
theory of GKM manifolds and describe how to translate geometric properties 
of these manifolds into combinatorial properties of their associated 
GKM graphs. 
In chapter two we define an \emph{abstract one-skeleton} to be a labeled 
graph $(\Gamma, \alpha)$ for which $\alpha$ satisfies certain simple axioms 
(axiomatizing properties of the GKM-skeleta discussed in chapter one.) We 
then define the \emph{cohomology ring}, $H(\Gamma, \alpha)$, to be, as 
above, the set of all maps, $f: V_{\Gamma} \to \SS(\fg^*)$ which satisfy 
the compatibility conditions \eqref{eq:n1.2} and prove that this ring is 
a free $\SS(\fg^*)$-module with $b_{2i}(\Gamma)$ generators in dimension 
$2i$. (Involved in the proof of this theorem are the graph-theoretical 
analogues of two basic theorems in equivariant symplectic geometry: the 
Kirwan surjectivity theorem and the ``blow-up-blow-down'' theorem of 
Brion-Procesi-Guillemin-Sternberg-Godinho. Both these theorems have to 
do with the concept of symplectic reduction, and a large part of chapter 
two will be concerned with defining this concept in the context of 
abstract one-skeleta.) 

Chapter three contains a number of applications. One of these is the 
theorem about edge-reflecting polytopes which we described above. Another 
is a ``realization'' theorem for abstract GKM-skeleta. This asserts that an 
abstract one-skeleton $(\Gamma, \alpha)$ is the GKM one-skeleton of a GKM 
manifold if and only is $\alpha$ satisfies certain integrality conditions. 
This is a joint result with Viktor Ginzburg, Yael Karshon and Sue Tolman 
and is closely related to the realization theorem proved by them in 
\cite{GKT}. 

A third application has to do with the theory of Schubert polynomials. 
In Section \ref{ssec:schubs} we show that, for the Grassmannian, 
$Gr^k(\CC^n)$, 
the canonical generators of $H(\Gamma, \alpha)$ predicted by our theory 
have an alternative description in terms of the Hecke algebra of divided 
difference operators and thus can be identified with the ``double 
Schubert polynomials'' of \cite{Bi}. (Together with Tara Holm we have 
generalized this to all partial flag varieties; for details, see \cite{GHZ}.)

We would like to express our thanks to several of our colleagues
for helping us to understand some of the key motivating examples
in this subject:  David Vogan for furnishing us with an
enlightening example of a GKM action of $T^2$ on $S^6$, Yael 
Karshon, Viktor Ginzburg and Sue Tolman for furnishing us with an equally
enlightening example of a GKM action of $T^2$ on the
$n$-fold equivariant ramified cover of $S^2 \times S^2$, Mark Goresky for
pointing out to us the connection between the GKM theory of toric 
varieties and Stanley-Reisner theory, Werner Ballmann for making 
us aware of the fact that, for the Grassmannian, GKM theory
reduces to studying an object which graph theorists call the
\emph{Johnson graph}, Sara Billey for helping us to understand
the tie-in between our theory of Thom classes for this graph and
the standard Schubert calculus, Rebecca Goldin and  Tara Holm for 
helping us work out the details of this example (in sections 
\ref{ssec:grass} and \ref{ssec:schubs})
and Ethan Bolker for his beautiful observation (see section 
\ref{ssec:combbetti}, Theorem \ref{th:2.5})
that the Betti numbers of a GKM one-skeleton are well-defined,
independent of the choice of an admissible orientation. 
Last, but not least, we would like to thank the referee of this paper 
for a superb refereeing job.

\section{\bf GKM manifolds}
\label{sec:gkmmanif}

\subsection{\bf The GKM one-skeleton}
\label{ssec:gkmgraph}

For each of the weights $\alpha_{i,p}$ on the list (\ref{eq:1.0}) 
let $H_i$ be the identity component of the kernel of the map
$$\exp{\xi} \in G \longrightarrow \exp{(\sqrt{-1}\alpha_{i,p}(\xi))}$$
and let $X_i$ be the connected component of $M^{H_i}$ containing $p$.

\begin{theorem} \label{th:xe}
$X_i$ is diffeomorphic to $S^2$ and the action of $G$ on $X_i$ is 
diffeomorphic to the standard rotation action of the circle $G/H_i$ on
$S^2$.
\end{theorem}

\begin{proof}
Consider the decomposition
$$T_pM = \bigoplus T_p^{\alpha_{i,p}}$$
of $T_pM$ into 2-dimensional weight spaces. Our assumption that the weights 
(\ref{eq:1.0}) are pairwise linear independent imply that
$$T_pX_i = T_p^{\alpha_{i,p}}$$
and hence that $X_i$ is two-dimensional. Since the only oriented two-manifold
with faithful $S^1$ actions are $S^2$ and $T^2$, $X_i$ has to be one of them.
However, the action of $S^1$ on $T^2$ is fixed point free, so $X_i$ is 
diffeomorphic to $S^2$. Finally, the fact that the action of $S^1$ is the 
standard $S^1$ action is standard.
\end{proof}

Thus $p$ is the point of intersection of $d$ embedded $G$-invariant 2-spheres.
These can be represented graphically as in the figure below:
\begin{figure}[h]
\begin{center}
\includegraphics{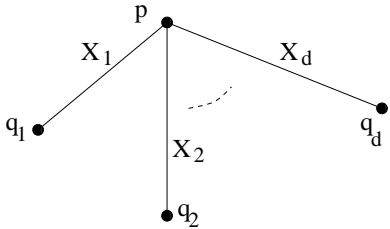}
\caption{A vertex of the GKM graph}
\end{center}
\end{figure}

\noindent Each of these 2-spheres joins $p$ to another fixed point $q_i$ 
and each $q_i$ is in turn the point of intersection of $d$ 2-spheres. 
One of these is $X_i$ and rejoins $q_i$ to $p$, but the others join $q_i$
to other fixed points, and at these points we can repeat the construction.
We will define {\em the GKM graph $\Gamma$} to be the graph we obtain by 
repeating this construction until we run out of fixed points. 

This graph can be defined more intrinsically as follows: the vertices of 
$\Gamma$ correspond to the fixed points of $M$,
an edge, $e$, of $\Gamma$ corresponds to a $G$-invariant embedded 
two-sphere, $X_e$, and joins the vertices that correspond to the two 
fixed points situated on $X_e$.
For an \emph{oriented} edge $e$, we will denote by $i(e)$ and $t(e)$ 
the initial and terminal vertices of $e$. In addition, we will denote 
by $\bar{e}$ the edge $e$ with its orientation reversed. 
Thus $i(\bar{e})=t(e)$ and $i(e)=t(\bar{e})$.

To keep track of the action of $G$ on this configuration of embedded 
$S^2$'s we will assign to each oriented edge $e$ the weight $\alpha_e$ of
the isotropy representation of $G$ on $T_{i(e)}X_e$. 
Denoting by $E_{\Gamma}$ the set of oriented edges of $E$, 
this gives us a map
$$\alpha : E_{\Gamma} \to \fg^*$$
which we will call {\em the axial function} 
of $\Gamma$. The pair $(\Gamma, \alpha)$ will be called the 
\emph{GKM one-skeleton} associated to the GKM manifold $M$.

Let $V_{\Gamma}$ be the set of vertices of $\Gamma$ and let
$$\pi : E_{\Gamma} \to V_{\Gamma}$$
be the fibration defined by $\pi(e) = i(e)$. A \emph{connection} on the 
bundle $(E_{\Gamma}, V_{\Gamma},\pi)$ is, by definition, a recipe for
transporting the fibers of $\pi$ along paths in $\Gamma$. In particular, 
a \emph{canonical connection} 
can be defined as follows. Let $e$ be an oriented edge of $\Gamma$ joining
the vertex $p=i(e)$ to the vertex $p'=t(e)$ and let $e_i$ and
$e_i'$, for $ i=1,..,d$, be the ``points'' (\emph{i.e.} oriented edges) 
on the fibers above $p$ and $p'$. By a  theorem of Klyashko (\cite{Kl}), 
the restriction to $X_e$ of the tangent bundle to 
$M$ splits equivariantly into a sum of line bundles
\begin{equation}
\label{eq:sumofline}
\bigoplus \LL_i, \qquad i=1,..,d
\end{equation}
and one can relabel the $e_i$'s and $e_i'$'s so that
$$(\LL_i)_p = T_pX_{e_i} \quad \mbox{ and } \quad 
  (\LL_i)_{p'} = T_{p'}X_{e_i'}$$
and from this one gets a canonical identification 
$e_i \longleftrightarrow e_i',$ 
\emph{i.e.} a canonical map
$$\theta_e : E_p \to E_{p'},$$
$E_p$ being the fiber $\pi^{-1}(p)$ above $p$ and $E_{p'}$ the fiber above
$p'$. 

Associated to the notion of connection is that of holonomy.
Consider a connection, $\theta$, on $\Gamma$ and fix 
$p \in V_{\Gamma}$. For each loop, $\gamma$, starting and 
ending at $p$ one gets a bijection, 
$\sigma_\gamma : E_p \to E_p$, by composing the maps corresponding
to the edges of $\gamma$. Let $Hol(\Gamma,\theta,p)$ be the
subgroup of the permutation group $\Sigma(E_p)$ generated by 
the elements of the form $\sigma_\gamma$ for all loops $\gamma$
based at $p$. If $p_1$ and $p_2$ can be connected by a path 
then the holonomy groups $Hol(\Gamma,\theta,p_1)$ and 
$Hol(\Gamma,\theta,p_2)$ are isomorphic by conjugacy; therefore 
we can define \emph{the holonomy group} 
$Hol(\Gamma,\theta)$ as being the group, $Hol(\Gamma,\theta,p)$,
for any point $p$. We will also say that $\theta$ has 
\emph{trivial holonomy}  
if $Hol(\Gamma_i,\theta)$ is trivial for each connected
component $\Gamma_i$ of $\Gamma$.

The following theorem lists some basic properties of the triple
$(\Gamma, \alpha, \theta)$. 

\begin{theorem}
\label{th:aximan}
\begin{enumerate}
\item \label{item1} 
	For every $p \in V_{\Gamma}$, the weights $\alpha_e$, $e \in E_p$, 
	are pairwise linearly independent.
\item For every $e \in E_{\Gamma}$, $(\theta_{e})^{-1} = \theta_{\bar{e}}$.
\item $\theta_e$ maps $e$ to $\bar{e}$.
\item \label{item4} $\alpha_{\bar{e}} = - \alpha_{e}$.
\item \label{item5}
	Let $p=i(e)$, $p'=t(e)$ and let $e_i \leftrightarrow e_i'$,
	$i=1,..,d$ be the map of $E_p$ onto $E_{p'}$ defining $\theta_e$.
	Then
\begin{equation}
\label{eq:axiomfive}
\alpha_{e_i'}= \alpha_{e_i} + c\alpha_e
\end{equation}
	for some constant $c=c_{i,e}$ depending on $i$ and $e$.
\end{enumerate}
\end{theorem}

\begin{proof}
The first four assertions are obvious. To prove the last one, let 
$H$ be the identity component of the kernel of the map
$$\exp{\xi} \in G \longrightarrow \exp{(\sqrt{-1}\alpha_e(\xi))}.$$
Each point of $X_e$ is an $H$-fixed point, so for each $x\in X_e$ one has 
an isotropy representation of $H$ on $T_xM$. The weights of this 
representation are independent of $x$, so they are the same at $p$ and $p'$.
\end{proof}

\begin{remark*}
It is easy to see that the constants $c_{i,e}$ in \eqref{eq:axiomfive}
are {\em integers}. In fact, let $c(\LL_i)$ be the Chern class of the 
line bundle $\LL_i$ in \eqref{eq:sumofline}. By the Atiyah-Bott-Berline-Vergne
localization theorem, the integral of $c(\LL_i)$ over $X_e$ is
\begin{equation}
\label{eq:cie}
c_{i,e}= \frac{\alpha_{e_i'} - \alpha_{e_i}}{\alpha_e} \quad ;
\end{equation}
hence $c_{i,e}$ is the Chern number of $\LL_i$, and so, in particular, 
an integer.
\end{remark*}

\subsection{\bf GKM theory for orbifolds}
\label{ssec:gkmorbif}

By \emph{``orbifolds''} 
we will mean orbifolds having a presentation of the 
form $M=X/K$, $K$ being a torus and $X$ being a manifold on which $K$ acts
in a faithful, locally free fashion. GKM theory for such orbifolds is
essentially the same as GKM theory for manifolds, the major difference
being that the $S^2$'s corresponding to the edges of the graph $\Gamma$ 
may now be orbifold $S^2$'s, that is are either \emph{tear-drops} or 
\emph{footballs}.
\begin{figure}[h]
\begin{center}
\includegraphics{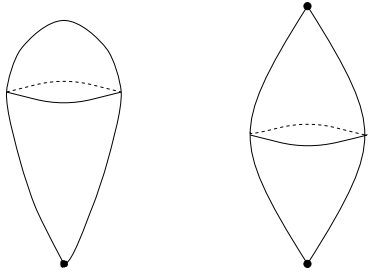}
\caption{A tear-drop and a football}
\label{fig:tear}
\end{center}
\end{figure}

One consequence of this is that the axiomatic properties of the axial 
function $\alpha$ are slightly more complicated. 
For each point $p=xK \in M$, with $x \in X$, let $m_p=\#K_x$. Then,
if $e \in E_{\Gamma}$, item \ref{item4} in Theorem \ref{th:aximan} has 
to be replaced by
\begin{equation}
\label{eq:131}
m_{i(\bar{e})} \alpha_{\bar{e}} = - m_{i(e)}\alpha_{e};
\end{equation}
however, the other properties of $\alpha$ described in Theorem \ref{th:aximan} 
are still true as stated.

\begin{example}
Let $M$ be the football of type $(m,n)$, that is the quotient:
$$(\CC^2-0)/\sim,$$
where
$$(z_1,z_2) \sim (\lambda^mz_1,\lambda^nz_2)$$
with $\lambda \in \CC -0$ and $m,n$ relatively prime positive integers.
Let $S^1$ act on $M$ by the action
$$e^{i\theta}[z_1:z_2] = [e^{i\theta}z_1: e^{i\theta}z_2].$$
The fixed points of this action are $p=[0:1]$ and $q=[1:0]$, and a 
coordinate system centered at $p$ is given by
$$z \in \CC \to [z:1].$$
In this coordinate system $z \sim \omega_n^k z$, $\omega_n$ being a 
primitive $n$-th root of unity, so, in particular, $m_p=n$. The action
of $S^1$ in this coordinate system is given by
$$e^{i\theta}[z:1] = [e^{i\theta}z:e^{i\theta}]=
[e^{i\frac{n-m}{n}\theta}z:1],$$
so the weight of the isotropy representation of $S^1$ on $T_pM$
is $\alpha_p= \frac{n-m}{n}$;
similarly $m_q=m$ and $\alpha_q=\frac{m-n}{m}$ so that
$$m_p \alpha_p = -m_q \alpha_q,$$
in confirmation of \eqref{eq:131}.

Notice, by the way, that the character associated with this 
weight, $e^{i\frac{n-m}{n}\theta}$, is taking values not in $S^1$ but in 
$S^1/\{\omega_n^k\}$. This is, of course, consistent with the fact that the 
linear action of $S^1$ on the coordinate system above is only an action 
modulo the identification $z \sim \omega_n^kz$.

\end{example}

\subsection{\bf Combinatorial Betti numbers}
\label{ssec:combbetti}

Let $M$ be a GKM manifold and $(\Gamma, \alpha)$ its GKM one-skeleton; 
we say that $\xi \in \fg$ is a \emph{polarizing vector} if 
$\alpha_e(\xi) \neq 0$ for all $e \in E_{\Gamma}$, and we denote by 
$\P$ \emph{the set of polarizing vectors}, \emph{i.e.}
\begin{equation}
\label{eq:polset}
\P = \{ \xi \in \fg ; \alpha_e(\xi) \neq 0 \mbox{ for all } 
e \in E_{\Gamma} \}.
\end{equation}
For a fixed $\xi \in \P$ define \emph{the index $\sigma_p$ of a vertex} 
$p \in V_{\Gamma}$ to be the number of edges $e \in E_p$ with 
$\alpha_e(\xi) < 0$. This definition clearly depends on the choice of $\xi$.
Let
\begin{equation}
\label{eq:bettino}
b_{2i}(\Gamma) = \# \{ p ; \sigma_p= i \}.
\end{equation}
We claim that this definition {\em doesn't} depend on $\xi$, in spite of 
the fact that $\sigma_p$ does, and we will call $b_{2i}(\Gamma)$ the 
\emph{(combinatorial) $2i$-th Betti number of} $\Gamma$.

\begin{theorem}
\label{th:2.5}
$b_{2i}(\Gamma)$ doesn't depend on $\xi$; it is a combinatorial 
invariant of $(\Gamma, \alpha)$.
\end{theorem}

\begin{proof}
Let $\P_i, \; i=1,...,N$, be the connected components of $\P$
and consider an $(n-1)$-dimensional wall separating two adjacent 
$\P_i$'s. This wall is defined by an equation of the form
\begin{equation}
\label{eq:wallbetti}
\alpha_{e}(\xi) =0
\end{equation}
for some $e \in E_{\Gamma}$. Let $p=i(e)$, $q=t(e)$,
and let's compute the changes in $\sigma_p$ and $\sigma_q$ as $\xi$
passes through this wall: 
Let $E_p= \{ e_i, \; i=1,...,d \}$ and 
$E_q= \{ e'_i, \; i=1,...,d \}$ 
(with $e_d=e$ and $e'_d=\bar{e}$).  
By item \ref{item5} of Theorem \ref{th:aximan} 
we can order the $e_i$'s 
so that, for $i \leq d-1$,
$$\alpha_{e_i} = \alpha_{e'_i} +c_i \alpha_{e} \; .$$
From item \ref{item1} of Theorem \ref{th:aximan} it 
follows that for every $i=1,...,d-1$,
$$\mbox{dim}\, (\,\mbox{ker}\, \alpha_{e} \cap \mbox{ ker}\, 
\alpha_{e_i}) =n-2$$
and therefore there exists $\xi_0$ such that 
$\alpha_{e}(\xi_0)=0$ but 
$\alpha_{e_i}(\xi_0) \neq 0 \neq  \alpha_{e'_i}(\xi_0)$, 
for all $i=1,...,d-1$.

Thus there exists a neighborhood $U$ of $\xi_0$ in $\fg$ such that
for $i=1,...,d-1$ and $\xi \in U$,  
$\alpha_{e_i}(\xi)$ and  $ \alpha_{e'_i}(\xi)$
have the same sign,
and this common sign doesn't depend on $\xi \in U$.
Such a neighborhood will intersect both regions created by the 
wall (\ref{eq:wallbetti}). Now suppose that $\xi \in U$ and that
$r$ of the numbers $\alpha_{e_i}(\xi), \; i=1,...,d-1,$ are 
negative. Since 
$\alpha_{\bar{e}}(\xi) = - \alpha_{e}(\xi)$, it
follows that for $\alpha_e(\xi)$ positive
$$\sigma_p =r \quad \mbox{ and } \quad \sigma_q =r+1$$
and for $\alpha_e(\xi)$ negative 
$$\sigma_p =r+1 \quad \mbox{ and } \quad \sigma_q =r.$$
In either case, as $\xi$ passes through the wall (\ref{eq:wallbetti}),
the Betti numbers don't change. 
\end{proof}

\begin{remark*}
When we change $\xi$ to $-\xi$, a vertex with index $k$ will now have 
index $d-k$, where $d$ is the valence of $\Gamma$. Since the Betti
numbers don't depend on $\xi$ it follows that
\begin{equation}
\label{eq:1.4}
b_{2(d-k)}= b_{2k}, \qquad \forall \; k=0,..,d.
\end{equation}
\end{remark*}

We will show in section \ref{ssec:examples} 
that these combinatorial Betti numbers
{\em may not be}, in general, equal to the Betti numbers of $M$. An important 
exception, however, is the following:

\begin{theorem}
\label{th:combbetti=topbetti}
If the action, $\tau$, of $G$ on $M$ is Hamiltonian then 
$b_{2i}(\Gamma) = b_{2i}(M).$
\end{theorem}

\begin{proof}
For $\xi \in \P$, the vector field $\xi_M$ is Hamiltonian and its 
Hamiltonian function, $f$, is a Morse function whose critical points are 
the fixed points of $\tau$. Moreover, the index of a critical point, $p$, 
is just $2\sigma_p$, so the number of critical points of index $2i$ is 
$b_{2i}(\Gamma)$; for more details see \cite{At}. 
\end{proof}

\begin{remark*}
In the Hamiltonian case, the relation \eqref{eq:1.4} is a direct consequence 
of Theorem \eqref{th:combbetti=topbetti} and Poincare duality.
\end{remark*}

\subsection{\bf Hamiltonian GKM-skeleta}
\label{ssec:hamiltgkm}

We will discuss in this section a few graph theoretic pathologies which
can't occur if $(\Gamma, \alpha)$ is the GKM one-skeleton 
of a Hamiltonian $G$-manifold.
As in the previous section, fix a vector $\xi \in \P$. 
For an unoriented edge $e$, of $\Gamma$, let 
$e^{+} \in E_{\Gamma}$ 
be the edge $e$, \emph{oriented} such that $\alpha_{e^+}(\xi) >0$ and 
let $e^{-} \in E_{\Gamma}$ be the edge $e$ with the opposite orientation.
We are thus defining an orientation, $o_{\xi}$ (which we will call 
the \emph{$\xi$-orientation}), of the edges of $\Gamma$, and it is clear that
this orientation depends only on the connected component  of $\P$ in 
which $\xi$ sits. On the other hand it 
is clear that different components will give rise to different orientations. 
(For instance, replacing $\xi$ with $-\xi$ reverses all the orientations.)
For an oriented edge $e$, let $e_0$ be the same edge, but unoriented;
we will say that $e$ \emph{points upward} if $e=e_0^+$ and that $e$ 
\emph{points downward} if $\bar{e}=e_0^+$.
We will say that $(\Gamma, \alpha)$ is \emph{$\xi$-acyclic} 
if the oriented 
graph $(\Gamma,o_{\xi})$ has no cycles and that $(\Gamma, \alpha)$ 
\emph{satisfies the no-cycle condition} if it is $\xi$-acyclic for at least 
one $\xi$.

\begin{definition}
Given $\xi \in \P$, a function $f: V_{\Gamma} \to \RR$ is called 
$\xi$-\emph{compatible} if, for every edge, $e$, of $\Gamma$,
$f(t(e^+)) >  f(i(e^+))$.
\end{definition}

If $f$ is injective and $\xi$-compatible, the
orientation, $o_{\xi}$, of $\Gamma$ associated with $\xi$ can't have cycles 
since $f$ has to be strictly increasing along any
oriented path. The converse is also true:

\begin{theorem} 
\label{th:posnocycle}
  If $(\Gamma, \alpha)$  is a $\xi$-acyclic GKM one-skeleton then 
   there exists an injective function, $f: V_{\Gamma} \to \RR$, 
   which is $\xi$-compatible.
\end{theorem}

\begin{proof} For a vertex $p$, define $f_0(p)$ to be the length 
of the longest oriented path with terminal vertex at $p$. Then $f_0$ is 
$\xi$-compatible and takes only integer values; 
a small perturbation of $f_0$ produces an injective function, $f$,  
which is still $\xi$-compatible.
\end{proof}

An important example of an acyclic GKM one-skeleton is the following:

\begin{theorem}
\label{th:hamnocycle}
If $(\Gamma, \alpha)$ is the GKM one-skeleton of a Hamiltonian 
$G$-manifold then it satisfies the no-cycle condition for all $\xi \in \P$.
\end{theorem}

\begin{proof}
Let $f$ be a Morse function as in the theorem of the previous section.
Then its restriction to $M^G$ is an $\xi$-compatible 
on the vertices, $V_{\Gamma}=M^G$, of $\Gamma$.
\end{proof}

We will next describe another type of pathology which can't occur if
$M$ is Hamiltonian. By Theorem \ref{th:combbetti=topbetti}, $M$ Hamiltonian
implies that the Betti numbers of $M$ are the same as the 
combinatorial Betti numbers of its GKM one-skeleton.
In particular, if $M$ is connected, $b_0(\Gamma)=1$. Thus, for each of the 
orientations above, there exists exactly one vertex, $p$, for which all the 
edges issuing from $p$ point upward.

This observation also applies to certain subgraphs of $\Gamma$. Recall from 
section \ref{ssec:gkmgraph} 
that $\Gamma$ is equipped with a connection, $\theta$. Let 
$\Gamma'$ be an $r$-valent subgraph of $\Gamma$, and for every vertex, $p$, 
of $\Gamma'$ let $E_p$ and $E_p'$ be the edges of $\Gamma$ and $\Gamma'$ 
issuing from $p$.

\begin{definition} \label{def:totgeod}
The subgraph $\Gamma'$ is called \emph{totally geodesic} with respect to 
$\theta$ if, for every oriented edge $e$ of $\Gamma'$ with $i(e)=p$ and 
$t(e)=q$, the holonomy map $\theta_e : E_p \to E_q$ maps $E_p'$ to $E_q'$.
\end{definition}

An important example of a totally geodesic subgraph of $\Gamma$ is the 
following. Let $\fh$ be a vector subspace of $\fg$ and let $\Gamma_{\fh}$
be the subgraph whose edges are the edges $e$ of $\Gamma$ for which 
$\alpha_{e^{\pm}} \in \fh^{\bot}$. Then all connected components of 
$\Gamma_{\fh}$ are totally geodesic.

It is easy to see that this graph is also a GKM graph:

\begin{theorem}
\label{th:gammah}
Let $H$ be the connected Lie subgroup of $G$ whose Lie algebra is $\fh$.
Then $\Gamma_{\fh}$ is the GKM graph of the $G$-space $M^H$.
\end{theorem}

If $M$ is Hamiltonian, so is $M^H$, so as a corollary of this theorem we 
obtain:

\begin{theorem}
\label{th:zerogammah}
If $M$ is Hamiltonian, every connected component of $\Gamma_{\fh}$
has zeroth Betti number equal to 1.
\end{theorem}

We will see in section \ref{ssec:examples} 
that the pathologies ruled out by Theorems
\ref{th:hamnocycle} and \ref{th:zerogammah} 
(the existence of cycles and the existence of totally
geodesic subgraphs with zeroth Betti number greater that one) can occur
``in nature'', that is, there are GKM manifolds for which both these
phenomena occur.

\subsection{\bf Reduction}
\label{ssec:reduction}

Let $M$ be a GKM manifold for which the action of $G$ on $M$ is Hamiltonian.
Let $H$ be a circle subgroup of $G$ with
$M^G = M^H$ and let $f : M \to \RR$ be its moment map. For every regular 
value $c$ of $f$, the reduced space
$$M_c = f^{-1}(c)/H$$
is a symplectic orbifold and the action on it of the quotient group $G_1=G/H$
is Hamiltonian (For Hamiltonian actions on orbifolds see \cite{LT}).

\begin{theorem}
\label{th:reduced}
The reduced space $M_c$ is a GKM orbifold for all regular values, $c$, of $f$ 
if and only if for every $p \in M^G$ the weights $\alpha_i=\alpha_{i,p}$ on 
the list (\ref{eq:1.0}) are \emph{three-independent}: 
that is, for every triple
of distinct values $i,j,k$, the weights, $\alpha_i, \alpha_j, \alpha_k$, 
are linearly independent.
\end{theorem}

\begin{proof}
We will prove the ``if'' part of this theorem by giving an explicit 
description of the one-skeleton of $M_c$. Let $(\Gamma, \alpha)$ be 
the GKM one-skeleton of $M$. Let $e$ be an oriented edge of $\Gamma$ with 
vertices, $p$ and $q$, for which
\begin{equation}
\label{eq:pcq}
f(p) < c < f(q)
\end{equation}
and let $X_e$ be the embedded two-sphere corresponding to $e$.  Then the 
reduction of $X_e$ with respect to $H$ at $c$ consists of a single 
$G_1$-fixed point, $p_e^c$, of $M_c$, and every $G_1$-fixed point is of this 
type. Thus the vertices of the GKM graph of $M_c$ are in one to one 
correspondence with the edges of $\Gamma$ satisfying (\ref{eq:pcq}).

What about the edges of this graph ? For the point, $p$, above, let's arrange 
the weight vectors on the list (\ref{eq:1.0}) such that 
$\alpha_d = \alpha_e$, 
and let $\alpha_i$ be any one of the remaining weight vectors. Let
$$H_i=\{ \exp{\xi} \in G ; \alpha_i(\xi)=\alpha_d(\xi)=0 \}$$
and let $W_i$ be the component of $M^{H_i}$ containing the point $p$.

\begin{lemma}
\label{lem:wi}
The dimension of $W_i$ is 4.
\end{lemma}

\begin{proof}
As in the proof of Theorem \ref{th:xe}, let 
$$T_pM = \bigoplus T_p^{\alpha_i}$$
be the decomposition of the tangent space to $M$ at $p$ into two-dimensional
weight spaces. The assumption that the $\alpha_i$'s are three-independent
implies that
$$T_pW_i = T_p^{\alpha_i} \oplus T_p^{\alpha_d}$$
and hence that $\dim{T_pW_i}=4$.
\end{proof}

Since $W_i$ is connected, its GKM graph, by Theorem \ref{th:hamnocycle}, 
consists of two
oriented chains, along each of which $f$ is strictly increasing. 
One of these chains contains the edge, $e$, and on the other chain there is
exactly one oriented edge, $e'$, whose vertices, $p'=i(e')$ and 
$q'=t(e')$, satisfy
$f(p') < c < f(q').$ 
Consider now the reduction of $W_i$ at $c$. This is a two-dimensional
symplectic orbifold on which the group $G_1$ acts faithfully and in a 
Hamiltonian fashion, so it has to be either a ``tear-drop'' or a 
``football'' (see section \ref{ssec:gkmorbif}). Moreover, it contains 
exactly two
$G_1$-fixed points, $p_c^e$ and $p_c^{e'}$. Thus, to summarize, we get the 
following description of the GKM graph, $\Gamma_c$, of $M_c$:
\begin{enumerate}
\item The vertices of this graph are the points $p_c^e$ corresponding to 
the oriented edges $e$ of $\Gamma$ which satisfy (\ref{eq:pcq}).
\item Let $p=i(e)$ and let $\alpha_{i,p}$ be a weight on the list 
(\ref{eq:1.0}) distinct from $\alpha_e$. Let $\fh_i$ be the 
codimension 2 subspace of $\fg$ defined by 
$\alpha_{i,p}(\xi) = \alpha_e(\xi) =0$ and let 
$\Gamma_i$ be the connected 
component of $\Gamma_{\fh_i}$ containing $p$. Then there are exactly
two oriented edges of $\Gamma_i$ satisfying (\ref{eq:pcq}): the edge $e$ and
another edge $e_i$.
\item The edges of $\Gamma_c$ meeting at $p_c^e$ are in one to one 
correspondence with the graphs $\Gamma_i$ and each of these edges joins 
the vertex $p_c^e$ to the vertex $p_c^{e_i}$.
\end{enumerate}

To complete this description of $\Gamma_c$ we must still describe the 
canonical connection on this graph and its axial function. This, however, we 
will postpone until later (see section \ref{ssec:redskel}).
\end{proof}

\subsection{\bf The flip-flop theorem}
\label{ssec:flipflop}

The flip-flop theorem describes how the orbifold, $M_c$, changes as $c$ goes 
through a critical value of $f$. Suppose there is exactly one critical point, 
$p \in M^G$, with $f(p)=c$. Let $W^+$ and $W^-$ be the unstable and stable
manifolds at $p$ of the gradient vector field associated with $f$. In a 
neighborhood of $p$, $W^+$ and $W^-$ can be identified with linear subspaces 
of the tangent space to $M$ at $p$. Namely, let $\xi$ be the element of
the Lie algebra of $H$ for which $\iota(\xi_M)\omega=df$, $\omega$ being the
symplectic form, and let 
$$T_pM = \bigoplus T_p^{\alpha_i}$$
be the decomposition of $T_pM$ into two-dimensional weight spaces. Then, in a 
neighborhood of $p$
\begin{equation}
\label{eq:tw+}
W^+ \simeq \bigoplus_{\alpha_i(\xi) > 0} T_p^{\alpha_i}
\quad \mbox{ and } \quad
W^- \simeq \bigoplus_{\alpha_i(\xi) < 0} T_p^{\alpha_i}.
\end{equation}
%
%and
%
%\begin{equation}
%\label{eq:tw-}
%W^- \simeq \bigoplus_{\alpha_i(\xi) < 0} T_p^{\alpha_i}.
%\end{equation}
%
In particular:

\begin{theorem}
\label{th:redweps}
For $\epsilon >0$ small enough, the reduced spaces, $W_{c+\epsilon}^+$ and
$W_{c-\epsilon}^-$, are the (twisted) projective spaces obtained by reducing
(\ref{eq:tw+}) 
%and (\ref{eq:tw-}) 
at $c+\epsilon$ and $c-\epsilon$ by the 
linear action of the circle group $H$.
\end{theorem}

The reduced spaces $W_{c+\epsilon}^+$ and $W_{c-\epsilon}^-$ are 
symplectic {\em sub-orbifolds} of $M_{c+\epsilon}$ and $M_{c-\epsilon}$ 
and the ``flip-flop'' theorem asserts:

\begin{theorem}
\label{th:flipflop}
The blow-up of  $M_{c+\epsilon}$ along $W_{c+\epsilon}^+$ is 
diffeomorphic as a $G_1$-manifold to the blow-up of 
$M_{c-\epsilon}$ along $W_{c-\epsilon}^-$.
\end{theorem}

\begin{remarks*}
\begin{enumerate}
\item The ``blowing-up'' referred to here is symplectic blow-up in the sense
of Gromov.
\item This theorem can be refined to describe how the symplectic structures 
of these two blow-ups are related; see \cite{GS2}.
\item This result is due to Guillemin-Sternberg and Godinho. 
An analogous result for 
complex manifolds (with G.I.T. reduction playing the role of symplectic
reduction) can be found in \cite{BP}.
\end{enumerate}
\end{remarks*}

To see how the GKM one-skeleton of $M_{c+\epsilon}$ is related to the GKM
skeleton of $M_{c-\epsilon}$ we must find out how GKM-skeleta are 
affected by blowing-up. Consider the simplest case of a 
\emph{blow-up}: Let $M$
be a GKM manifold, let $p$ be a point of $M^G$, let
$$T_pM = \bigoplus_{i=1}^d T_p^{\alpha_i}$$
be the decomposition of $T_pM$ into weight spaces and let $X_i$, $i=1,..,d$
be the embedded GKM 2-spheres at $p$ with 
$$T_pX_i = T_p^{\alpha_i}.$$
As a complex vector space $T_pM$ is $d$-dimensional, and each of these weight 
space is one-dimensional.

Now blow-up $M$ at $p$. As an abstract set, this blow-up is a disjoint union
of the projective space $\CC P(T_pM)$ and $M- \{ p \}$. The action of
$G$ on the first of these sets has exactly $d$ fixed points: a fixed point,
$p_i$, corresponding to each subspace $T_p^{\alpha_i}$ of $T_pM$; and 
each pair of fixed points $p_i$ and $p_j$ are joined by an embedded 
two-sphere, 
the projective line in $\CC P(T_pM)$ corresponding to the subspace 
$T_p^{\alpha_i} \oplus T_p^{\alpha_j}$ of $T_pM$. In addition, each of the 
two-spheres, $X_i$, is unaffected when we blow it up at $p$, but instead of
joining $p$ to a fixed point $q_i$ in $M-\{ p \}$, it now joins $p_i$ to
$q_i$. 

For the blow-up of $M$ along a $G$-invariant symplectic submanifold, 
$W^{2r}$, the story is essentially the same. As an abstract set the 
blow-up is the disjoint union of the projectivized normal bundle of $W$ and 
$M-W$. Thus each fixed point $p$ of $G$ in $W$ gets replaced, in the blow-up,
by $d-r$ new fixed points in the projectivized normal space to $W$ at $p$.
If $\Gamma$ is the GKM graph of $M$ and $\Gamma_1$ is the GKM graph of $W$,
then, just as above, the GKM graph of the blow-up is obtained from
$\Gamma$ and $\Gamma_1$ by replacing each vertex of $\Gamma_1$ by a complete
graph on $d-r$ vertices, one vertex for each edge of $\Gamma-\Gamma_1$ 
at $p$ (see Figure \ref{fig:blow} in section \ref{ssec:blowconstr}).

This description is particularly simple if $M$ is $M_{c-\epsilon}$
and $W$ is $W_{c-\epsilon}^-$. By Theorem \ref{th:redweps}, 
$W_{c-\epsilon}^-$ is just a 
twisted projective space of dimension $r-1$, $r$ being the index of the 
fixed point $p$, so its graph is the complete graph on $r$ vertices, 
$\Delta_r$. Hence, after blowing-up, it gets replaced by the graph 
$\Delta_r \times \Delta_{d-r}$. Similarly, the GKM graph of
$W_{c+\epsilon}^+$ is $\Delta_{d-r}$ and when we blow-up $M_{c+\epsilon}$
along $W_{c+\epsilon}^+$, it gets replaced by 
$\Delta_{d-r} \times \Delta_r$. Thus, as one passes through the critical 
value $c$, the following scenario takes place:
\begin{enumerate}
\item $\Delta_r$ gets blown-up to $\Delta_r \times \Delta_{d-r}$.
\item $\Delta_r \times \Delta_{d-r}$ gets ``flip-flopped'' to
$\Delta_{d-r} \times \Delta_r$.
\item $\Delta_{d-r} \times \Delta_r$ gets blown-down to $\Delta_{d-r}$.
\end{enumerate}

To complete the description of this transition we must still describe how 
this flip-flop process affects the connections and the axial functions
on these graphs. This, too, we will postpone until later, to section 
\ref{ssec:passage}.

\subsection{\bf Equivariant cohomology}
\label{ssec:thecohring}

Let $H_G(M)$ be \emph{the equivariant cohomology ring of} $M$ with complex 
coefficients. 
From the inclusion map $i : M^G \to M$ one gets a transpose map 
in cohomology
\begin{equation}
\label{eq:coh1}
i^* : H_G(M) \to H_G(M^G)
\end{equation}
and we will describe in this section some simple necessary conditions for 
an element of $H_G(M^G)$ to be in the image of this map. Since $M^G$ is
a finite set
\begin{equation}
\label{eq:coh2}
H_G(M^G) = \bigoplus H_G(\{p\}) ,  \quad p \in M^G;
\end{equation}
however, $H_G(\{p\})$ is the polynomial ring, $\SS(\fg^*)$, so the right side
of (\ref{eq:coh2}) is the ring
\begin{equation}
\label{eq:coh3}
\bigoplus^N \SS(\fg^*), \quad N=\#M^G.
\end{equation}
It is useful to keep track of the fact that each summand of (\ref{eq:coh3})
corresponds to a fixed point by identifying (\ref{eq:coh3}) with the ring
\begin{equation}
\label{eq:coh4}
\mbox{Maps}(V_{\Gamma}, \SS(\fg^*)).
\end{equation}

Let $e \in E_{\Gamma}$ and let $\fg_e^*$ be the quotient of $\fg^*$
by the one-dimensional subspace $\{ c\alpha_e; c \in \CC \}$. 
From the projection $\rho_e : \fg^* \to \fg_e^*$ one gets an epimorphism 
of rings
\begin{equation}
\label{eq:coh5}
\rho_e : \SS(\fg^*) \to \SS(\fg_e^*).
\end{equation}
(Since $\fg_e^* = \fg_{\bar{e}}^*$, we will use the notations 
$\fg_e^*$ and $\rho_e$ for unoriented edges, as well). 

\begin{theorem}
\label{th:imro}
A necessary condition for an element, $\phi$, of the ring (\ref{eq:coh4})
to be in the image of the map (\ref{eq:coh1}) is that for every edge, $e$, of 
$\Gamma$ it satisfies the compatibility condition
\begin{equation}
\label{eq:coh6}
\rho_e\phi_p = \rho_e\phi_q,
\end{equation}
$p$ and $q$ being the vertices of $e$ and $\phi_p$ and 
$\phi_q$ the elements of 
$\SS(\fg^*)$ assigned to them by the map $\phi : V_{\Gamma} \to \SS(\fg^*)$.
\end{theorem}

\begin{proof}
The right and left hand sides of \eqref{eq:coh6} are the pull-backs to 
$p$ and $q$ of the image of $\phi$ under the map $H_G(M) \to H_K(M)$, 
where $K=\exp{ker \; \alpha_e}$ and $\fk=\fg_e$. Since $p$ and $q$ 
belong to the same connected component of $M^K$, the pull-backs coincide.
\end{proof}

Let us denote by $H(\Gamma,\alpha)$ (or simply by $H(\Gamma)$ when the 
choice of $\alpha$ is clear) the subring of (\ref{eq:coh4}) 
consisting of those elements which satisfy the compatibility condition
(\ref{eq:coh6}); we will call $H(\Gamma, \alpha)$ the 
\emph{cohomology ring} of $(\Gamma, \alpha)$.
By the theorem above the map (\ref{eq:coh1}) factors 
through $H(\Gamma,\alpha)$ to give a ring homomorphism
\begin{equation}
\label{eq:coh8}
i^* : H_G(M) \to H(\Gamma,\alpha).
\end{equation}
This homomorphism also has a bit of additional structure. The constant maps
of $V_{\Gamma}$ into $\SS(\fg^*)$ obviously satisfy the condition
(\ref{eq:coh6}), so that $\SS(\fg^*)$ is a subring of $H(\Gamma,\alpha)$.
Also, from the constant map $M \to pt$ one gets a transpose map 
$H_G(pt) \to H_G(M)$, mapping $\SS(\fg^*)$ into $H_G(M)$ and it is easy 
to see that (\ref{eq:coh8}) is a morphism of $\SS(\fg^*)$-modules.
One of the main theorems of \cite{GKM} asserts that the homomorphism 
(\ref{eq:coh8}) is frequently an {\em isomorphism}. More explicitly, 
recall that if $K$ is a subgroup of $G$ there is a forgetfulness map
$H_G(M) \to H_K(M)$ and, in particular, for $K= \{e\}$, there is a map
\begin{equation}
\label{eq:coh9}
H_G(M) \to H(M).
\end{equation}

\begin{definition} 
$M$ is \emph{equivariantly formal} if (\ref{eq:coh9}) is onto.
\end{definition}

There are many alternative equivalent definitions of equivariant formality. 
For instance, for every compact $G$-manifold
\begin{equation}
\label{eq:coh10}
\sum \dim H^i(M) \geq \sum \dim H^i(M^G),
\end{equation}
and $M$ is equivariantly formal if and only if the inequality is equality.
Thus for GKM manifolds, $M$ is equivariantly formal if the sum of its Betti 
numbers is equal to the cardinality of $M^G$. In other words:

\begin{theorem}
\label{th:eqfor}
For a GKM manifold equivariant formality is equivalent to
\begin{equation}
\label{eq:coh11}
\sum b_i(M) = \sum b_i(\Gamma).
\end{equation}
\end{theorem}

For instance, if the action of $G$ on $M$ is Hamiltonian then 
$b_i(M) =  b_i(\Gamma)$ so (\ref{eq:coh11}) is trivially satisfied. A less
trivial example of (\ref{eq:coh11}) is the action of the Cartan 
subgroup of $G_2$ on the 6-sphere $G_2/SU(3)$. For this example we will
show in section \ref{ssec:examples} that $b_2(\Gamma)=b_4(\Gamma)=1$; 
$b_0(M)=b_6(M)=1$ and all the other Betti numbers are zero. 
The theorem of Goresky-Kottwitz-MacPherson which we alluded to above 
asserts:

\begin{theorem}
\label{th:gkm}
If $M$ is equivariantly formal, the map (\ref{eq:coh8}) is a bijection.
\end{theorem}

In other words, if $M$ is equivariantly formal, the equivariant cohomology
ring of $M$ is isomorphic to the cohomology ring $H(\Gamma,\alpha)$
of the GKM one-skeleton, $(\Gamma, \alpha)$. 
Recently, a number of relatively simple proofs have 
been given of this theorem: for example, a proof of Berline-Vergne \cite{BV}
based on localization ideas, and, in the Hamiltonian case, a very simple 
Morse theoretic proof by Tolman-Weitsman \cite{TW2}. 
Theorem \ref{th:eqfor} is, as we mentioned above, just one of many 
alternative criteria for equivariant formality. Another is:

\begin{theorem}
\label{th:eqfor2}
$M$ is equivariantly formal if, as $\SS(\fg^*)$-modules,
$$H_G(M) \simeq H(M) \otimes \SS(\fg^*).$$
\end{theorem}

Thus if $(\Gamma, \alpha)$ is the GKM one-skeleton of $M$, one gets from this 
the following result:

\begin{theorem}
\label{th:free}
If $M$ is equivariantly formal, $H(\Gamma,\alpha)$ is a free module
over $\SS(\fg^*)$ with $b_{2i}(M)$ generators in dimension $2i$.
\end{theorem}

One of the questions which we will address in the second part of this paper
is: ``When is the graph theoretical analogue of this theorem true with
the $b_{2i}(M)$'s replaced by the $b_{2i}(\Gamma)$'s ?'' From the
examples in section \ref{ssec:examples} we will see that even for 
GKM-skeleta 
this theorem is not true with the $b_{2i}(M)$'s replaced by the 
$b_{2i}(\Gamma)$'s. However, we will show that one can make this 
substitution providing $\Gamma$ has the properties described in Theorems 
\ref{th:zerogammah} and \ref{th:reduced}.

\subsection{\bf The Kirwan map}
\label{ssec:kirwmap}

Let $M$ be a compact Hamiltonian $G$-manifold, $H$ a circle subgroup of $G$ 
and $f : M \to \RR$ the $H$-moment mapping. If $c$ is a regular value 
of $f$ then the reduced space 
$$M_c = f^{-1}(c)/H$$
is a Hamiltonian $G_1$-space, with $G_1=G/H$, and one can define a 
morphism in cohomology
\begin{equation}
\label{eq:kir1}
\K_c : H_G(M) \to H_{G_1}(M_c)
\end{equation}
as follows. Let $Z=f^{-1}(c)$. Since $c$ is a regular value of $f$, 
the action of $H$ on $Z$ is locally free, so there is an isomorphism in 
cohomology (cf. \cite[sec. 4.6]{GS3}):
$$H_G(Z) \to H_{G_1}(M_c)$$
and the map (\ref{eq:kir1}) is just the composition of this with the 
restriction map
$$H_G(M) \to H_G(Z).$$
The homomorphism (\ref{eq:kir1}) is called \emph{the Kirwan map}, and a 
fundamental result of Kirwan (see \cite{Ki}) is:

\begin{theorem}
\label{th:kirw}
The map (\ref{eq:kir1}) is surjective.
\end{theorem}

One way of proving this theorem is to use the flip-flop theorem 
of section \ref{ssec:flipflop}: 
Let $c_1$ and $c_2$ be regular values of $f$ and 
suppose that there is just one critical point, $p$, of $f$ with
$c_1 < f(p) < c_2$. Assume by induction that Kirwan's theorem is true for 
$c_1$ and prove it for $c_2$. The flip-flop theorem says that $M_{c_2}$
is obtained from $M_{c_1}$ by a blow-up followed by a blow-down, 
and to see what effect these operations have on cohomology one makes 
use of the following theorem \cite{McD} :

\begin{theorem}
\label{th:dm}
Let $M$ be a compact Hamiltonian $G$-manifold and $W$ a $G$-invariant
symplectic submanifold of $M$. If $\beta : M^{\#} \to M$ is the symplectic
blow-up of $M$ along $W$ and $W^{\#} = \beta^{-1}(W)$ is its singular locus,
then there is a short exact sequence in cohomology
\begin{equation}
\label{eq:shexseq}
0 \to H_G(M) \to H_G(M^{\#}) \to H_G^{\natural}(W^{\#}) \to 0,
\end{equation}
the first arrow being $\beta^*$, the second being restriction and 
$H_G^{\natural}(W^{\#})$ being the quotient, $H(W^{\#})/\beta^*H(W)$.
\end{theorem}

Suppose now that $M$ satisfies the hypotheses of Theorem \ref{th:reduced}. 
Then both $M$ and $M_c$ are GKM spaces. Let $(\Gamma, \alpha)$ and 
$(\Gamma_c,\alpha_c)$ be their GKM-skeleta. By Theorem
\ref{th:gkm}, $H_G(M)$ is isomorphic to $H(\Gamma,\alpha)$ and 
$H_{G_1}(M_c)$ is isomorphic to $H(\Gamma_c,\alpha_c)$, so, from 
(\ref{eq:kir1}) we get a Kirwan map 
$$\K_c : H(\Gamma,\alpha) \to H(\Gamma_c,\alpha_c).$$
We will show that there is a purely graph theoretical description of this 
map: Recall that an element $\phi$ of $H(\Gamma,\alpha)$ is a map
of $V_{\Gamma}$ to $\SS(\fg^*)$ which, for every edge $e \in E_{\Gamma}$, 
satisfies the compatibility condition (\ref{eq:coh6}), 
$p=i(e)$ and $q=t(e)$ being the vertices of $e$. 
Suppose that $f(p)<c<f(q)$. Then $e$ corresponds to a 
vertex $p_c^e$ of $\Gamma_c$. Moreover, if $\fh$ is the Lie algebra of $H$
then $\fg_1 =\fg/\fh$ so that there is a map $\fg_1^* \to \fg^*$ which can 
be composed with the map $\rho_e : \fg^* \to \fg_e^*$ to give a bijection,
$\fg_1^* \to \fg_e^*$, and an inverse bijection, $\fg_e^* \to \fg_1^*$. 
This, in turn, induces an isomorphism of rings 
$$\gamma_e : \SS(\fg_e^*) \to \SS(\fg_1^*)$$
and hence we get an element $\gamma_e\rho_e\phi_p = \gamma_e\rho_e\phi_q$ 
of $\SS(\fg_1^*)$.

\begin{theorem}
\label{th:kirate}
The value of $\K_c\phi$ at the vertex $p_c^e$ of $\Gamma_c$ is 
$\gamma_e\rho_e\phi_p$.
\end{theorem}

\begin{proof}
Let $a$ be the element of $H_G(M)$ whose restriction to $M^G$ is $\phi$.
Let $X_e$ be the embedded two-sphere in $M$ corresponding to $e$ and let
$a_e\in H_G(X_e)$  be the restriction of $a$ to $X_e$. Then the one-point
manifold $\{ p_c^e\}$ is the reduction of $X_e$ at $c$ with respect to $H$.
Therefore it suffices to check that $\gamma_e\rho_e\phi_p$ is the image
of $a_e$ under the Kirwan map
$$H_G(X_e) \to H_{G_1}(\{p_c^e\}).$$
\end{proof}

\subsection{\bf Examples}
\label{ssec:examples}

We will describe in this section two examples of GKM manifolds for which the 
Betti numbers $b_{2i}(M)$ {\em don't} coincide with the combinatorial Betti
numbers $b_{2i}(\Gamma)$.

\begin{example} $\mathbf{G_2/SU(3)}$: 

\begin{figure}[h]
\begin{center}
\includegraphics{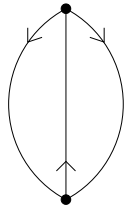}
\caption{$T^2$ action on $S^6$}
\label{fig:s6}
\end{center}
\end{figure}

This space is topologically just the standard
6-sphere. Moreover, if $p$ is the identity coset, the isotropy 
representation of $SU(3)$ on $T_p$ is the standard representation of 
$SU(3)$ on $\CC^3$ and hence the complex structure on $T_p$ given by the 
identification $T_p \simeq \CC^3$ extends to a $G_2$ invariant almost 
complex structure on $S^6$. Let $T^2$ be the Cartan subgroup of $G_2$. We will 
show that the action of $T^2$ on $S^6$ is a GKM action, determine the
GKM one-skeleton and compute its combinatorial Betti numbers. 

Recall that $G_2$ is by definition the group of automorphisms of the 
non-associative eight dimensional algebra of Cayley numbers and that there 
is an \emph{intrinsic} description of the almost complex structure on $S^6$
which makes use of algebraic proprieties of the Cayley numbers.
(For details see \cite[pp. 139-140]{KN}.) In this description an element 
of the Cayley numbers is identified with an element, 
$x=(z_1,z_2,w_1,w_2)$, of $\CC^4$ and $S^6$ is
realized as the unit sphere in the real subspace, 
$\overline{z_1} = -z_1$, that is
$$S^6 = \{ x \in \CC^4 \; ; \; 
|z_1|^2+|z_2|^2 + |w_1|^2 + | w_2|^2 =1 , \overline{z_1} = -z_1 
\}. $$

Let $\alpha$ and $\beta$ be basis vectors for the weight lattice of $T^2$. 
Then the action of $T^2$ on the Cayley algebra defined by 
$$e^{ i\theta} \cdot (z_1,z_2,w_1,w_2) \longmapsto 
(z_1, e^{ i(\alpha+\beta)(\theta)}z_2, e^{- i\alpha(\theta)}w_1,
e^{ i\beta(\theta)}w_2) $$
is an action by automorphisms (see \cite{Ja}), and it clearly leaves 
$S^6$ fixed, which implies that the induced action of $T^2$ on $S^6$ 
preserves the almost complex structure. Let $p=(i,0,0,0)$ and 
$q=(-i,0,0,0)$ be the fixed points of this induced action. If we identify
$T_pS^6$ with
$$\{(0,z_2,w_1,w_2) \; ; \; z_2,w_1,w_2 \in \CC \} \subset \CC^4$$
then the almost complex structure at $p$ is 
$$J_p(0,z_2,w_1,w_2) = (0,iz_2,iw_1,-iw_2).$$
Hence, identifying $T_pS^6$ with $\CC^3$ by
$$(0,z_2,w_1,w_2) \longmapsto (z_2,w_1,\overline{w_2}),$$
we deduce that the weights of the induced representation 
of $T^2$ on $T_pS^6$ are $\alpha+\beta, -\alpha$ and $-\beta$. Similarly,
for the representation of $T^2$ on $T_qS^6$, the weights are
$-\alpha -\beta, \alpha$ and  $\beta$. 

The GKM graph $\Gamma$ of this $T^2$-space will consist of the two vertices,
$p$ and $q$, linked by three edges (see Figure \ref{fig:s6}), 
labeled by the weights 
$\alpha+\beta, -\alpha$ and $-\beta$ and along every edge the 
connection swaps the remaining two edges. For every
$\xi \in \P$, the oriented graph $(\Gamma, o_{\xi})$ has cycles and its 
combinatorial Betti numbers are $b_0=b_3=0$, $b_1=b_2=1$.

\begin{remark*}
Note that, by Theorem \ref{th:eqfor}, $G_2/SU(3)$ is equivariantly formal, 
so, in spite of the fact that the Betti numbers don't coincide with the 
combinatorial Betti numbers, still $H(\Gamma,\alpha) = H_{T^2}(S^6)$.
\end{remark*}
\end{example}

\begin{example} 
{\bf The $n$-fold equivariant ramified cover of $S^2 \times S^2$:} 

In the previous example, every $\xi$-orientation of the GKM graph 
had cycles. We will next describe and example of a GKM manifold whose 
GKM graph 
does have an acyclic $\xi$-orientation but for which the combinatorial 
Betti numbers are different from the topological Betti numbers. 
The 4-manifold
$$W=S^2 \times S^2 = \CC P^1 \times \CC P^1$$
is a toric variety whose moment polytope $\Box_4$ is  the square in 
$\RR^2$ with vertices at (1,1), (-1,1), (-1,-1) and (1,-1).

Let $\phi: W \to \Box_4$ be the moment map and let 
$\psi : \RR^2 \to \RR^2$ be the map 
$$(x,y)=x+iy \to (x+iy)^n.$$
The pre-image of $\Box_4$ under this map is a regular curved 
polygon $\Box_{4n}$
with $4n$ sides. The fiber product of $W$ and $\Box_{4n}$
\begin{equation}
\label{eq:exm}
M= \{ (p,z) \in W \times \Box_{4n} ; \phi(p)=\psi(z) \} 
\end{equation}
is a connected compact manifold, with maps
$$\pi : M \to W, \quad \pi(p,z)=p \quad
\mbox{ and } \quad
 \gamma:M \to \Box_{4n}, \quad \gamma(p,z) = z.$$
Moreover, if we let $T^2$ act trivially on $\Box_{4n}$ and act on $W$ by
its given action, we get from (\ref{eq:exm}) an action of $T^2$ on $M$ 
which makes the ``fiber product'' diagram
\begin{equation}\label{eq:diagram}
\begin{CD}
M @>{\gamma}>> \Box_{4n} \\
@V{\pi}VV     @VV{\psi}V \\
W @>{\phi}>> \Box_4
\end{CD}
\end{equation}
$T^2$ equivariant. Let $M^{\#} = M - \gamma^{-1}(0)$ and 
$W^{\#} = W - \phi^{-1}(0)$.

\begin{lemma}
The map $\pi : M^{\#} \to W^{\#}$ is an $n$-to-1 covering map.
\end{lemma}

\begin{proof}
This follows from (\ref{eq:diagram}) and the fact that
$\psi : \Box_{4n}-\{ 0 \} \to \Box_4 - \{0 \}$ 
is an $n$-to-1 covering map.
\end{proof}

\begin{corollary}
There is a $T^2$-invariant complex structure on $M^{\#}$.
\end{corollary}

Let $\Box_{4n}^0$ be the interior of $\Box_{4n}$ and let
$M^0=\gamma^{-1}(\Box_{4n}^0)$. 

\begin{lemma}
There is a $T^2$-equivariant diffeomorphism 
$M^0 \to T^2 \times \Box_{4n}^0$ 
which intertwines $\gamma$ and the map
$pr_2 : T^2 \times \Box_{4n}^0 \to \Box_{4n}^0$.
\end{lemma}

\begin{proof}
Let $\Box_{4}^0$ be the interior of $\Box_4$ and let 
$W^0=\phi^{-1}(\Box_{4}^0)$. Since $W$ is a toric variety, there is a 
$T^2$-equivariant diffeomorphism 
$$W^0 \to T^2 \times \Box_4^0$$
which intertwines $\phi$ and the map 
$pr_2 : T^2 \times \Box_{4n}^0 \to \Box_{4n}^0$, so the lemma follows 
from the description (\ref{eq:exm}) of $M$.
\end{proof}

Let us use this diffeomorphism to pull back the complex structure on
$M^0 \cap M^{\#}$ to $T^2 \times (\Box_{4n}^0-\{0\})$. We will show that by
modifying this structure, if necessary, on a small neighborhood of
$T^2 \times \{0\}$, we can extend it to a $T^2$-invariant almost complex
structure on $T^2 \times \Box_{4n}^0$. Since the tangent bundle of
$T^2 \times \Box_{4n}^0$ is trivial, a $T^2$-invariant almost complex 
structure on $T^2 \times \Box_{4n}^0$ is simply a map
\begin{equation}
\label{eq:exj}
J : \Box_{4n}^0 \to GL(4,\RR)_+/ GL(2,\CC),
\end{equation}
so to prove this assertion we must show that the map
\begin{equation}
\label{eq:exj0}
J_0 : \Box_{4n}^0-\{0\} \to GL(4,\RR)_+/ GL(2,\CC)
\end{equation}
associated with the complex structure on 
$T^2 \times (\Box_{4n}^0- \{0 \})$ can be modified slightly on a small
disk $D$ about the origin so that it is extendible over $\Box_{4n}^0$.
However, 
$GL(4,\RR)_+/ GL(2,\CC)$ is homotopy equivalent to the 2-sphere
$SO(4)/U(2)$ so the restriction of $J_0$ to  $\partial D$ is a map 
$J_0 : S^1 \to S^2$ and since $S^2$ is simply connected, this map extends 
over the interior.

Pulling this almost complex structure back to $M^0$ we conclude:

\begin{theorem}
\label{th:acsm}
There exists a $T^2$-invariant almost complex structure on $M$.
\end{theorem}

We will now show that $M$ is a GKM manifold and compute its combinatorial 
Betti numbers. From the fact that $W$ is a toric variety one easily proves:

\begin{lemma}
\label{lem:1skw}
The one-skeleton of $W$ is $W-W^0$ and is the union of four two-spheres
with GKM graph $\Gamma_4 =\partial \Box_4$.
Moreover, its axial function is the function that assigns to the 
edges of $\Gamma_4$ the weights 
$\alpha_1=(1,0), \alpha_2=(0,1), \alpha_3=(-1,0)$ and $\alpha_4=(0,-1)$
(starting from the bottom edge and proceeding counter-clockwise).
\end{lemma}

Since $T^2$ acts freely on $M^0$ and since 
$\pi : M-M^0 \to W-W^0$ 
is a $n$-to-1 covering map, the lemma implies:

\begin{theorem}
\label{th:1skm}
The one-skeleton of $M$ is $M-M^0$ and it is the union of $4n$ two-spheres
with GKM graph $\Gamma_{4n} = \partial \Box_{4n}$. Its
axial function is the pull-back of the axial function of $W$ by the map
$\psi: \partial \Box_{4n} \to \partial \Box_4$.
\end{theorem}

In particular, $b_{2i}(\Gamma_{4n}) = nb_{2i}(\Gamma_4)$ so that
$b_0(\Gamma_{4n})=b_4(\Gamma_{4n})=n$ and $b_2(\Gamma_{4n})=2n$.

By a simple Mayer-Vietoris type computation, with $M, M_0$ and $M-M_0$, 
it is easy to compute the ``honest'' Betti numbers of $M$ and to show that 
$$b_0(M)=b_4(M)=1 \quad \mbox{ and } \quad b_2(M)=4n-2$$
and also that the odd Betti numbers are zero. Thus, in particular, by
Theorem \ref{th:eqfor}, $M$ is equivariantly formal and 
$H_G(M)= H(\Gamma,\alpha)$.

\end{example}

\subsection{\bf Edge-reflecting polytopes}
\label{ssec:polytopes}

Let $\Delta$ be an edge-reflecting polytope and let $\Gamma$ be its
one-skeleton: the graph consisting of the vertices and edges of $\Delta$.
The edge-reflecting property enables one to define a connection on $\Gamma$
as follows: Let $p$ and $p'$ be adjacent vertices of $\Delta$ and $e$ the
edge joining $p$ to $p'$. 
If $e_i$ is an edge joining $p$ to a vertex $q_i \neq p'$ 
then there exists, by the edge-reflecting property, a unique edge $e_i'$,
joining $p'$ to another vertex $q_i' \neq p$, such that $p,p',q_i$ and 
$q_i'$ are collinear. The correspondence $e_i \leftrightarrow e_i'$ and
$e \leftrightarrow \bar{e}$ defines a bijective map
$$\theta_e : E_p \to E_{p'}$$
and the collection of these maps is, by definition, a connection on $\Gamma$.

We can also define an axial function 
$$\alpha : E_{\Gamma} \to \RR^n$$
by attaching to each oriented edge $e$ the vector
$$\alpha_e = \overrightarrow{pq},$$
where $p=i(e)$ and $q=t(e)$ are the endpoints of $e$.

The triple $(\Gamma, \theta, \alpha)$ doesn't quite satisfy the properties 
described in Theorem \ref{th:aximan}. 
It does satisfy the first four of them but it only 
satisfies a somewhat weaker version of the fifth, namely:
\begin{equation}
\label{eq:axiom5}
\alpha_{e_i'} = \lambda_{i,e}\alpha_{e_i} + c_{i,e}\alpha_e \mbox{ with } 
\lambda_{i,e} > 0.
\end{equation}

\begin{proof}
Condition (\ref{eq:axiom5}) is just a restatement of the assumption that 
$e_i, e_i'$ and $e$ are coplanar; the positivity of $\lambda_{i,e}$ is a 
consequence of the convexity of $\Delta$. If $\lambda_{i,e}$ weren't positive,
$e$ would be in the {\em interior} of the intersection of $\Delta$ with the 
plane spanned by $e_i$ and $e_i'$.
\end{proof}

\begin{remarks*}
\begin{enumerate}
\item We will call $(\Gamma, \alpha)$ the \emph{GKM one-skeleton} of $\Delta$. 

\item For edge-reflecting polytopes, item \ref{item1} 
of Theorem \ref{th:aximan} can be 
replaced by the much stronger statement:

\begin{quote}
{\em For every $p \in V_{\Gamma}$, the vectors 
$\alpha_e \in E_p$ are $n$-independent: for every sequence 
$1\leq i_1 < i_2 < ... < i_n \leq d$, the vectors 
$\alpha_{i_1}, \alpha_{i_2}, ..., \alpha_{i_n}$ are linearly independent.}
\end{quote}

\item Moreover, if $(\Gamma, \alpha)$ is the GKM one-skeleton of an
edge-reflecting polytope, it satisfies both the ``no-cycle'' condition
of Theorem \ref{th:hamnocycle} 
and the ``zeroth Betti number'' condition of Theorem
\ref{th:zerogammah}. 
\end{enumerate}
\end{remarks*}

\subsection{\bf Grassmannians as GKM manifolds}
\label{ssec:grass}

Let $G$ be the $n$-torus $(S^1)^n$ and $\tau_0$ the representation of $G$
on $\CC^n$ given by 
$$\tau_0(e^{i\theta})z = ( e^{i\theta_1}z_1,...,e^{i\theta_n}z_n).$$
We will denote by $v_i$, $i=1,..,n$, the standard basis vectors of $\CC^n$
and by $\alpha_i$, $ i=1,..,n$, the weights of $\tau_0$ associated with 
these basis vectors. Thus, identifying $\fg$ with $\RR^n$, 
\begin{equation}
\label{eq:grass1}
\alpha_i(\xi) = \xi_i, \quad \mbox{ for } \xi = (\xi_1,..,\xi_n) \in \fg.
\end{equation}
From the action, $\tau_0$, we get an induced action, $\tau$, of $G$ on the
Grassmannian $Gr^k(\CC^n)$. We will prove that this is a GKM action by 
proving :

\begin{theorem}
\label{th:grass=gkm}
The fixed points of $\tau$ are in one-to-one correspondence with the 
$k$-element subsets of $\{ 1,..,n\}$. For the fixed point $p=p_S$, 
corresponding to the subset $S$, the isotropy representation of $G$
on $T_p$ has weights 
\begin{equation}
\label{eq:grass2}
\alpha_j - \alpha_i \quad, \qquad i \in S , \; j \in S^c.
\end{equation}
\end{theorem}

\begin{proof}
The fixed point $p_S$ corresponding to $S$ is the subspace $V_S$ of 
$\CC^n$ spanned by $\{ v_i \; ; \; i \in S \}$. Therefore the tangent 
space at $p_S$ is
\begin{equation}
\label{eq:grass3}
\mbox{Hom}_{\CC}(V_S, V_{S^c}),
\end{equation}
with basis vectors $\{v_j \otimes v_i^* \; ; \; i \in S \; , j \in S^c\}$
and the weights associated with these basis vectors are (\ref{eq:grass2}).
\end{proof}

We leave the following as an easy exercise:

\begin{theorem}
\label{th:3indep}
If $n \geq 4$ the weights (\ref{eq:grass2}) are 3-independent.
\end{theorem}

In particular $\tau$ is a GKM action as claimed. Let $\Gamma$ be its GKM 
graph. The following theorem gives a description of the $G$-invariant 
2-spheres which correspond to the edges of this graph:

\begin{theorem}
\label{th:2spheres}
Let $S$ and $S'$ be $k$-element subsets of $\{ 1,..,n\}$ with 
$\# (S\cap S') = k-1$. Define $S_1 = S\cap S'$ and $S_2=S \cup S'$ and let
$X_{S,S'}$ be the set of all $k$-dimensional subspaces $V$ of $\CC^n$ such
that 
\begin{equation}
\label{eq:grass4}
V_{S_1} \subset V \subset V_{S_2}.
\end{equation}
Then $X_{S,S'}$ is a $G$-invariant embedded 2-sphere containing $p_S$ and
$p_{S'}$ and all $G$-invariant embedded 2-spheres containing $p_S$ are of 
this form.
\end{theorem}

\begin{proof}
From the identification of $X_{S,S'}$ with the projective space
$\CC P(V_{S_2}/V_{S_1})$ one sees that $X_{S,S'}$ is an embedded 2-sphere.
Moreover, since $V_S$ and $V_{S'}$ satisfy (\ref{eq:grass4}), this
sphere contains $p_S$ and $p_{S'}$. To prove the last assertion note that
the tangent space to $X_{S,S'}$ at $p_S$ is
\begin{equation}
\label{eq:grass5}
\mbox{Hom}_{\CC}(V_S/V_{S_1}, V_{S_2} / V_S).
\end{equation}
Thus, if $\{ i \} = S-S_1$ and $\{ j \} = S_2 -S$, this tangent space has
$v_i^* \otimes v_j$ as basis vector with weight $\alpha_j - \alpha_i$. 
Thus the tangent spaces to these spheres account for all the weights on the 
list (\ref{eq:grass2}).
\end{proof}

From the result above we get the following description of the graph, 
$\Gamma$:

\begin{theorem}
\label{th:grassGamma}
The vertices of the graph, $\Gamma$, are in one-to-one correspondence with
the $k$-element subsets, $S$, of $\{1,..,n\}$ via the map $S \to p_S$; 
two vertices $p_S$ and $p_{S'}$ are adjacent if $\#(S \cap S')=k-1$.
\end{theorem}

The graph we just described is called \emph{the Johnson graph} 
and is a familiar object in graph theory; 
see for instance \cite{BCN}.
The axial function, $\alpha$, and the connection, $\theta$, 
are  easy to decipher
from the results above: Let $e$ be an oriented edge joining the vertex 
$p_{S}=i(e)$ to the vertex $p_{S'}=t(e)$. Then, by (\ref{eq:grass5}), 
\begin{equation}
\label{eq:grass6}
\alpha_e = \alpha_j -\alpha_i,
\end{equation}
with $\{i \}= S-S'$ and $\{j \} = S'-S$; and these identities determine the
axial function, $\alpha$. As for the connection, $\theta$, we 
note that since the axial function, $\alpha$, has the 3-independence 
property of Theorem \ref{th:3indep}, there is a unique connection on
$\Gamma$ which is {\em compatible} with $\alpha$ in the sense that 
$\alpha$ and $\theta$ satisfy the properties of Theorem \ref{th:aximan}.
Thus all we have to do is to produce a connection which satisfies these
hypotheses, and we leave it to the reader to check that the following 
connection does: Let $p=p_S$ and $p'=p_{S'}$ be adjacent vertices 
with $\{i \}= S-S'$ and $\{j \} = S'-S$; and let $e$ be the oriented edge 
joining $p$ to 
$p'$. By Theorem \ref{th:grass=gkm}, the set of edges, $E_p$, 
can be identified
with the set of pairs, $(i,j) \in S \times S^c$, and $E_{p'}$ can be 
identified with the set of pairs, $(i',j') \in S' \times (S')^c$.
Define $\theta_e : E_p \to E_{p'}$ to be the map that sends:
\begin{equation}
\label{eq:grass7}
\begin{cases}
(k,l) \text{ to } (k,l) & \text{ if $k \neq i$ and $l \neq j$} \\
(i,l) \text{ to } (j,l) & \text{ if $l \neq j$} \\
(k,j) \text{ to } (k,i
) & \text{ if $k \neq i$} \\
(i,j) \text{ to } (j,i) & 
\end{cases}
\end{equation}
and let $\theta$ be the connection consisting of all these maps.

We will next discuss some Morse theoretic properties of the Johnson graph.
Since the Grassmannian is a co-adjoint orbit of $SU(n)$, 
$(\Gamma, \alpha)$ has the 
no-cycle property described in Theorem \ref{th:hamnocycle}. However, it
is also easy to verify this directly: For every fixed point $p=p_S$ let
\begin{equation}
\label{eq:grass8}
\alpha_S = \sum_{i \in S} \alpha_i
\end{equation}
and note that if $e$ is an oriented edge that joins $p_{S}$ to $p_{S'}$
then
\begin{equation}
\label{eq:grass9}
\alpha_e = \alpha_{S'} - \alpha_S.
\end{equation}

As in section \ref{ssec:combbetti}, let $\P$ be the set of polarizing 
elements of $\fg$: $\xi \in \P$ if and only if $\alpha_e(\xi) \neq 0$
for all $e \in E_{\Gamma}$. By (\ref{eq:grass1}) and (\ref{eq:grass6})

\begin{equation}
\label{eq:grass10}
\xi = (\xi_1,...,\xi_n) \in \P \Longleftrightarrow \xi_i \neq \xi_j
\end{equation}
and by (\ref{eq:grass9}) it is clear that if $\xi \in \P$ then the function
\begin{equation}
\label{eq:grass11}
\phi^{\xi} : V_{\Gamma} \to \RR \quad , \quad 
\phi^{\xi}(p_S) =  \alpha_S(\xi)
\end{equation}
is $\xi$-compatible.

A particularly apposite choice of $\xi$ is $\xi_i=  i$, $i=1,..,n$. 
We claim that for this choice of $\xi$ we have:

\begin{theorem}
\label{th:selfindex}
The function 
$\phi=\phi^{\xi}$ is \emph{self-indexing} modulo an additive constant:
\begin{equation}
\label{eq:grass12}
\phi(p_S) = index(p_S) + \frac{k(k+1)}{2}.
\end{equation}
\end{theorem}

\begin{proof}
The index of $p_S$ is the number of edges $e \in E_{p_S}$ with
$\alpha_e(\xi) < 0$; alternatively, it is the number of pairs 
$(i,j) \in S \times S^c$ with $\alpha_j(\xi) -\alpha_i(\xi) < 0$,
which is the same as $j-i <0$. Let $i_1 < i_2< ... < i_k$ be the
elements of $S$. The number of elements $j \in S^c$ with 
$j < i_1$ is $i_1-1$; the number of elements $ j \in S^c$ with $j< i_2$ 
is $i_2-2$ and so on; therefore the number of pairs 
$(i,j) \in S \times S^c$ with $j <i$ is 
$(i_1 + ... + i_k) - k(k+1)/2 = \phi(p_S) - k(k+1)/2$.
\end{proof}

We will conclude this description of the Johnson graph by saying a few 
words about the cohomology ring $H(\Gamma,\alpha)$. Let's introduce a 
partial ordering on $V_{\Gamma}$ by decreeing that for adjacent vertices
$p$ and $p'$ 
\begin{equation}
\label{eq:grass13}
p \prec p' \Longleftrightarrow \phi(p) < \phi(p')
\end{equation}
and, more generally, for any pair of vertices $p$ and $p'$, 
$p \prec p'$ if there exists a sequence of adjacent vertices
$$p=p_0 \prec p_1 \prec ... \prec p_r=p'.$$
We will prove in section \ref{ssec:generators} (as a special case of a 
more general theorem) that $H(\Gamma,\alpha)$ is a free module over
$\SS(\fg^*)$, with generators $\{ \tau_p \; ; \; p \in V_{\Gamma} \}$
uniquely characterized by the following two properties:
\begin{enumerate}
\item $\frac{1}{2}deg(\tau_p) = index(p)$
\item The support of $\tau_p$ is contained in the set
\begin{equation}
\label{eq:grass14}
F_p = \{ q \in V_{\Gamma} \; ; \; p \prec q \}.
\end{equation}
\end{enumerate}

To reconcile this result with classical results of Kostant, Kumar and others 
on the cohomology ring of the Grassmannian, we will also give in 
section \ref{ssec:schubs} an alternative description of $\tau_p$ in terms
of the Hecke algebra of divided difference operators; for this we will need 
an alternative description of the ordering (\ref{eq:grass13}). One 
property of the Johnson graph which we haven't yet commented on is that it 
is a {\em symmetric} graph. Given two pairs of adjacent vertices 
$(p,p')$ and $(q,q')$, one can find a permutation $\sigma \in S_n$ with
$\sigma (p) = p'$ and $\sigma(q)=q'$. We claim that the partial ordering
(\ref{eq:grass13}) is equivalent to the so-called \emph{Bruhat order} 
on $V_{\Gamma}$ (see \cite{Hu}):

\begin{theorem}
\label{th:Bruhatorder}
If $p$ and $p'$ are vertices of $\Gamma$, then $p \prec p'$ if and only if 
there exists a sequence of elementary reflections
\begin{equation}
\label{eq:grass15}
\sigma_i \quad : \quad i \leftrightarrow i+1
\end{equation}
with $i=i_1,...,i_m$ such that 
\begin{equation}
\label{eq:grass16}
m = index(p') - index(p),
\end{equation}
\begin{equation}
\label{eq:grass17}
p' = \sigma_{i_m} \circ ... \circ \sigma_{i_1}(p)
\end{equation}
and such that $\phi$ is strictly increasing along the sequence of
adjacent vertices
\begin{equation}
\label{eq:grass18}
p_k =  \sigma_{i_k} \circ ... \circ \sigma_{i_1}(p), \quad k =1,...,m.
\end{equation}
\end{theorem}

For a proof of this see for instance \cite{GHZ}.

\section{\bf Abstract one-skeleta}
\label{sec:abstract}

\subsection{\bf Abstract one-skeleta}
\label{sec:skeletons}

If one strips the manifold scaffolding from GKM theory, one gets a 
graph-like object which we will call an abstract one-skeleton.
Let $\fg^*$ be an arbitrary $n$-dimensional vector space. 

\begin{definition}
\label{def:abs1skel}
An \emph{abstract one-skeleton} 
is a triple consisting of a $d$-valent graph, $\Gamma$ 
(with $V_{\Gamma}$ as vertices and $E_{\Gamma}$ as
oriented edges), a connection, $\theta$, on the ``tangent bundle'' 
of $\Gamma$
$$\pi : E_{\Gamma} \to V_{\Gamma}, \qquad \pi(e) = i(e) 
 \;\;\; \mbox{ (the initial vertex of} e)$$
and an axial function 
$$\alpha : E_{\Gamma} \to \fg^*$$
satisfying the axioms:
\begin{enumerate}[{A}1]
\item  \label{axiom1} 
For every $p \in V_{\Gamma}$, the vectors 
$\{ \alpha_e \; ; \; e \in E_p = \pi^{-1}(p) \}$ 
are pairwise linearly independent.
\item \label{axiom2}
If $e$ is an oriented edge of $\Gamma$, 
and $\bar{e}$ is the same edge with its orientation reversed, 
there exist positive numbers, $m_{e}$ and $m_{\bar{e}}$, such that
\begin{equation}
\label{eq:axiom2}
m_{\bar{e}}\alpha_{\bar{e}} = - m_{e}\alpha_{e}.
\end{equation}
\item \label{axiom3} Let $e\in E_{\Gamma}$, $p=i(e)$ and $p'=t(e)$. 
Let $e_i$, $i=1,..,d$, be the elements of $E_p$ and $e_i'$, $i=1,..,d$, their
images with respect to $\theta_e$ in $E_p'$. Then
\begin{equation}
\label{eq:axiom3}
\alpha_{e_i'} = \lambda_{i,e} \alpha_{e_i} + c_{i,e} \alpha_e
\end{equation}
with $\lambda_{i,e}> 0$ and $c_{i,e} \in \RR$.
\end{enumerate} 
\end{definition}

\begin{remarks*}
\begin{enumerate}

\item We will denote an abstract one-skeleton by $(\Gamma, \alpha)$.

\item Axioms A1-A3 imply that $\theta_e(e) = \bar{e}$.

\item There is a natural notion of equivalence for axial functions: Let
$\theta$ be a connection on a graph $\Gamma$ and let $\alpha$ and $\alpha'$
be axial functions. We will say that $\alpha$ and $\alpha'$ are 
\emph{equivalent axial functions} if for every oriented edge $e$,
\begin{equation}
\label{eq:eqaxi}
\alpha_{e}' = \lambda_{e}\alpha_{e}, 
\mbox{ with } \lambda_{e} > 0.
\end{equation}

\item We can always replace an axial function, $\alpha$, by an equivalent
axial function for which the constants $m$ in \eqref{eq:axiom2} 
are 1, {\em i.e.} we can assume
\begin{equation} \label{eq:axiom2`}
\alpha_{\bar{e}} = - \alpha_{e}.
\end{equation}

\item One can define the Betti numbers, $b_{2i}(\Gamma)$, and the 
cohomology ring, $H(\Gamma,\alpha)$, of an abstract one-skeleton exactly 
as in sections
\ref{ssec:combbetti} and \ref{ssec:thecohring}. (It is easy to check, by 
the way, that in our proof of the well-definedness of $b_{2i}(\Gamma)$
(Theorem \ref{th:2.5}) 
we can replace item 5 of Theorem \ref{th:aximan} by the somewhat weaker
hypothesis \eqref{eq:axiom3})

\item It is also clear that the definition of $b_{2i}(\Gamma)$ and of
$H(\Gamma,\alpha)$ is unchanged if we replace $\alpha$ by an equivalent 
axial function, $\alpha'$.

\item Let $\fg^*$ be, as in the first part of this paper, the dual of the Lie 
algebra of an $n$-dimensional torus, $G$. Suppose that for every 
$e \in E_{\Gamma}$, $\alpha_e$ is an element of the weight lattice of $G$.
We will say that the abstract one-skeleton $(\Gamma,\alpha)$ is
an \emph{abstract GKM one-skeleton} if the 
$m$'s in \eqref{eq:axiom2} and the 
$\lambda$'s in \eqref{eq:axiom3} are all equal to 1, that is
\begin{equation} \label{eq:axiom2gkm}
\alpha_{\bar{e}} = - \alpha_{e}.
\end{equation}
and
\begin{equation}
\label{eq:axiom3gkm}
\alpha_{e_i'} = \alpha_{e_i} + c_{i,e} \alpha_e
\end{equation}
and the $c_{i,e}$'s are integers. We will show in section 
\ref{ssec:realization} that every abstract GKM one-skeleton is actually the 
GKM one-skeleton of a GKM manifold.

\end{enumerate}
\end{remarks*}

\begin{definition}
\label{def:3indep}
We will say that an axial function, $\alpha$, is \emph{three-independent}
if, for every $p \in V_{\Gamma}$, the vectors 
$\{ \alpha_e;  e \in E_p \}$, are
3-independent in the sense of Theorem~\ref{th:reduced}.
\end{definition}

It is clear that if $\alpha$ and $\alpha'$ are equivalent axial functions 
and one of them is three-independent the other is as well. The hypothesis of
three-independence will be frequently evoked in this chapter. 
It will enable us 
to blow-up and blow-down abstract one-skeleta and, by mimicking 
Theorem \ref{th:reduced}, to define an analogue of symplectic reduction
for abstract one-skeleta. 
It also rules out the existence of 2-cycles in $\Gamma$
(such as the three 2-cycles exhibited in Figure \ref{fig:s6}.)

\begin{proposition} \label{2cycles}
If $\alpha$ is three-independent every pair of adjacent vertices in $\Gamma$ 
is connected by a \emph{unique} unoriented edge.
\end{proposition}

\begin{proof}
Suppose that there are two distinct oriented edges, $e$ and $e_1$, 
from $p$ to $p'$; let $e'=\theta_{e_1}(e) \in E_{p'}$. 
Since $ \alpha_{\bar{e}} = - \alpha_e$ and 
$\alpha_{e'} = \lambda \alpha_e + c \alpha_{e_1}$, 
with $\lambda >0$, it 
follows that $e' \neq \bar{e}$. Thus the vectors $\alpha_{e'}, 
\alpha_{\bar{e}}$, and $\alpha_{\bar{e_1}}$ are distinct and coplanar, 
which contradicts the three-independence of $\alpha$ at $p'$.
\end{proof}

Another useful consequence of three-independence is the following.

\begin{proposition}
If $\fh$ is a codimension 2 subspace of $\fg$, the graph $\Gamma_{\fh}$
is 2-valent.
\end{proposition}

Finally, if $\alpha$ is three-independent, the compatibility conditions between
$\theta$ and $\alpha$ imposed by Axiom A3  \emph{determine}
$\theta$.

\begin{proposition}
The connection, $\theta$, is the only connection on $\Gamma$ satisfying 
\eqref{eq:axiom3}.
\end{proposition}

The GKM-theorem asserts that if $M$ is a GKM manifold and is equivariantly 
formal then $H_G(M)$ is isomorphic to $H(\Gamma,\alpha)$. In particular, 
$H(\Gamma,\alpha)$ is a free module over the ring $\SS(\fg^*)$ with 
$b_{2i}(M)$ generators in dimension $2i$. In the second part of this paper
we will attempt to ascertain:
`` To what extent is this theorem true for abstract one-skeleta with
$b_{2i}(M)$ replaced by $b_{2i}(\Gamma)$?''. The examples we've encountered 
in the first part of this paper (section \ref{ssec:examples}) 
already give us some inkling of what to expect: this 
assertion is {\em unlikely} to be true if for some admissible orientation
of $\Gamma$ (see section \ref{ssec:hamiltgkm}) there exist oriented closed 
paths, or if for some subspace $\fh$ of $\fg$, the totally geodesic
subgraph $\Gamma_{\fh}$ of $\Gamma$ has fewer connected components than
predicted by its combinatorial Betti number.
This motivates the following definition.
\begin{definition} \label{def:nca}
The abstract one-skeleton $(\Gamma,\alpha)$ is \emph{non-cyclic} if
\begin{enumerate}[{NCA}1]
\item \label{axiom:nca1}
For some vector $\xi$ in the set \eqref{eq:polset}, 
$(\Gamma, \alpha)$ is $\xi$-acyclic, 
{\em i.e.} the oriented graph $(\Gamma, o_{\xi})$ has no closed paths.
\item \label{axiom:nca2}
For every codimension 2 subspace, $\fh$, of $\fg$ and for every connected 
component, $\Gamma_0$, of $\Gamma_{\fh}$
\begin{equation}
\label{eq:axiomnca2}
b_0(\Gamma_0) =1.
\end{equation}
\end{enumerate}
\end{definition}
A reminder: For the definition of $o_{\xi}$ see section 
\ref{ssec:hamiltgkm}. Also recall that by Theorem \ref{th:posnocycle},
$\xi$-acyclicity implies the existence of a function 
$f : V_{\Gamma} \to \RR$ which is $\xi$-compatible.

In the remaining of Section \ref{sec:abstract}, by one-skeleton 
we will mean an \emph{abstract} one-skeleton.

\subsubsection{\bf Examples}
\label{sec:examples}

\begin{example} 
{\bf The complete one-skeleton:} 
\label{ssec:complete}

In this example the vertices of $\Gamma$ are the elements of the 
$N$-element set $V=\{ p_1,\dots , p_N \}$ and each pair of elements, 
$(p_i,p_j),\; i \neq j$,
is joined by an edge. We will denote by $e=p_ip_j$ the oriented edge 
that joins $p_i=i(e)$ to $p_j=t(e)$. Thus the set of oriented edges is 
just the set
$$ \{ p_ip_j \;; \; 1 \leq i,j \leq N, \; i \neq j \}$$
and its fiber over $p_i$ is
$$ E_i = \{ p_ip_j \; ; \; 1 \leq j \leq N, \; i\neq j \}.$$
A connection, $\theta$, is defined by maps,
$\theta_{ij} :E_i \to E_j$, where
$$\theta_{ij}(p_ip_k) = \left\{ 
\begin{array}{lcl} p_jp_i & \mbox{ if } & k =j \\
                   p_jp_k & \mbox{ if } & k \neq i,j 
\end{array} \right. $$
(Note that this connection is invariant under all permutations of
the vertices.)

Let $\tau : V \to \fg^* $ be any function such that $\tau_1,..,\tau_N$
are 3 independent; then the function, $\alpha$, given by
\begin{equation}
\label{eq:axcompl}
\alpha_{p_ip_j} =  \tau_j - \tau_i
\end{equation}
is an axial function compatible with $\theta$. We will call 
$\tau : V \to \SS^1(\fg^*)$ \emph{the generating class of } $\Gamma$. 

The following theorem describes the additive structure of the cohomology 
ring of $(\Gamma, \alpha)$. If $\dim{(\fg^*)}=n$ let
\begin{equation}
\label{eq:4.2}
\lambda_k = \lambda_{k,n} = \left\{ \begin{array}{cc}
\dim \SS^k(\fg^*) & , \mbox{ if } k \geq 0 \\
0               & , \mbox{ if } k < 0 
\end{array} \right.
\end{equation}

\begin{theorem}
\label{th:cohcompl}
If $(\Gamma, \alpha)$ is the complete one-skeleton with $N$ vertices and 
generating class $\tau$, then
\begin{equation*}
H^{2m}(\Gamma,\alpha) \simeq \bigoplus_{k=0}^{N-1} 
\SS^{m-k}(\fg^*)\, \tau^k 
\end{equation*}
for every $m \geq 0$. In particular,
\begin{equation}
\label{eq:4.3}
\dim H^{2m}(\Gamma,\alpha) = \sum_{k=0}^{N-1} \lambda_{m-k}.
\end{equation}
\end{theorem}

\begin{proof}

The generating class $\tau \in H^2(\Gamma,\alpha)$ satisfies the relation
\begin{equation}
\label{eq:4.4}
\tau^N = \sigma_1(\tau_1,...,\tau_N) \tau^{N-1} -
\sigma_2(\tau_1,...,\tau_N) \tau^{N-2} + ... \; ,
\end{equation}
where $\sigma_k(\tau_1,...,\tau_N) \in \SS^k(\fg^*)$ 
is the $k^{th}$ symmetric polynomial in $\tau_1,..,\tau_N$.

We will show that every element $f \in H^{2m}(\Gamma,\alpha)$ can be 
written uniquely as
\begin{equation}
\label{eq:4.5}
f= \sum_{k=0}^{N-1} f_{m-k}\tau^k ,
\end{equation}
with $f_{m-k} \in \SS^{m-k}(\fg^*)$ if $k \leq m$ and $f_{m-k}=0$
if $k>m$.

For $m=0$ the statement is obvious. Assume $m >0$ and let 
\begin{equation}
\label{eq:4.7}
g_m  = 
(-1)^{N+1} \tau_1 \cdots \tau_N \sum_{i=1}^N
\frac{f(p_i)}{\tau_i \prod_{j \neq i} (\tau_i - \tau_j)} .
\end{equation}
A priori, $g_m$ is an element in the field of fractions of $S(\fg^*)$.
Since $f \in H(\Gamma,\alpha)$, $\tau_i -\tau_j$ divides $f(p_i)-f(p_j)$
for all $ i \neq j$, and hence all the factors in the denominator of 
$g_m$ will be canceled so $g_m \in \SS^m(\fg^*)$. 
Moreover, from 
(\ref{eq:4.7}) follows that $f(p_i)-g_m \equiv 0$ on $\tau_i \equiv 0$;
therefore there exists $h_i \in \SS^{m-1}(\fg^*)$ such that 
\begin{equation*}
f(p_i)=g_m + \tau_ih_i \quad , \quad \forall \; i=1,..,N.
\end{equation*}
Since
\begin{equation*}
f(p_i)-f(p_j) =(\tau_i -\tau_j) h_i + \tau_j (h_i -h_j), 
\quad \forall i \neq j,
\end{equation*}
it follows that $\tau_i -\tau_j$ divides $h_i - h_j$ for all $i \neq j$,
that is, the function $h : V \to \SS^{m-1}(\fg^*)$, given by $h(p_i)=h_i$,
satisfies the compatibility conditions and is therefore an element of
$H^{2(m-1)}(\Gamma,\alpha)$. Thus
\begin{equation}
\label{eq:4.8}
f=g_m +\tau h ,
\end{equation}
with $h \in H^{2(m-1)}(\Gamma,\alpha)$. From the induction hypothesis, $h$ 
can be uniquely written as 
\begin{equation}
\label{eq:4.9}
h=h_{m-1}+h_{m-2}\tau + ... +h_{m-N}\tau^{N-1}.
\end{equation}
Introducing (\ref{eq:4.9}) in (\ref{eq:4.8}) and using (\ref{eq:4.4})
we deduce that $f$ can be written in the form (\ref{eq:4.5}) and the 
uniqueness follows from the non-degeneracy of the Vandermonde determinant 
with entries $(\tau_i^k)$.
If $m<k$ then the corresponding $f_{m-k}$ would have negative degree; 
so the only possibility is that it is zero.
\end{proof}

\begin{remark*}
This example is associated with a simple (but very important) 
GKM action, the action of $T^N$ on $\CC P^{N-1}$.
\end{remark*}

\end{example}

\begin{example}
{\bf Sub-skeleta:}
\label{ssec:subsk}

Let $\Gamma_0$ be an $r$-valent subgraph of $\Gamma$ which is totally 
geodesic in the sense of Definition \ref{def:totgeod}. Then the restriction 
of $\theta$ and $\alpha$ to $\Gamma_0$ define a connection, $\theta_0$, and 
and axial function, $\alpha_0$, on $\Gamma_0$; and we will call 
$(\Gamma_0, \alpha_0)$ a \emph{sub-skeleton} of $(\Gamma, \alpha)$. 

Associated to sub-skeleton is the notion of normal holonomy.
Let $\Gamma_0$ be a totally geodesic subgraph 
of $\Gamma$ and $\theta_0$ the 
connection on $\Gamma_0$ induced by $\theta$. For a vertex 
$p \in V_{\Gamma_0}$, let $E^0_p$ be the fiber of $E_{\Gamma_0}$ over $p$
and $N_p = E_p - E^0_p$. For every loop $\gamma$ in $\Gamma_0$
based at $p$, the map $\sigma_\gamma$ preserves the decomposition
$E_p = E^0_p \cup N_p$ and thus induces a permutation $\sigma^0_\gamma$
of $N_p$. Let $Hol^{\bot}(\Gamma,\Gamma_0,p)$ be the subgroup of
$\Sigma(N_p) \subset \Sigma(E_p)$
generated by the 
permutations $\sigma^0_\gamma$, for all loops $\gamma$ included in 
$\Gamma_0$ and based at $p$. Again, if $p_1$ and $p_2$ are connected by
a path in $\Gamma_0$, then $Hol^{\bot}(\Gamma,\Gamma_0,p_1)$ and 
$Hol^{\bot}(\Gamma,\Gamma_0,p_2)$
are isomorphic by conjugacy; so, if $\Gamma_0$ is connected, we can define 
\emph{the normal holonomy group}, $Hol^{\bot}(\Gamma_0,\Gamma)$,
of $\Gamma_0$ in $\Gamma$ to be
$Hol^{\bot}(\Gamma,\Gamma_0,p)$ for any $p \in \Gamma_0$. 
We will also say that $\Gamma_0$ has
\emph{trivial normal holonomy in $\Gamma$} if 
$Hol^{\bot}(\Gamma_{0i},\Gamma)$ is trivial for every connected component
$\Gamma_{0i}$ of $\Gamma_0$.
\end{example}

\begin{example}
{\bf Product one-skeleta:}
\label{ssec:prodsk}

Let $(\Gamma_i,\alpha_i)$, $i=1,2$, be a $d_i$-valent one-skeleta. 
The vertices of the product graph 
$\Gamma =\Gamma_1 \times \Gamma_2$, are the 
pairs, $(p,q)$, $p \in V_{\Gamma_1}$ and $q \in V_{\Gamma_2}$;
and two vertices, $(p,q)$ and $(p',q')$, are joined by an edge if either
$p=p'$ and $q$ and $q'$ are joined by an edge in $\Gamma_2$ or
$q=q'$ and $p$ and $p'$ are joined by an edge in $\Gamma_1$. 
(If $p$ is joined to $p'$ by several edges, each of these edges will 
correspond to an edge joining $(p,q)$ to $(p',q)$. We will, however, be a 
bit careless about this fact in the paragraph below and denote (oriented) 
edges by the pairs of adjacent vertices they join.)

If $\theta_1$ and $\theta_2$ are
connections on $\Gamma_1$ and $\Gamma_2$, {\em the product
connection}, 
$\theta$, on $\Gamma_1 \times \Gamma_2$ is defined by 
$$\theta_{(p,q;p',q)}  =  \theta_{p,p'} \times (Id)_{q,q} 
\quad \mbox{ and } \quad
\theta_{(p,q;p,q')}  =  (Id)_{p,p} \times \theta_{q,q'}$$
and one can construct an axial function on $\Gamma$ compatible with 
$\theta$ by defining 
\begin{equation*}
\alpha((p,q),(p',q')) = \left\{ \begin{array}{ll}
\alpha_1(p,p') & \mbox{ if } q=q' \mbox{ and } (p,p') \in E_1 \\ 
\alpha_2(q,q') & \mbox{ if } p=p' \mbox{ and } (q,q') \in E_2
\end{array} \right. 
\end{equation*}
Then  $(\Gamma,\alpha)$, is a $d_1+d_2$-valent one-skeleton 
which we will call {\it the direct product of $\Gamma_1$ and $\Gamma_2$}.
\end{example}

\begin{example}
{\bf Twisted products:}
\label{ssec:twprod}

One can extend the results we've just described to twisted
products. Let $\Gamma_0$ and $\Gamma'$ be two graphs and 
$\psi : E_{\Gamma_0} \to Aut(\Gamma')$ a map such that 
$$\psi(e) = (\psi(\bar{e}))^{-1}$$
for every oriented edge $e$. We will define 
\emph{the twisted product of $\Gamma_0$ and $\Gamma'$},
$$\Gamma =\Gamma_0 \times_{\psi} \Gamma', $$
as follows. The set of vertices of this new graph is 
$V_{\Gamma} = V_{\Gamma_0} \times V_{\Gamma'}$.

Two vertices $p_1=(p_{01},p_1')$ and $p_2=(p_{02},p_2')$ are joined 
by an edge iff 
\begin{enumerate}
\item $p_{01}=p_{02}$ and $p_1', p_2'$ are joined by an edge of $\Gamma'$
(these edges will be called {\em vertical}) or
\item $p_{01}$ is joined with $p_{02}$ by an edge, $e\in E_{\Gamma_0}$ and 
$p_2'=\psi(e) (p_1')$ 
(these edges will be called {\em horizontal}).
\end{enumerate}
For a vertex $p = (p_0,p') \in V_{\Gamma}$ 
we denote by $E_p^h$ the set of oriented horizontal edges issuing from 
$p$ and by $E_p^v$ the set of oriented vertical edges issuing from $p$. 
Then the projection
$$\pi : V_{\Gamma} = V_{\Gamma_0} \times V_{\Gamma'} \to V_{\Gamma_0}$$ 
induces a 
bijection $d\pi_{p} : E_p^h \to E_{p_0}$. 
 
Let $q_0 \in V_{\Gamma_0}$ and $\Omega(\Gamma_0, q_0)$ be 
{\em the fundamental group of} $\Gamma_0$, that is, the set of all 
loops based at $q_0$. Every such loop, $\gamma$, induces an element 
$\psi_{\gamma} \in Aut(\Gamma')$ by composing the automorphisms 
corresponding to its edges. Now let $G$ be the subgroup of $Aut(\Gamma')$ 
which is the image of the morphism
$$\Psi : \Omega(\Gamma_0, q_0) \to Aut(\Gamma'), 
\qquad \Psi(\gamma)=\psi_{\gamma}.$$

If $\theta_0$ is a connection on $\Gamma_0$ and $\theta'$ is 
a $G$-invariant connection on $\Gamma'$ we get 
a connection on the twisted product as follows:

If we require this connection 
to take horizontal edges to horizontal edges 
and vertical edges to vertical edges, we must decide how these 
horizontal and vertical components are related at adjacent points of $V$. 
For the horizontal component, the relation is simple: By assumption, the map
$$d\pi_{p} : E_{p}^h \to E_{p_0}, \qquad (\mbox{ where } p=(p_0,p'))$$
is a bijection; so if $p_1$ and $p_2$ are adjacent points on the same fiber
above $p_0$,
$$ \theta_{p_1p_2} = d\pi_{p_2}^{-1} \circ d\pi_{p_1} : 
E_{p_1}^h \stackrel{\sim}{\to} E_{p_2}^h \simeq E_{p_0} \; ,$$
and if $p_1=(p_{01},p_1')$ and $p_2=(p_{02},p_2')$ 
are adjacent points on different fibers,
$$ \theta_{p_1p_2} = d\pi_{p_2}^{-1} \circ \theta_{p_{01}p_{02}} \circ
d\pi_{p_1},$$
so that the following diagram commutes
\begin{displaymath}
\begin{CD}
E_{p_1}^h  @>>> E_{p_2}^h \\
@V{\simeq}VV    @VV{\simeq}V  \\
E_{p_{01}} @>{\theta_{p_{01}p_{02}}}>> E_{p_{02}}  
\end{CD}
\end{displaymath}

For the vertical component the relation is nearly as simple. Let 
$p_0 \in V_{\Gamma_0}$ and $\gamma$ a path in $\Gamma_0$ from 
$q_0$ to $p_0$. Then $\psi$ induces an automorphism 
$$\Psi_{\gamma} : \pi^{-1}(q_0) \to \pi^{-1}(p_0) \simeq V_{\Gamma'},$$
and we use $\Psi_{\gamma}$ to define a connection $\theta_{p_0}'$ 
on $\pi^{-1}(p_0)$.
Since $\theta'$ is $G$-invariant, this new connection is actually 
independent of $\gamma$. 

Thus if $p_1$ and 
$p_2$ are adjacent points on the fiber above $p_0$, the connection on
$\Gamma$ along $p_1p_2$ is induced on vertical edges by $\theta_{p_0}'$. 
On the other hand, if $p_1$ and $p_2$ are adjacent points on different
fibers and $p_{0i}=\pi(p_i)$, then $\psi_{p_{01}p_{02}}$ induces a map  
$E_{p_1}^v \to E_{p_2}^v$ 
that will be the connection on vertical edges.

Note that if $\Psi$ is trivial then the twisted product is actually 
a direct product.

The axial functions on $\Gamma$ which we will encounter in section
\ref{sec:blowup} {\em won't} as a rule be of the product form described in 
Example \ref{ssec:prodsk}. They will, however, satisfy 
\begin{equation}
\label{eq:1.11}
\alpha(e) = \alpha_0(e_0), \quad e_0=d\pi_{p}(e),
\end{equation}
at all $p \in V_{\Gamma}$ and $e \in E_{p}^h$, $\alpha_0$ being
a given axial function on $\Gamma_0$.
\end{example}

\begin{example}
{\bf Fibrations:}
\label{ssec:fibrations}

An important example of twisted products is given by fibrations. 
Let $\Gamma$ and $\Gamma_0$ be graphs of valence $d$
and $d_0$ and let $V$ and $V_0$ be their vertex sets.

\begin{definition}
\label{def:morph}
A \emph{morphism} of $\Gamma$ into $\Gamma_0$ is a map
$f : V \to V_{0}$ 
with the property that if $p$ and $q$ are adjacent 
points in $\Gamma$ then  either $f(p)=f(q)$ or $f(p)$ and $f(q)$ are 
adjacent points in $\Gamma_{0}$.
\end{definition}

Let $p \in V$, $p_0=f(p)$, let $E_{p}^v$ be the set of oriented edges, 
$pq$, with $f(p) = f(q)$, and 
let $E_{p}^h = E_{p} - E_p^v$. (One can regard 
$E_{p}^h$ as the ``horizontal'' component of $E_{p}$ and 
$E_{p}^v$ as the ``vertical'' component.) 
By definition there is a map
\begin{equation}
\label{eq:1.5}
df_{p} : E_{p} - E_{p}^v \to E_{p_0}
\end{equation}
which we will call \emph{the derivative of} $f$ {\em at} $p$.

\begin{definition}
A morphism, $f$, is a {\em submersion at }$p$ if the map, $df_{p}$,
is bijective. If $df$ is bijective
at all points of $V$ then we will simply say that $f$ is a 
{\em submersion}.
\end{definition}

\begin{theorem}
\label{th:submersion}
Let $\Gamma$ and $\Gamma_0$ be connected and let $f$ be a submersion.
Then:
\begin{enumerate}
\item $f$ is surjective ;
\item For every $p \in V_0$, the set $V_{p}=f^{-1}(p)$ is the vertex
	set of a subgraph of valence $r=d-d_0$ ;
\item If $p,q \in V_0$ are adjacent, there is a canonical bijective map
\begin{equation}
\label{eq:1.6}
K_{p,q} : V_{p} \to V_{q} 
\end{equation}
defined by 
\begin{equation}
\label{eq:1.7}
K_{p,q}(p') = q' \quad \Longleftrightarrow \quad p' \mbox{ and } q' 
\mbox{ are adjacent} ;
\end{equation}
\item In particular, the cardinality of $V_{p}$ is the same for all $p$.
\end{enumerate}
\end{theorem}

We will leave the proofs of these assertions as easy exercises. Note by 
the way that the map (\ref{eq:1.6}) satisfies 
\begin{equation*}
K_{p,q}^{-1} = K_{q,p}.
\end{equation*}

\begin{definition}
The submersion $f$ is a {\em fibration} if the map (\ref{eq:1.6}) 
preserves adjacency: two points in $V_{p}$ are adjacent if and only if
their images in $V_{q}$ are adjacent.
\end{definition}

Let $f$ be a fibration,  $p_0$ be a base point in $V_0$ 
and $p$ any other point. Let $\Gamma_{p_0}$ and 
$\Gamma_p$ be the subgraphs of $\Gamma$ (of valence $r$) whose vertices are 
the points of $V_{p_0}$ and $V_p$. For every path
\begin{equation*}
p_0 \rightarrow p_1 \rightarrow \dots \rightarrow p_N = p
\end{equation*}
in $\Gamma_0$ joining $p_0$ to $p$, there is a holonomy map
\begin{equation}
\label{eq:1.8}
K_{p_{N-1},p_N} \circ \dots \circ K_{p_1,p_0} : V_{p_0} \to V_p
\end{equation}
which preserves adjacency. Hence all the graphs $\Gamma_p$ are isomorphic
(and, in particular, isomorphic to $\Gamma' :=\Gamma_{p_0}$). 
Thus $\Gamma$ can be 
regarded as a twisted product of $\Gamma_0$ and $\Gamma'$. 
Moreover, for every closed path 
\begin{equation*}
\gamma : \quad p_0 \rightarrow p_1 \rightarrow \dots 
\rightarrow p_N = p_0
\end{equation*}
there is a holonomy map
\begin{equation}
\label{eq:1.9}
K_{\gamma} : V_{p_0} \to V_{p_0},
\end{equation}
and it is clear that if this map is the identity for all $\gamma$,
this twisted product is a direct product, that is 
$$\Gamma \simeq \Gamma_0 \times \Gamma' $$

From (\ref{eq:1.9}) one gets a homomorphism of the fundamental group of
$\Gamma_0$ into $Aut(\Gamma')$. Let $G$ be its image. 
Given a connection on
$\Gamma_0$ and a $G$ invariant connection on $\Gamma'$, one 
gets, by the construction described in the previous example,
a connection on $\Gamma$. 
\end{example}

From now on, unless specified otherwise, we will assume that 
$(\Gamma, \alpha)$ is 3-independent. Also, frequently we will
refer to the one-skeleton $(\Gamma, \alpha)$ simply as $\Gamma$.

\subsection{\bf The blow-up operation}
\label{sec:blowup}

\subsubsection{\bf The blow-up of a one-skeleton}
\label{ssec:blowconstr}

Let $(\Gamma,\alpha)$ be a $d$-valent one-skeleton,
and let $\Gamma_0$ be a sub-skeleton of valence $d_0$ and covalence 
$s=d-d_0$. We will define in this
section a new $d$-valent one-skeleton, $(\Gamma^{\#}, \alpha^{\#})$, 
which we will call 
the blow-up of $\Gamma$ along $\Gamma_0$; we will also define a 
blowing-down map 
$$\beta : V_{\Gamma^{\#}} \to V_{\Gamma},$$
which will be a morphism of graphs in the sense of  
Definition \ref{def:morph}. The singular locus of this blowing-down 
map can be described as a twisted product of $\Gamma_0$ and a 
complete one-skeleton
on $r$ vertices; and $\Gamma^{\#}$ itself will be obtained from this
singular locus by gluing it to the complement of $\Gamma_0$ in $\Gamma$.
Here are the details:

\begin{figure}[h]
\begin{center}
\includegraphics{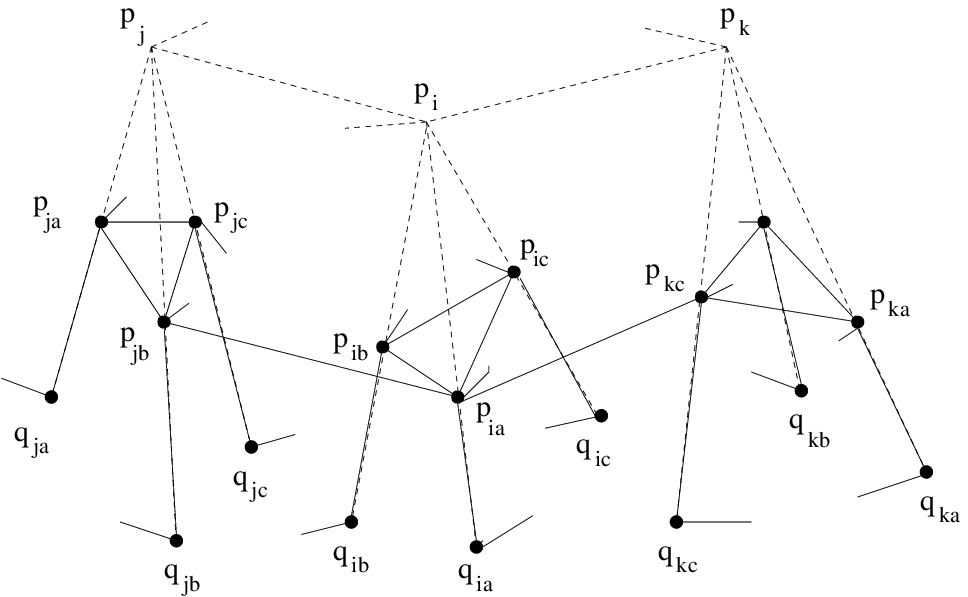}
\caption{Blow-up}
\label{fig:blow}
\end{center}
\end{figure}

Let $V_0$ and $V$ be the vertices of $\Gamma_0$ and $\Gamma$ and let
$V_2 = V - V_0$. For each $p_i \in V_0$ let 
$\{q_{ia} \; ; \; a=d_0+1, ..,d \}$ be  the 
set of points in $V_2$ which are adjacent to $p_i$. 
Define a new set of vertices, 
$N_{p_i} = \{ p_{ia} \; ; \; a=d_0+1,..,d \}$, 
with one new vertex, $p_{ia}$, for each edge $p_{i}q_{ia}$.
For each pair of
adjacent points, $p$ and $q$, in $V_0$, the holonomy map 
$\theta_{p,q} : E_p \to E_q$ induces a map 
\begin{equation*}
K_{p,q} : N_p \to N_q.
\end{equation*}
Let $V_1$ be the disjoint union of the $N_p$'s and let $f : V_1 \to V_0$ be 
the map which sends $N_p$ to $p$. We will make $V_1$ into a graph, 
$\Gamma_1$, by decreeing that two points, $p'$ and $q'$ of $V_1$, are
adjacent iff $f(p')=f(q')$ or $p=f(p')$ and $q=f(q')$ are adjacent and
$q'=K_{p,q}(p')$. It is clear that this notion of adjacency defines a 
graph, $\Gamma_1$, of valence $d-1$, and that $f$ is a fibration in the 
sense of Example \ref{ssec:fibrations}. Let $p_0$ be a base point in $V_0$. 
The subgraph $\Gamma'= f^{-1}(p_0)$ 
is a complete graph on $s$ vertices; therefore
we can equip it with the connection described in Example
\ref{ssec:complete}. This connection is invariant under {\em all} the
automorphisms of $\Gamma'$, so we can, as in Example \ref{ssec:fibrations}, 
take its twisted product with the connection on $\Gamma_0$ 
to get a connection on $\Gamma_1$. 

To define an axial function on $\Gamma_1$ which is compatible with this
connection we will have to assume that the axial function, $\alpha$, on 
$\Gamma$ satisfies a GKM hypothesis of type \eqref{eq:axiom3gkm}.
Fortunately however:
\begin{enumerate}
\item We will only have to make this assumption for the edges of $\Gamma$ 
normal to $\Gamma_0$, that is, we will only have to assume that $\alpha$
satisfies the condition
\begin{equation}
\label{eq:condition2}
\alpha_{p_iq_{ia}} - \alpha_{p_jq_{jb}} \mbox{ is a multiple of } 
\alpha_{p_ip_j}
\end{equation}
for every edge, $p_ip_j$, of $\Gamma_0$, 
where $p_jq_{jb}= \theta_{p_ip_j}(p_iq_{ia})$.

\item In the blow-up-blow-down construction in section \ref{ssec:passage}
in which we will apply the construction which we are about to describe,
the hypotheses \eqref{eq:condition2} are satisfied.
\end{enumerate}

Consider positive numbers $(n_{ia})$ 
such that 
\begin{equation}
\label{eq:condition1}
n_{ia}=n_{jb}
\end{equation}
if $ p_{ia}$ and  $p_{jb}$ are joined by an horizontal edge.
We can define an axial function, $\alpha'$, on $\Gamma_1$ as follows.

On horizontal edges of $\Gamma_1$, which are of the form
$p_{ia}p_{jb}$, 
 we will require that $\alpha'$ be
defined by (\ref{eq:1.11}), that is 
\begin{equation*}
\alpha_{p_{ia}p_{jb}}' = \alpha_{p_ip_j}.
\end{equation*}

On vertical edges, that are of the form $e'=p_{ia}p_{ib}$, 
we define $\alpha'$ by
\begin{equation*}
\alpha_{p_{ia}p_{ib}}' = \alpha_{p_iq_{ib}} -
\frac{n_{ib}}{n_{ia}}\alpha_{p_iq_{ia}}.
\end{equation*}

We will now define $\Gamma^{\#}$. Its vertices will be the set
\begin{equation*}
V^{\#} = V_1 \sqcup V_2 \; ,
\end{equation*}
and we define adjacency in $V^{\#}$ as follows:

\begin{enumerate}
\item Two points, $p'$ and $q'$, in $V_1$, are adjacent if 
	they are adjacent in $\Gamma_1$.
\item Two points, $p$ and $q$, in $V_2$, are adjacent if 
	they are adjacent in $\Gamma$.
\item Consider a point $p'=p_{ia} \in N_{p_{i}} \subset V_1$. 
	 By definition, $p'$ corresponds to a point $q_{ia} \in V_2$,
	which is adjacent to $p_i \in V$.
	Join $p'=p_{ia}$ to $q_{ia}$.
\end{enumerate}

Then $\Gamma^{\#}$ is a graph of valence $d$ and the blowing-down map 
\begin{equation*}
\beta : V^{\#} \to V 
\end{equation*}
is defined to be equal to $f$ on $V_1$ and to the identity map on $V_2$.

We define an axial function, $\alpha^{\#}$, by letting 
$\alpha^{\#}= \alpha'$ on edges of type 1 and $\alpha^{\#}=\alpha$
on edges of type 2. Thus it remains to define $\alpha^{\#}$ on edges 
of type 3. Let $p_{ia}q_{ia}$ such an edge. Then we define
\begin{equation*}
\alpha_{q_{ia}p_{ia}}^{\#} =\alpha_{q_{ia}p_i}
\quad \mbox{and} \quad
\alpha_{p_{ia}q_{ia}}^{\#} =\frac{1}{n_{ia}}\alpha_{p_iq_{ia}}.
\end{equation*}

We define a connection, $\theta^{\#}$, on $\Gamma^{\#}$, 
by letting $\theta^{\#}$ be equal to $\theta'$ on edges of type 1 
and equal to $\theta$ on edges of type 2. 
Thus it remains to define 
$\theta^{\#}$ along edges of type 3. Let $p_i$ be a vertex of $V_0$ and 
$q_{ia}\in V_2$ an adjacent vertex. Then there is a holonomy map
\begin{equation}
\label{eq:2.1}
\theta_{p_iq_{ia}} : E_{p_i} \to E_{q_{ia}}.
\end{equation}
Moreover, one can identify $E_{p_i}$ with $E_{p_{ia}}$ as follows. 
If $p_j$ is a vertex of $\Gamma_0$ adjacent to $p_i$ then there 
is a unique vertex, $p_{jb}$, sitting over $p_j$ in $\Gamma_1$ 
and adjacent in 
$\Gamma_1$ to $p_{ia}$, by (\ref{eq:1.7}). If $q_{ib}$ 
is a vertex of $\Gamma$ 
adjacent to $p_i$ but not in $\Gamma_0$, then, by definition, 
it corresponds to an element, $p_{ib}$, of $N_{p_i}$. 
Thus we can join it to $p_{ia}$ by an edge 
of type 2, or, if $q=q_{ia}$, by an edge of type 3. Composing the map 
(\ref{eq:2.1}) with this identification of $E_{p_{ia}}$ with $E_{p_i}$,
 we get a holonomy map 
\begin{equation*}
\theta_{p_{ia}q_{ia}}^{\#} : E_{p_{ia}} \to E_{q_{ia}}.
\end{equation*}

Then $(\Gamma^{\#}, \alpha^{\#})$, is a $d$-valent one-skeleton, 
called \emph{the blow-up of $\Gamma$ along $\Gamma_0$}. There 
exists a  \emph{blowing-down} map $\beta : \Gamma^{\#} \to \Gamma$, 
obtained by collapsing all $p_{ia}$'s to the corresponding $p_i$. 
The pre-image of $\Gamma_0$ under $\beta$, called 
\emph{the singular locus of} $\beta$, is a $(d-1)$-valent sub-skeleton of
$\Gamma^{\#}$. A particularly important case occurs when $\Gamma_0$ has
trivial normal holonomy in $\Gamma$. In this case the singular locus 
is naturally isomorphic to the direct product of $\Gamma_0$ with 
$\Gamma'$, which is a complete one-skeleton in $s=d-d_0$ vertices.

Define $\tau : V_{\Gamma^{\#}} \to \SS^1(\fg^*)$ by 
\begin{equation}
\label{eq:thomclass}
\tau(v)= \left\{ \begin{array}{ccc}
0 & \mbox{ if } & v  \in V_{\Gamma^{\#}} - V_{\Gamma_0^{\#}} \\
\frac{\textstyle 1}{\textstyle n_{ia}}
\alpha_{p_iq_{ia}} & \mbox{ if } & v=p_{ia} 
\end{array} \right. 
\end{equation}
Then $\tau \in H^{2}(\Gamma^{\#}, \alpha^{\#})$ and will be called
\emph{the Thom class of $\Gamma_0^{\#}$ in $\Gamma^{\#}$}.

\subsubsection{\bf The cohomology of the singular locus}
\label{ssec:cohsingloc}

Let $\Gamma$ be a one-skeleton, $\Gamma_0$ a  sub-skeleton
of covalence $s$, $\Gamma^{\#}$ the blow-up of $\Gamma$ along $\Gamma_0$
and $\Gamma_0^{\#}$ the singular locus of $\Gamma^{\#}$, as defined in the 
preceding section. Let 
\begin{equation}
\label{eq:4.10}
\beta : \Gamma_0^{\#} \to \Gamma_0
\end{equation}
be the blowing-down map. One has an inclusion, 
$H(\Gamma) \to H(\Gamma^{\#})$, and an element $f \in H(\Gamma^{\#})$ 
is the image of an element of
$H(\Gamma)$ if and only if it is constant on the fibers of $\beta$.
Let $\tau \in H^2(\Gamma^{\#})$ be the Thom class of $\Gamma_0^{\#}$ in
$\Gamma^{\#}$ and $\tau_0 \in H^2(\Gamma_0^{\#})$ be the restriction 
of $\tau$ to $\Gamma_0^{\#}$.

\begin{lemma}
\label{lem:cohsingloc}
Every element $f \in H^{2m}(\Gamma_0^{\#})$ can be written uniquely
as
\begin{equation}
\label{eq:4.11}
f= \sum_{k=0}^{s-1} \tau_0^kf_{m-k},
\end{equation}
with $f_{m-k} \in H^{2(m-k)}(\Gamma_0)$ if $k \leq m$ and 0 otherwise.
\end{lemma}

\begin{proof}

For $ p \in V_0$ let
$N_p=\beta^{-1}(p)$ be the fiber over $p$. 
By definition, $N_p$ is a complete one-skeleton
with $s$ vertices  for which a generating class is the restriction of $\tau_0$
to $N_p$. If $f \in H^{2m}(\Gamma_0^{\#})$ then $h$, 
the restriction of $f$ to
$N_p$, is an element of $H^{2m}(N_p)$ and hence, by (\ref{eq:4.5}),
\begin{equation*}
h = \sum_{k=0}^{s-1}f_{m-k}(p) \tau_0^k,
\end{equation*}
where $f_{m-k}(p) \in \SS^{m-k}(\fg^*)$ if $k \leq m$ and is 0 otherwise.
To get (\ref{eq:4.11}) we need to show that the maps 
$f_k : V_0 \to \SS^k(\fg^*)$ are in $H^{2k}(\Gamma_0)$.

Let $p_i,p_j \in V_0$ be joined by an edge. If $q_{ia}, a=d_0+1,..,d$ 
are the neighbors of $p_i$ not in $\Gamma_0$, 
the connection along the edge 
$p_ip_j$ transforms the edges, $p_iq_{ia}$, into edges, $p_jq_{ja}$, 
$a=d_0+1,..,d$, modulo some relabeling.

Then $p_{ia}$ and $p_{ja}$ are joined by an edge in $\Gamma_0^{\#}$, 
which implies that $\alpha_{p_ip_j}$ divides $f(p_{ia})-f(p_{ja})$, and hence
that 
\begin{equation*}
f(p_{ia}) -f(p_{ja}) \equiv 0 \mod{(\tau_0(p_{ia}) - \tau_0(p_{ja}))}.
\end{equation*}
From
\begin{eqnarray}
f(p_{ia})-f(p_{ja}) & = &  \sum_{k=0}^{s-1} 
f_{m-k}(p_j)(\tau^k(p_{ia}) - \tau^k(p_{ja}))  \nonumber \\
& + &  \sum_{k=0}^{s-1}
(f_{m-k}(p_i)-f_{m-k}(p_j))\tau^k(p_{ia}) \nonumber
\end{eqnarray}
we deduce that for every $a=d_0+1,..,d$
\begin{equation}
\label{eq:4.12}
\sum_{k=0}^{s-1} (f_{m-k}(p_i)-f_{m-k}(p_j))\tau^k(p_{ia}) \equiv 0
\pmod{\alpha_{p_ip_j}}.
\end{equation}
Since $\tau(p_{ia}) - \tau(p_{ib})$ is not a multiple of $\alpha_{p_ip_j}$
for $a \neq b$ (recall that $(\Gamma,\alpha)$ 
is assumed to be 3-independent; see the comment at the end of 
Section \ref{sec:skeletons}), 
the relations (\ref{eq:4.12}) imply that every term
$f_{m-k}(p_i)-f_{m-k}(p_j)$ is a multiple of $\alpha_{p_ip_j}$, 
which means that $f_{m-k} \in H^{2(m-k)}(\Gamma_0)$ if $k \leq m$ or is 0 
otherwise.
\end{proof}

\subsubsection{\bf The cohomology of the blow-up}
\label{ssec:cohblow}

We can now determine the additive structure of the cohomology ring of the 
blown-up one-skeleton. The following identity is a graph theoretic version of 
the exact sequence (\ref{eq:shexseq}).

\begin{theorem}
\label{th:cohblow}
\begin{equation*}
H^{2m}(\Gamma^{\#}) \simeq H^{2m}(\Gamma) \oplus 
\bigoplus_{k=1}^{s-1} H^{2(m-k)}(\Gamma_0).
\end{equation*}
\end{theorem}

\begin{proof}
We will show that
every element $f \in H^{2m}(\Gamma^{\#})$ can be written uniquely as
\begin{equation}
\label{eq:4.13}
f= g + \sum_{k=1}^{s-1} \tau^k f_{m-k},
\end{equation}
with $g \in H^{2m}(\Gamma)$, $f_{m-k} \in H^{2(m-k)}(\Gamma_0)$, if
$1 \leq k \leq m$ and 0, if $k > m$.

The restriction, $h$,  of $f$ to $\Gamma_0^{\#}$,  is an element of 
$H^{2m}(\Gamma_0^{\#})$, and, therefore, from Lemma 
\ref{lem:cohsingloc} it follows that
\begin{equation*}
h =f_m + \sum_{k=1}^{s-1} \tau_0^{k}f_{m-k}.
\end{equation*}
But then
\begin{equation*}
g=f-\sum_{k=1}^{s-1} \tau_0^{k}f_{m-k}
\end{equation*}
is constant along fibers of $\beta$, implying that $g \in H^{2m}(\Gamma)$.
Hence $f$ can be written as in (\ref{eq:4.13}). If 
\begin{equation*}
f= g' + \sum_{k=1}^{s-1} \tau^k f_{m-k}'
\end{equation*}
is another decomposition of $f$, then $g-g'$ is supported on
$\Gamma_0^{\#}$ and, therefore, 
\begin{equation*}
0= g-g' + \sum_{k=1}^{s-1} \tau_0^k (f_{m-k} -f_{m-k}'),
\end{equation*}
which, from the uniqueness of (\ref{eq:4.11}), implies that $g=g'$
and $f_{m-k}=f_{m-k}'$ for all $k$'s.
\end{proof}

\subsection{\bf Reduction}
\label{sec:reduction}

\subsubsection{\bf The reduced one-skeleton}
\label{ssec:redskel}

Let $(\Gamma,\alpha)$ be a $d$-valent 
non-cyclic (in the sense of Definition  \ref{def:nca})
one-skeleton and let 
$\phi : V_{\Gamma} \to \RR$ be an injective 
function which is $\xi$-compatible for some polarizing vector $\xi$. 
The image of $\phi$ will be called the set of \emph{critical values} 
of $\phi$ and its complement in $\RR$ the set of 
\emph{regular values}.

For each regular value, $c$, we will construct a new 
$(d-1)$-valent one-skeleton 
$(\Gamma_c, \alpha^c)$. This new one-skeleton will 
be called \emph{the reduced one-skeleton of $(\Gamma, \alpha)$ at} $c$. 
The construction we are about to describe is motivated by the geometric
description of reduction in Theorem \ref{th:reduced}. (In the remaining of 
this section we will use the notation $(p,q)$ for an \emph{unoriented} edge
joining $p$ and $q$, and the notation $pq$ for an \emph{oriented} edge 
with initial vertex $p$ and terminal vertex $q$.)

Consider the \emph{cross-section} of $\Gamma$ at $c$, consisting 
of all edges $(p_0,p_i)$ of $\Gamma$ such that 
$\phi(p_i) < c < \phi(p_0)$; to each such edge 
we associate a vertex $v_i$ of a new graph,
denoted by $\Gamma_c$. Let $r$ be the index of $p_0$ and 
$s=d-r$.

\begin{figure}[h]
\begin{center}
\includegraphics{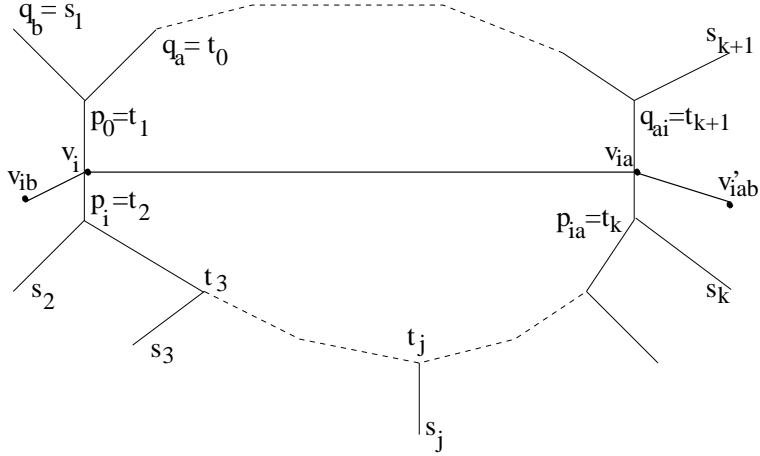}
\caption{Reduction}
\label{fig:red}
\end{center}
\end{figure}

Let the other $d-1$ oriented edges issuing from $p_0$ be denoted by
$p_0q_a$, $a =1,..,d$, $a \neq i$, and let 
$\fh_a$ be the annihilator in $\fg$
of the 2-dimensional linear subspace of $\fg^*$ generated by
$\alpha_{p_0p_i}$ and $\alpha_{p_0q_a}$. The connected
component of $\Gamma_{\fh_a}$ that contains $p_0,p_i$ and $q_a$
has Betti number equal to 1; therefore there exists exactly one more 
edge $e_{ia}=(p_{ia},q_{ai})$ in this component that crosses 
the $c$-level of $\phi$, 
that is, for which $\phi(p_{ia}) < c < \phi(q_{ai})$. If $v_{ia}$ is the 
vertex of $\Gamma_c$ corresponding to $e_{ia}$, we add an edge connecting
$v_i$ to $v_{ia}$. The axial function on the oriented edge $v_iv_{ia}$ is
\begin{equation}
\label{eq:3.2}
\alpha_{v_iv_{ia}}^c = \alpha_{p_0q_a} -
\frac{\alpha_{p_0q_a}(\xi)}{\alpha_{p_0p_i}(\xi)} \alpha_{p_0p_i}
\end{equation}
This axial function takes values in $\fg_{\xi}^*$, 
the annihilator of $\xi$ in $\fg^*$.

The reduced one-skeleton has two connections: an ``up'' connection and a 
``down'' connection. Along the edge $v_iv_{ia}$,, the ``down'' connection
of the reduced one-skeleton is defined as follows. 
Let $v_iv_{ib}$ be another edge at $v_i$, corresponding to an
edge $p_0q_b$. Let $s_1 = q_b$ and let 
$t_2s_2,..., t_{k+1}s_{k+1}$ be the edges obtained by transporting
$t_1s_1$ along the path $t_1, t_2,...,t_{k+1}$. The edge 
$t_{k+1}s_{k+1}$ 
corresponds to a neighbor, $v_{iab}'$, of $v_{ia}$ and we will define the
``down'' connection on the edge $v_i v_{ia}$ by requiring that it sends
$v_iv_{ib}$ to  $v_{ia}v_{iab}'$.
The ``up'' connection is defined similarly except that instead of transporting 
$t_1s_1$ along the bottom path in Figure \ref{fig:red}, we transport 
it along the top path from $t_1$ to $t_{k+1}$.

\begin{theorem}
\label{th:4.1}
The reduced one-skeleton at $c$ is a $(d-1)$-valent one-skeleton. If
$(\Gamma,\alpha)$ is $l$-independent then $(\Gamma_c, \alpha^{c})$ is 
$(l-1)$-independent.
\end{theorem}

\begin{proof}

Let $v_i$ be a vertex of $\Gamma_c$, corresponding to the edge
$e=(p_0,p_i)$ of $\Gamma$ and let $v_{ia}$ be a neighbor of $v_i$,
corresponding to the edge $e_{ia}=(p_{ia},q_{ai})$ and obtained as above 
by using the edge $e_a=p_0q_a$. 
Let $t_0 = q_a, t_1 = p_0,t_2=p_i, ..., t_k = p_{ia}, t_{k+1}= q_{ai}$ 
be the path that connects $q_0$ and 
$q_{ai}$, crosses the $c$-level and is contained in the 2-dimensional
sub-skeleton of $\Gamma$ generated by $p_0, p_i$ and $q_a$
(see Figure \ref{fig:red}).
 
It is clear that $\alpha_{v_iv_{ia}}^c$ is a positive multiple of 
\begin{equation*}
\iota_{\xi} ( \alpha_{t_1t_0} \wedge \alpha_{t_1t_2}),
\end{equation*}
where $\iota_{\xi}$ is the interior product with $\xi$. 
Axiom \eqref{eq:axiom3} implies that, for every $j$, 
\begin{equation*}
\alpha_{t_jt_{j-1}} \wedge \alpha_{t_jt_{j+1}} \mbox{ and }
\alpha_{t_{j-1}t_{j-2}} \wedge \alpha_{t_{j-1}t_j}
\end{equation*}
are positive multiples of each others and, hence, that
\begin{equation*}
\iota_{\xi}( \alpha_{t_jt_{j-1}} \wedge \alpha_{t_jt_{j+1}}) 
\mbox{ is a positive multiple of }
\iota_{\xi}( \alpha_{t_{j-1}t_{j-2}} \wedge \alpha_{t_{j-1}t_{j}}).
\end{equation*}
Therefore $\alpha_{v_iv_{ia}}^c$ is a negative multiple of 
$\alpha_{v_{ia}v_i}^c$, so $\alpha^c$ satisfies axiom A\ref{axiom2}
of Definition \ref{def:abs1skel}.

We will show that $(\Gamma_c,\alpha^c)$ satisfies 
axiom A\ref{axiom3} of Definition \ref{def:abs1skel}. Note that 
\begin{equation*}
\alpha_{v_iv_{ib}}^c =
\frac{\iota_{\xi} ( \alpha_{t_1t_2} \wedge \alpha_{t_1s_1})} 
{\alpha_{t_1t_2}(\xi)}
\quad \mbox{ and that } \quad
\alpha_{v_{ia}v_{iab}'}^c =
\frac{\iota_{\xi} ( \alpha_{t_{k+1}t_k} \wedge \alpha_{t_{k+1}s_{k+1}})}
{\alpha_{t_{k+1}t_k}(\xi)}.
\end{equation*}
A direct computation shows that
\begin{equation*}
\frac{\iota_{\xi} ( \alpha_{t_jt_{j+1}} \wedge \alpha_{t_js_j})}
{\alpha_{t_jt_{j+1}}(\xi)} -
\frac{\iota_{\xi} ( \alpha_{t_jt_{j-1}} \wedge \alpha_{t_js_j})}
{\alpha_{t_jt_{j-1}}(\xi)} 
\end{equation*}
is a multiple of $\alpha_{v_iv_{ia}}^c$; if
$\alpha_{t_{j+1}s_{j+1}} = \lambda_j \alpha_{t_js_j} + 
c_j \alpha_{t_jt_{j+1}}$ 
 then 
\begin{equation*}
\frac{\iota_{\xi} ( \alpha_{t_{j+1}t_j} \wedge \alpha_{t_{j+1}s_{j+1}})}
{\alpha_{t_{j+1}t_j}(\xi)} =
\lambda_j
\frac{\iota_{\xi} ( \alpha_{t_jt_{j+1}} \wedge \alpha_{t_js_j})}
{\alpha_{t_jt_{j+1}}(\xi)},
\end{equation*}
and eliminating the intermediary terms, we see that
\begin{equation}
\label{eq:3.4}
\alpha_{v_{ia}v_{iab}'}^c - \lambda \alpha_{v_iv_{ia}}^c 
\mbox{ is a multiple of }
\alpha_{v_iv_{ia}}^c,
\end{equation}
with $\lambda=\lambda_k\cdots \lambda_1> 0$.

Axiom A\ref{axiom1} of Definition \ref{def:abs1skel}, 
as well as the statement about the 
$(l-1)$- independence, follows at once from the fact that 
the value of $\alpha^c$ at the vertex $v_i$ of $\Gamma^c$ is a linear 
combination of {\it two} values of $\alpha$ at a vertex $p_0$ of
$\Gamma$, one of which is fixed, namely $\alpha_{p_0p_i}$. 
\end{proof}

\subsubsection{\bf Passage over a critical value}
\label{ssec:passage}

In this section we will describe what happens to the 
reduced one-skeleton at $c$ as $c$ varies; it is clear that if 
there is no critical value between $c$ and $c'$
then the two reduced one-skeleta are identical. 
Suppose, therefore, that there exists exactly one critical value in the 
interval $(c,c')$, and that it is attained at the vertex $p_0$. 
Let $r$ be the index of $p_0$ and $s=d-r$.

\begin{theorem}
\label{th:4.2}
$(\Gamma_{c'}, \alpha_{c'})$ is obtained from 
$(\Gamma_c, \alpha_c)$ by a blowing-up of $\Gamma_c$ along a
complete sub-skeleton with $r$ vertices followed by a 
blowing-down along a complete sub-skeleton with $s$ vertices.
\end{theorem}

\begin{proof}
 
The modifications from $\Gamma_c$ to $\Gamma_{c'}$ 
are due to the edges that cross one level but not the other one; 
but these are exactly the edges with one end-point $p_0$. 
Let $p_0p_i, \; i=1,...,r$ the edges of $\Gamma$ with initial vertex $p_0$ 
that point downward and $p_0q_a ,\; a=r+1,...,d$ the edges that
point upward. Let $v_i$ be the vertex of $\Gamma_c$ associated
to $(p_0,p_i)$ and $w_a$ the vertex of $\Gamma_{c'}$ associated
to $(p_0,q_a)$, for $i=1,..,r$ and $a=r+1,..,d$. 
The $v_i$'s are the vertices of a complete sub-skeleton 
$\Gamma_c'$ of $\Gamma_c$ and the $w_a$'s are the vertices of 
a complete sub-skeleton $\Gamma_{c'}'$ of $\Gamma_{c'}$, 
$\Gamma_c'$ having trivial normal holonomy with respect to the 
``up'' connection on $\Gamma_c$ and $\Gamma_{c'}'$ having trivial
normal holonomy with respect to the ``down'' connection on $\Gamma_{c'}$.
(This can be shown as follows: 
The \emph{``normal bundle''} to $v_i$ is the same for all $i'$s, namely
it can be identified with the set of edges $p_0q_a$. Moreover,
the holonomy map associated with $v_iv_j$ is by definition 
just the identity map on this set of edges.)

\begin{figure}[h]
\begin{center}
\includegraphics{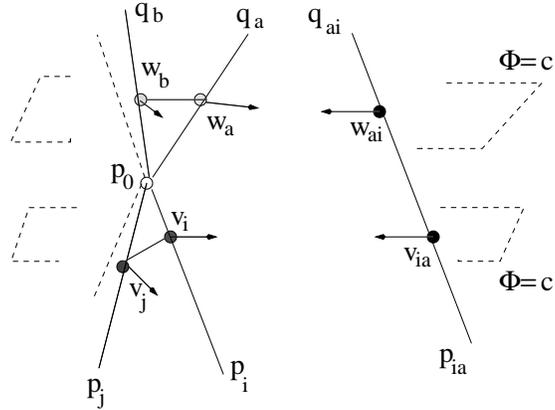}
\caption{Evolution}
\end{center}
\end{figure}

Consider $\fh_{ia}$, the annihilator of the 2-dimensional subspace
of $\fg^*$ generated by $\alpha_{p_0p_i}$ and $\alpha_{p_0q_a}$;
the connected component of $\Gamma_{\fh_{ia}}$ that contains 
$p_0, p_i$ and $q_a$ will contain exactly one edge 
$e_{ia}=(p_{ia},q_{ia})$ that crosses both the $c$-level and the $c'$-level.
To this edge will correspond a neighbor $v_{ia}$ of $v_i$ in 
$\Gamma_c$ and
a neighbor $w_{ai}$ of $w_a$ in $\Gamma_{c'}$.

Let $\mu>0$ and for all $i=1,..,r$, $a=r+1,..,d$, define 
\begin{eqnarray}
\label{eq:nias}
n_{ia} & = & \mu \alpha_{p_0q_a}(\xi) > 0 \\
\label{eq:nais}
n_{ai} & = & -\mu \alpha_{p_0p_i}(\xi) > 0 .
\end{eqnarray}
Denote
\begin{equation*}
\tau_i  = 
 - \frac{\alpha_{p_0p_i}}{\mu\alpha_{p_0p_i}(\xi)} 
=
\frac{\alpha_{p_0p_i} }{n_{ai}} \quad \mbox{ and } \quad 
\tau_a  = 
  \frac{\alpha_{p_0q_a}}{\mu\alpha_{p_0q_a}(\xi)} =
\frac{\alpha_{p_0q_a}}{n_{ia}} \; .
\end{equation*}
Then we have
\begin{eqnarray}
\alpha_{v_iv_j}^c  & = &  \alpha_{p_0p_j} -
\frac{\alpha_{p_0p_j}(\xi)}{\alpha_{p_0p_i}(\xi)} \alpha_{p_0p_i} =
n_{aj}(\tau_j-\tau_i)  \nonumber \\
\alpha_{w_aw_b}^{c'} &  = & \alpha_{p_0q_b} -
\frac{\alpha_{p_0q_b}(\xi)}{\alpha_{p_0q_a}(\xi)} \alpha_{p_0q_a} =
n_{ib}(\tau_b-\tau_a) \nonumber \\
\alpha_{v_iv_{ia}}^c  & = & \alpha_{p_0q_a}-
\frac{\alpha_{p_0q_a}(\xi)}{\alpha_{p_0p_i}(\xi)}\alpha_{p_0p_i} =
n_{ia}(\tau_a + \tau_i) \nonumber \\
\alpha_{w_aw_{ai}}^{c'}  & = &  \alpha_{p_0p_i}-
\frac{\alpha_{p_0p_i}(\xi)}{\alpha_{p_0q_a}(\xi)}\alpha_{p_0q_a} =
n_{ai} (\tau_i + \tau_a) \nonumber .
\end{eqnarray}

Let $\Gamma_c^0$ be the sub-skeleton of $\Gamma_c$ with vertices 
$\{ v_1,...,v_r \}$. Along the edge $v_iv_j$, the connection 
transports the edge $v_iv_{ia}$ to $v_jv_{ja}$. Noting that
$n_{ia}=n_{ja}$ and that
\begin{equation*}
\alpha_{v_jv_{ja}}^c - \alpha_{v_iv_{ia}}^c = \frac{n_{ia}}{n_{ja}} 
\alpha_{v_iv_j}^c,
\end{equation*}
{\em i.e.} that $\Gamma_c$ satisfies \eqref{eq:condition2} 
on edges of $\Gamma_c$ normal to $\Gamma_c^0$, we can define the blow-up 
$\Gamma_c^{\#}$ of
$\Gamma_c$ along $\Gamma_c^0$, by means of the positive numbers 
$n_{ia}$, $i=1,..,r$ and $a=r+1,..,d$.

The singular locus, $\Gamma_0^{\#}$, is, as a graph, a product of two
complete graphs, $\Gamma_r \times \Gamma_s$; each vertex $z_{ia}$
corresponds to a pair $(v_i,w_a)$  
and the blow-down map $\beta : \Gamma_0^{\#} \to \Gamma_0^c$ sends $z_{ia}$
to $v_i$.

There are edges connecting $z_{ia}$ with $v_{ia}$, $z_{ia}$ with $z_{ib}$
and $z_{ia}$ with $z_{ja}$, for all distinct $i,j =1,..,r$ and 
$a,b = r+1,..,d$. The values of the axial function $\alpha^{\#}$ on these
edges are 
\begin{eqnarray}
\alpha_{z_{ia}v_{ia}}^{\#} & = & \frac{1}{n_{ia}}
\alpha_{v_iv_{ia}}^c = \tau_a + \tau_i = 
-\frac{\mu}{n_{ia}n_{ai}} 
\iota_{\xi}(\alpha_{p_0p_i} \wedge \alpha_{p_0q_a})
\nonumber \\
\alpha_{z_{ia}z_{ib}}^{\#} & = & \alpha_{v_iv_{ib}}^c - 
\frac{n_{ib}}{n_{ia}}\alpha_{v_iv_{ia}}^c = 
n_{ib}(\tau_b - \tau_a) = \nonumber \\
& = &\frac{\mu}{n_{ai}}
\iota_{\xi}( \alpha_{p_0q_a} \wedge \alpha_{p_0q_b}) \nonumber \\
\alpha_{z_{ia}z_{ja}}^{\#} & = &  \alpha_{v_iv_j}^c =
n_{aj}(\tau_j - \tau_i) = -\frac{\mu}{n_{ai}}
\iota_{\xi}( \alpha_{p_0p_i} \wedge \alpha_{p_0p_j}) \nonumber
\end{eqnarray}
While $\alpha_{z_{ia}v_{ia}}^{\#}$ and $\alpha_{z_{ia}z_{ib}}^{\#}$ are
not collinear (since $\tau_a, \tau_b$ and $\tau_i$ are independent), and
neither are $\alpha_{z_{ia}v_{ia}}^{\#}$ and 
$\alpha_{z_{ia}z_{ja}}^{\#}$,  it may happen that 
$\alpha_{z_{ia}z_{ib}}^{\#}$ and $\alpha_{z_{ia}z_{ja}}^{\#}$ are collinear.
We can, however, circumvent this problem by means of the following lemma
(which we will leave as an easy exercise).

\begin{lemma}
Let $\omega_1, \omega_2, \omega_3,\omega_4 \in \fg^*$ be 3-independent
and suppose that for some $\xi \in \fg$,
$\iota_{\xi}(\omega_1 \wedge \omega_2)$ and 
$\iota_{\xi}(\omega_3 \wedge \omega_4)$ are collinear. Then the 2-planes
generated by $\{\omega_1, \omega_2\}$ and $\{\omega_3, \omega_4\}$
intersect in a line. Moreover, if $\{ \omega_0 \}$ is a basis for this 
line then $\omega_0(\xi)=0$.
\end{lemma}

Therefore $\alpha_{z_{ia}z_{ib}}^{\#}$ and $\alpha_{z_{ia}z_{ja}}^{\#}$ 
are collinear precisely when $\xi$ belongs to a pre-determined hyperplane; 
so we can insure 2-independence for $\Gamma_c^{\#}$, by avoiding a 
finite number of such hyperplanes.

\begin{definition}
\label{def:generic}
An element $\xi \in \fg$ is called \emph{generic} 
for the one-skeleton, $(\Gamma,\alpha)$, if
for every vertex, $p$, and every quadruple of distinct 
edges $e_1, e_2, e_3$ and $e_4$ in $E_p$, 
the vectors 
$\iota_{\xi}(\alpha_{e_1} \wedge \alpha_{e_2})$ and   
$\iota_{\xi}(\alpha_{e_3} \wedge \alpha_{e_4})$
are linearly independent.
\end{definition}

If $\Gamma$ has valence 3 then every element is generic. 
In general, for every element, $\xi$, of $\P$ and every neighborhood of 
$\xi$ in $\P$, there exists a generic element, $\xi'$, in that neighborhood, 
such that the orientations $o_{\xi}$ and $o_{\xi'}$ are the same, and 
the reduced skeleta corresponding to $\xi$ and $\xi'$ have isomorphic 
underlying graphs.

We now return to the proof of Theorem \ref{th:4.2}. For a generic 
$\xi \in \fg$, the blow-up $\Gamma_c^{\#}$ can still be defined.
Note that 
\begin{equation*}
\frac{1}{n_{ia}}\alpha_{v_iv_{ia}}^c =
\frac{1}{n_{ai}}\alpha_{w_aw_{ai}}^{c'} \quad \mbox{ and that }
\quad \alpha_{v_iv_{ib}}^c - 
\frac{n_{ib}}{n_{ia}}\alpha_{v_iv_{ia}}^c
=\alpha_{w_aw_b}^{c'}.
\end{equation*}
These relations imply that $\Gamma_c^{\#}$ is the same as 
the blow-up, $\Gamma_{c'}^{\#}$, of $\Gamma_{c'}$ along $\Gamma_{c'}^0$ 
using the $n_{ai}$'s. Therefore for generic $\xi$, the 
passage from $\Gamma_c$ to $\Gamma_{c'}$ is equivalent to a blow-up
from $\Gamma_c$ to $ \Gamma^{\#}= \Gamma_c^{\#}$ 
followed by a blow-down from 
$ \Gamma^{\#} = \Gamma_{c'}^{\#}$ to $\Gamma_{c'}$.
\end{proof}

\subsubsection{\bf The changes in cohomology}
\label{ssec:redchanges}

We will now describe how the cohomology changes as one passes from
$\Gamma_c$ to $\Gamma_{c'}$. For this we will use a variant of
Theorem \ref{th:cohblow}. This theorem itself can't be applied directly 
since the reduced one-skeleta might not be 3-independent. 
However, the proof of Lemma \ref{lem:cohsingloc} is valid up to
assertion \eqref{eq:4.12} without this assumption and beyond this 
point it suffices to assume that $\xi$ is generic.
Combining Theorems \ref{th:cohblow} and 
\ref{th:cohcompl} we conclude
\begin{eqnarray}
\label{eq:4.14}
\dim H^{2m}(\Gamma^{\#}) & = & \dim H^{2m}(\Gamma_c) +
\sum_{k=1}^{s-1} \dim H^{2(m-k)} ( \Gamma_r) = \nonumber \\
 & = &  \dim H^{2m}(\Gamma_c) + 
\sum_{k=1}^{s-1} \sum_{l=0}^{r-1} \lambda_{m-k-l} = \nonumber \\
 & = &  \dim H^{2m}(\Gamma_c) + 
\sum_{k=1}^{s-1} \sum_{l=1}^{r-1} \lambda_{m-k-l} +
\sum_{k=1}^{s-1} \lambda_{m-k} . \nonumber
\end{eqnarray}
Therefore 
\begin{equation*}
\dim H^{2m}(\Gamma^{\#})  =   \dim H^{2m}(\Gamma_c) + 
\sum_{k=1}^{s-1} \sum_{l=1}^{r-1} \lambda_{m-k-l} +
\sum_{k=1}^{s-1} \lambda_{m-k}
\end{equation*}
and 
\begin{equation*}
\dim H^{2m}(\Gamma^{\#})  =   \dim H^{2m}(\Gamma_{c'}) + 
\sum_{l=1}^{r-1} \sum_{k=1}^{s-1} \lambda_{m-k-l} +
\sum_{l=1}^{r-1} \lambda_{m-l} , 
\end{equation*}
which imply that
\begin{equation}
\label{eq:6.1}
\dim H^{2m}(\Gamma_{c'}) = \dim H^{2m}(\Gamma_c) +
\sum_{k=1}^{s-1} \lambda_{m-k} - \sum_{k=1}^{r-1} \lambda_{m-k}.
\end{equation}
(N.B. All $\lambda_j$'s in the displays above are $\lambda_{j,n-1}$'s,
since the dimension of $\fg_{\xi}^*$ is $n-1$.) Since
\begin{equation}
\label{eq:lambdann}
\lambda_{a,n} = \sum_{j=0}^{a} \lambda_{j,n-1},
\end{equation}
the equality (\ref{eq:6.1}) can be written as
\begin{equation}
\label{eq:6.1new}
\dim H^{2m}(\Gamma_{c'}) = \dim H^{2m}(\Gamma_c) +
\lambda_{m-r,n}- \lambda_{m-s,n}.
\end{equation}

\subsection{\bf The additive structure of $H(\Gamma,\alpha)$}
\label{sec:addstructure}

\subsubsection{\bf Symplectic cutting}
\label{ssec:cutting}

In this section we will use the results above to draw some conclusions about 
$H(\Gamma,\alpha)$ itself. This we will do by mimicking, in our graph
theoretic setting, the 
\emph{symplectic cutting} operation of E. Lerman (\cite{Le}).

Let $L$ be the ``edge'' graph, with two vertices labeled 0 and 1 and one 
edge connecting them. Consider $L^+ =(L, \alpha^{+})$, with the 
axial function given by
$\alpha_{01}^+ = {\bf 1},\alpha_{10}^+ = -{\bf 1}$ and 
$L^- = (L, \alpha^{-})$, with 
the axial function $\alpha_{01}^- = -{\bf 1},\alpha_{10}^- = {\bf 1}$. Here
$\alpha^{\pm} : E_L \to \RR^* \simeq \RR$. For $\RR$ we 
have the basis $\{ 1 \}$ and for its dual $\RR^*$ the basis 
$\{ {\bf 1} \}$. For both these axial functions, $ 1 \in \RR$ is
polarizing. Finally, let $\phi_0^{\pm} : V_L \to \RR$
be given by $ \phi_0^{\pm}(0)=0, \phi_0^{\pm}(1) = \pm 1$. Then
$\phi_0^{\pm}$ is $o_1$-compatible for $\alpha^{\pm}$.

Let $(\Gamma,\alpha)$ be a one-skeleton which is 3-independent and 
non-cyclic in the sense of Definition \ref{def:nca}. Let 
$\phi : V_{\Gamma} \to \RR$ be $\xi$-compatible for 
some $\xi \in \P$ and choose  
$a > \phi_{max} - \phi_{min} > 0$.
For $c \in  (\phi_{min}, \phi_{max})$
let $(\Gamma_c,\alpha^c)$ be the
reduced one-skeleton of $\Gamma$ at $c$.

Consider the product one-skeleton 
$(\Gamma \times L^+, \alpha \times \alpha^+)$, with 
$$\alpha \times \alpha^+ : E_{\Gamma \times L} \to 
\fg^* \oplus \RR^* \simeq (\fg \oplus \RR)^*.$$
This one-skeleton is also 3-independent and non-cyclic, and the function 
$$\Phi^{+}(p,t) = \phi(p)+a\phi_0^{+}(t)=\phi(p)+at$$
is $(\xi,1)$-compatible.

Define
$\Gamma_{\phi \leq c}$ to be the reduced one-skeleton of 
$(\Gamma \times L^+, \alpha \times \alpha^{+})$ 
at $\Phi^{+} = c$. The vertices of $\Gamma_{\phi \leq c}$
correspond to two types of edges of $\Gamma \times L^+$:
\begin{enumerate}
\item $((p_i,0),(p_j,0))$, with 
	$\phi(p_i) < c < \phi(p_j)$
\item $((p_i,0),(p_i,1))$ with $\phi(p_i) < c$
\end{enumerate}
\begin{figure}[h]
\begin{center}
\includegraphics{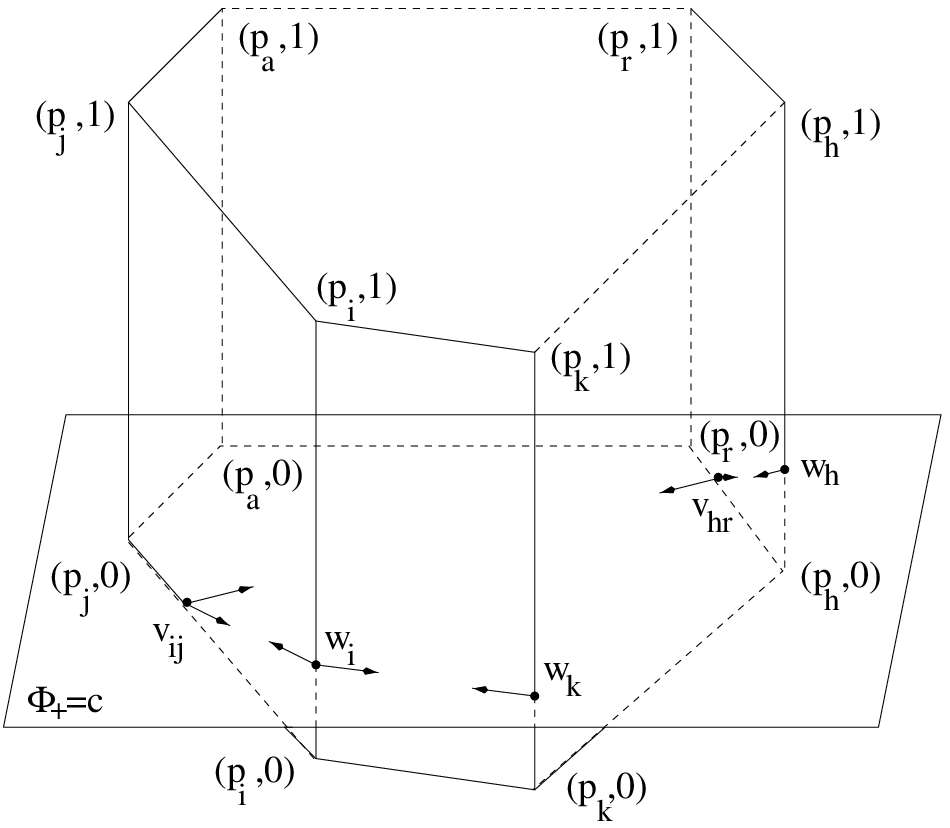}
\caption{Cutting}
\end{center}
\end{figure}
Let $v_{ij}$ be the vertex of $\Gamma_{\phi \leq c}$
corresponding to an edge $((p_i,0),(p_j,0))$ and $w_i$ the  vertex
corresponding to an edge $((p_i,0),(p_i,1))$. 

The neighbors of $w_i$ are of two types:
\begin{enumerate}
\item $w_k$, if $(p_i,p_k)$ is an edge of $\Gamma$ and 
	$\phi(p_k) < c$;
\item $v_{ij}$, if $(p_i,p_j)$ is an edge of $\Gamma$ and 
	$\phi(p_j) > c$.
\end{enumerate}
As for neighbors of $v_{ij}$, apart from $w_i$, they are precisely the
neighbors of $v_{ij}$ in the reduced one-skeleton $\Gamma_c$.
This $\Gamma_c$ sits inside $\Gamma_{\phi \leq c}$ as the subgraph
with vertices $v_{ij}$.

Using (\ref{eq:3.2}) we deduce that the axial function of 
$\Gamma_{\phi \leq c}$,
which we will denote by $\beta^{+}$, is given by: 
\begin{eqnarray}
\beta_{w_iw_k}^{+} &=& \alpha_{p_ip_k} - \alpha_{p_ip_k}(\xi)\cdot {\bf 1}
\nonumber \\
\beta_{v_{ij}w_i}^{+} & = &
- \frac{1}{\alpha_{p_jp_i}(\xi)} \alpha_{p_jp_i} +  {\bf 1} 
\nonumber \\
\beta_{v_{ij}v_{hr}}^{+} & = &  \alpha_{p_jp_a}-
\frac{\alpha_{p_jp_a}(\xi)}{\alpha_{p_jp_i}(\xi)}\alpha_{p_jp_i} =
\alpha_{v_{ij}v_{hr}}^c \nonumber
\end{eqnarray}
The axial function $\beta^{+}$ takes values in 
$(\fg \oplus \RR)_{(\xi,1)}^* \subset \fg^* \oplus \RR^*$.
However, there is a natural isomorphism 
$\fg^* \rightarrow (\fg \oplus \RR)_{(\xi,1)}^* $, given by 
\begin{equation}
\label{eq:sigma}
 \sigma \longrightarrow 
(\sigma, -\sigma(\xi)\cdot {\bf 1}),
\end{equation}
so we can regard $\beta^{+}$ as taking values in $\fg^*$, and, as such,
it is given by:
\begin{eqnarray}
\beta_{w_iw_k}^{+} &=& \alpha_{p_ip_k} \nonumber \\
\beta_{v_{ij}w_i}^{+} & = &
- \frac{1}{\alpha_{p_jp_i}(\xi)} \alpha_{p_jp_i} \nonumber \\
\beta_{v_{ij}v_{hr}}^{+} & = & \alpha_{v_{ij}v_{hr}}^c \nonumber
\end{eqnarray}

Similarly one can define $\Gamma_{\phi \geq c}$ as the reduced one-skeleton 
of $(\Gamma \times L^-, \alpha \times \alpha^-)$ at $\Phi^{-}=c$, where 
$$\Phi^{-}(p,t) = \phi(p)+a\phi_0^{-}(t)=\phi(p)-at$$
is $(\xi,1)$-compatible.
Note that if $\xi$ is generic then $(\xi, 1)$ is generic as well.

\subsubsection{\bf The dimension of $H(\Gamma,\alpha)$}
\label{ssec:dimension}

We will now apply (\ref{eq:6.1}) several times 
to suitable chosen one-skeleta to get the following result, which, 
in some sense, implies the main results of this article:

\begin{theorem}
\label{th:addstruct}
Let $(\Gamma,\alpha)$ be a $d$-valent one-skeleton 
which is 3-independent and non-cyclic. Then 
\begin{equation}
\label{eq:dimhg}
\dim H^{2m}(\Gamma,\alpha) = \sum_{k=0}^d b_{2k}(\Gamma) \lambda_{m-k,n} \; ,
\end{equation}
the $\lambda'$s being defined by (\ref{eq:4.2}).
\end{theorem}

\begin{proof}
Let $\xi \in \fg$ be a generic element of $\P$, 
let $\phi : V_{\Gamma} \to \RR$ be $\xi$-compatible 
and let $\Gamma_{\phi \leq c}$ be the one-skeleton defined 
in the previous section.

If there is only one vertex $(p,t)$ with $\Phi^{+}(p,t) < c_0$, 
this one-skeleton is $\Gamma_{d+1}$, the 
complete one-skeleton with $d+1$ vertices and if 
$$\phi_{max} < c_1 < a + \phi_{min}$$
this one-skeleton is just $(\Gamma, \alpha)$.
Therefore, by studying the change in the cohomology of 
$\Gamma_{\phi \leq c}$ as $c$ varies, we can
determine the additive structure of $H(\Gamma, \alpha)$ from the additive 
structure of $H(\Gamma_{d+1})$.

Let $c_0 < a < b < c_1$ such that there is exactly one vertex 
$p \in V_{\Gamma}$ with $a < \Phi^{+}(p,t) < b$.  If the index
of $p$ in $\Gamma$ is $\sigma(p)=r$ then the index of $(p,0)$ in 
$\Gamma \times L^{+}$ is also $r$. Note that since the zeroth Betti
number of $\Gamma$ is 1, $r$ can't be 0. Also, in this case,
$1 \leq s=d+1-r < d+1$. Thus we can apply 
(\ref{eq:6.1}) to obtain
\begin{equation*}
\dim H^{2m}(\Gamma_{\phi \leq b}) = 
\dim H^{2m}(\Gamma_{\phi \leq a}) +
\sum_{k=1}^{d-r} \lambda_{m-k} - \sum_{k=1}^{r-1} \lambda_{m-k}.
\end{equation*}
Adding together these changes we get
\begin{eqnarray}
\dim H^{2m}(\Gamma, \alpha) & = & \dim H^{2m}(\Gamma_{d+1}) + 
\sum_{\sigma(p) > 0} 
\left( \sum_{k=1}^{d-\sigma(p)}\lambda_{m-k} - 
\sum_{k=1}^{\sigma(p)-1} \lambda_{m-k} \right) \nonumber = \\
& = & \lambda_m + 
\sum_{\sigma(p) \geq 0}
\left( \sum_{k=1}^{d-\sigma(p)}\lambda_{m-k} - 
\sum_{k=1}^{\sigma(p)-1} \lambda_{m-k} \right) \nonumber
\end{eqnarray}

The minimum value for $k$ is 1 and the maximum is $d$; $\lambda_{m-k}$
appears in the first sum when $\sigma(p) \leq d-k$ and in the second one when
$\sigma(p) \geq k+1$. Therefore 
\begin{equation*}
\dim H^{2m}(\Gamma, \alpha)  =  
\lambda_m + \sum_{k=1}^d \left(
\sum_{l=0}^{d-k} b_{2l}(\Gamma) - 
\sum_{l=k+1}^{d} b_{2l}(\Gamma) \right) \lambda_{m-k}. 
\end{equation*}

Because of the relations $b_{2d-2l}(\Gamma) = b_{2l}(\Gamma)$ 
(see \eqref{eq:1.4}) 
the expression in bracket reduces to $b_{2k}(\Gamma)$ and therefore
$$\dim H^{2m}(\Gamma, \alpha) = 
\sum_{k=0}^d b_{2k}(\Gamma) \lambda_{m-k,n} \;.$$
\end{proof}

\subsubsection{\bf Generators for $H(\Gamma,\alpha)$}
\label{ssec:generators}

We can sharpen the result above by constructing a set of generators 
for $H(\Gamma,\alpha)$ with nice support conditions.
Let $\phi : V_{\Gamma} \to \RR$ be $\xi$-compatible for 
$\xi \in \P$. For $p \in V_{\Gamma}$, let $F_p$ be the 
\emph{flow-out of $p$}, 
that is the set of vertices of the oriented graph 
$(\Gamma, o_{\xi})$ that can be reached by a positively oriented path 
starting from $p$.

\begin{theorem}
\label{th:suppthom}
If $p \in V_{\Gamma}$ is a vertex of index $r$, then 
there exists an element $\tau_p \in H^{2r}(\Gamma, \alpha)$, 
with the following properties:
\begin{enumerate}
\item \label{item:i} $\tau_p$ is supported on $F_p$
\item \label{item:ii} $\tau_p(p)= \prod \alpha_{e}$,
the product over edges $e \in E_p$ with $\alpha_{e}(\xi)<0$.
\end{enumerate}
\end{theorem}

\begin{proof}
We first recall a construction used in \cite[sec. 2.9]{GZ}. 
For every regular value $c \in \RR$, 
let $H_c(\Gamma,\alpha)$ be the subring of 
those maps $f \in H(\Gamma,\alpha)$ that are supported on the set 
$\phi \geq c$. Now consider regular values 
$c,c'$ such that there is exactly 
one vertex, $p$, satisfying $c < \phi(p) < c'$. 
Let $r=\sigma(p)$ be the index of $p$ and let 
$\alpha_1,..,\alpha_r$ be the
values of the axial function on the edges pointing down from $p$.
Consider the restriction map
$$H_c^{2m}(\Gamma,\alpha) \to \SS^{m}(\fg^*), \quad f \to f(p).$$
The image of this map is contained in 
$\alpha_1 \cdots \alpha_r \SS^{m-r}(\fg^*)$, 
and the kernel is $H_{c'}^{2m}(\Gamma,\alpha)$,
so we have an exact sequence
\begin{equation}
\label{eq:6.2}
0 \to H_{c'}^{2m}(\Gamma,\alpha) \to H_c^{2m}(\Gamma,\alpha) \to
\alpha_1 \cdots \alpha_r \SS^{m-r}(\fg^*) .
\end{equation}
Therefore 
\begin{equation}
\label{eq:6.3}
\dim H_c^{2m}(\Gamma,\alpha) - \dim H_{c'}^{2m}(\Gamma,\alpha) \leq 
\lambda_{m-r}.
\end{equation}
So when we go from $c < \phi_{min}$ to $c' > \phi_{max}$
and add together the inequalities (\ref{eq:6.3}), we get the inequality
\begin{equation}
\label{eq:6.4}
\dim H^{2m}(\Gamma, \alpha) \leq \sum_{r=0}^d b_{2r}(\Gamma) \lambda_{m-r}.
\end{equation}
But we proved that (\ref{eq:6.4}) is actually an equality, so all the 
inequalities (\ref{eq:6.3}) are equalities, which means that the right 
arrow in (\ref{eq:6.2}) is surjective. In particular, when $m=r(=\sigma(p))$, 
there exists an element $\tau_p' \in H_c^{2 \sigma(p)}(\Gamma, \alpha)$ with
\begin{equation*}
\tau_p'(p)=\alpha_1 \cdots \alpha_r = 
\prod_{e \in E_p, \alpha_{e}(\xi)<0} \alpha_{e},
\end{equation*}
verifying the second condition.
However, to get  $\tau_p'$ to be supported on $F_p$ we will need to 
modify it, and this we will do inductively, as follows:

If $p$ is the vertex where $\phi$ takes its maximum, 
$F_p=\{ p \}$ and the first condition 
is automatically verified. Assume now that we have constructed an element
$\tau_q$ satisfying both conditions for all 
vertices, $q$, above $p$. We will show that a $\tau_p$ exists for $p$ as well. 
Suppose such an element doesn't exist and 
let $p_1$ be the highest possible first vertex 
not in $F_p$ where $\tau_p'$ is non zero, for all choices of $\tau_p'$.
At all the neighbors of $p_1$ below $p_1$, the value of $\tau_p'$ is zero and 
therefore $\tau_p'(p_1)$ can be written as 
\begin{equation*}
\tau_p'(p_1) = g \prod_{e \in E_{p_1}, \alpha_{e}(\xi) < 0} \alpha_{e},
\end{equation*}
for some $g$. If $\sigma(p_1) > \sigma(p) =\deg{\tau_p'}$ then this is 
possible 
if and only if $\tau_p'(p_1)=0$, which contradicts the choice of $p_1$.
So we must have $\sigma(p_1) \leq \sigma(p)$. From the induction hypotheses 
we know
that there exists a $\tau_{p_1}$ with the required properties.
If we replace $\tau_p'$ with $\tau_p''=\tau_p'-g\tau_{p_1}$ then this new
element will still satisfy the second condition, 
will be zero at all vertices 
not in $F_p$ that are below $p_1$ but will also be 0 at $p_1$, which 
contradicts the ``maximality'' of $\tau_p'$. 

Therefore an element $\tau_p$ must exists for $p$ as well. 
\end{proof}

The method of proof above also gives us the following uniqueness result.

\begin{theorem}
\label{th:unictaup}
Suppose that for every point $q \in F_p$, different from $p$, the index of
$q$ is strictly greater than the index of $p$. Then the class, $\tau_p$, is 
unique.
\end{theorem}

We will show that these classes generate $H(\Gamma, \alpha)$ as a module 
over $\SS(\fg^*)$.

\begin{theorem}
\label{th:freemodule}
If $\{\tau_p\}_{p\in V_{\Gamma}}$ satisfy the 
hypotheses of Theorem \ref{th:suppthom} then 
\begin{equation*}
H^{2m}(\Gamma, \alpha) = \bigoplus_{\sigma(p) \leq m} 
\SS^{m-\sigma(p)}(\fg^*) \tau_p.
\end{equation*}
In particular, $H(\Gamma,\alpha)$ is a free $\SS(\fg^*)$-module with the 
$\tau_p$'s as generators.
\end{theorem}

\begin{proof}
We will show that
every element $f \in H^{2m}(\Gamma, \alpha)$
can be written uniquely as
\begin{equation}
\label{eq:6.5}
f = \sum_{\sigma(p) \leq m} h_p\tau_p,
\end{equation}
where $h_p \in \SS^{m-\sigma(p)}(\fg^*)$.

Let $p_0,p_1,..,p_N$ be the vertices of $\Gamma$, ordered so that
$$\phi(p_0) < \phi(p_1) < ... < \phi(p_N).$$
Then $\tau_{p_0}(p_0) = 1$ and, if we let
$h_{p_0} = f(p_0)$,
$$f_0=f-h_{p_0}\tau_{p_0}, $$
vanishes at $p_0$. Now, suppose 
$$f_k = f - \sum_{\stackrel{\scriptstyle i < k}
{\sigma(p_i) \leq m}} h_{p_i} \tau_{p_i} 
\in H^{2m}(\Gamma,\alpha),$$
is supported on $\{ p_{k},..,p_N \}$. Then 
$$f_k(p_{k}) = h \prod_{\alpha_{p_{k},e}(\xi) < 0} \alpha_{p_{k},e}.$$
If $\sigma(p_{k}) > m$, then $f_k(p_{k}) =0$ 
and if $\sigma(p_{k}) \leq m$ let 
$$f_{k+1} = f_k - h_{p_{k}}\tau_{p_{k}} \in H^{2m}(\Gamma, \alpha).$$
Then $f_{k+1}$ is supported on $\{p_{k+1},..,p_N \}$. Proceeding 
inductively we conclude that 
$$f_N = f -\sum_{\sigma(p) \leq m} h_p\tau_p$$
is zero at all vertices, that is, that (\ref{eq:6.5}) holds.
\end{proof}

\subsection{\bf The Kirwan map}
\label{sec:kirwanmap}

\subsubsection{\bf The Kirwan map}
\label{ssec:kirwanmap}

Let $(\Gamma,\alpha)$ be as in Theorem \ref{th:addstruct} 
a one-skeleton which is 3-independent and non-cyclic
and let $\phi :V_{\Gamma} \to \RR$
be an injective function which is 
$\xi$-compatible for some $\xi \in \P$. 
Assume that the conditions of Theorem \ref{th:addstruct}
are satisfied.

Let
${\mathcal F}^{2k}(\Gamma_c, \alpha^c)$ be the set of all maps
\begin{displaymath}
  f:V_{\Gamma_c} \to \SS^k (\g^*_{\xi}) \, .
\end{displaymath}
The sum
$${\mathcal F}(\Gamma_c, \alpha^c) = 
\bigoplus {\mathcal F}^{2k}(\Gamma_c, \alpha^c)$$
is a graded ring under point-wise multiplication and by 
Theorem \ref{th:kirate} one gets a map
$$\K_c : H(\Gamma,\alpha) \to {\mathcal F}(\Gamma_c,\alpha^c).$$

\begin{theorem}
\label{th:Kirwan map}
$\K_c$ maps $H(\Gamma,\alpha)$ to $H(\Gamma_c, \alpha^c)$.
\end{theorem}

\begin{proof}

All we need to show is that the image of the map, $\K_c$, is indeed in
$H(\Gamma_c, \alpha^c)$, that is that $\K_c(f)$ satisfies the
compatibility conditions (\ref{eq:coh6}) for the reduced one-skeleton.

Let $\{ x,y_1,...,y_{n-1} \}$ be a basis of $\fg^*$ such that 
$x(\xi)=1$ and $\{ y_1,...,y_{n-1} \}$ is a basis of $\fg_{\xi}^*$.
Let $\alpha = \alpha(\xi)(x-\beta(y)) \in \fg^*$ such that
$\alpha(\xi) \ne 0$. Then the map 
$$\rho_{\alpha}: \SS(\fg^*) \rightarrow \SS(\fg_{\xi}^*)$$
given by the 
identification $\fg_{\xi}^* \simeq \fg^*/ \RR\alpha$ will send 
$x$ to $x-\alpha/\alpha(\xi)\in \fg_{\xi}^*$ and $y_j$ to $y_j$.
Therefore $\rho_{\alpha}$ will send a polynomial 
$P(x,y) \in \SS(\fg^*)$ to the polynomial 
$P(x-\alpha/\alpha(\xi),y) = P(\beta(y),y)\in \SS(\fg_{\xi}^*)$.

Now let $f \in H(\Gamma,\alpha)$. With the notations in 
Figure \ref{fig:red} we will show that
$\K_c(f)$ satisfies (\ref{eq:coh6}) for the edge $v_iv_{ia}$: 
\begin{equation}
\label{eq:7.1}
\K_c(f)(v_{ia}) - \K_c(f)(v_i) \equiv 0 \pmod{\alpha_{v_iv_{ia}}^c}
\quad \mbox{ in } \SS(\fg_{\xi}^*).
\end{equation}

For each $j$ we have 
\begin{equation*}
\rho_{\alpha_{t_jt_{j+1}}}(f_{t_j}) = 
f_{t_j}( x- \alpha_{t_jt_{j+1}}/ \alpha_{t_jt_{j+1}}(\xi), y )
\end{equation*}
and therefore the difference 
$$\rho_{\alpha_{t_jt_{j+1}}}(f_{t_j}) 
-  \rho_{\alpha_{t_jt_{j-1}}}(f_{t_j})$$
is the same as
$$f_{t_j}( x- \frac{\alpha_{t_jt_{j+1}}}{ \alpha_{t_jt_{j+1}}(\xi)}, y ) 
- f_{t_j}( x- \frac{\alpha_{t_jt_{j-1}}}{ \alpha_{t_jt_{j-1}}(\xi)}, y ),$$
which is a multiple (in $S(\fg_{\xi}^*)$) of 
\begin{equation*}
(x-\frac{\alpha_{t_jt_{j+1}}}{ \alpha_{t_jt_{j+1}}(\xi)} )
- (x-\frac{\alpha_{t_jt_{j-1}}}{ \alpha_{t_jt_{j-1}}(\xi)} )
\end{equation*}
and, thus, of $\alpha_{v_iv_{ia}}^c$.
Hence (\ref{eq:7.1}) follows from the fact that
\begin{equation*}
\K_c(f)(v_{ia}) - \K_c(f)(v_i) = \sum_{j=1}^k
\left( \rho_{\alpha_{t_jt_{j+1}}}(f_{t_j}) - 
\rho_{\alpha_{t_jt_{j-1}}}(f_{t_j}) \right).
\end{equation*}
\end{proof}

\subsubsection{\bf The kernel of the Kirwan map}
\label{ssec:kernel}

The following theorem, which  describes the kernel of the map above, 
is the combinatorial analogue of a result of Tolman and Weitsman 
(\cite{TW}).

\begin{theorem}
\label{th:kerkirw}
The kernel of the map 
$\K_c : H(\Gamma, \alpha) \to H(\Gamma_c, \alpha^c)$
consists of those elements $f \in H(\Gamma,\alpha)$ which can be written 
as a sum
$f=h_+ + h_-,$
with $h_{\pm} \in H(\Gamma,\alpha)$ such that $h_+$ is supported on 
$\phi >c$ and $h_-$ is supported on $\phi < c$.
\end{theorem}

\begin{proof}
If $f$ is in the kernel of $\K_c$ then 
$\K_c(f)(v) = 0$ for every edge $(p,q)$ with $\phi(p) < c < \phi(q)$, where
$v \in V_{\Gamma}$ is the vertex that corresponds to this edge of $\Gamma$.
Since 
$$0 = \K_c(f)(v) = \rho_{\alpha_{pq}} (f_p)  \;\; ,$$
it follows that $f_p$ is divisible by $\alpha_{pq}$. Similarly $f_q$ 
is divisible by $\alpha_{pq}$.

Consider now the maps, $h_{\pm}$, of $V_{\Gamma}$ into  $\SS(\fg^*))$ 
defined by
$$h_-(p)=\left\{ \begin{array}{ccc} 
f(p), & \mbox{ if } & \phi(p) < c \\
0, & \mbox{ if } & \phi(p) > c \end{array} \right. ,
\qquad
h_+(q)=\left\{ \begin{array}{ccc} 
0, & \mbox{ if } & \phi(q) < c \\
f(q), & \mbox{ if } & \phi(q) > c \end{array} \right. \; .
$$
It is clear that $f=h_+ + h_-$ and that $h_{\pm} \in H(\Gamma,\alpha)$.
\end{proof}

By a slight modification of the proof of Theorem \ref{th:freemodule},
one can prove the following:

\begin{corollary}
\label{cor:6.3.1}
The dimension of $ker\{\K_c : H^{2m}(\Gamma) \to H^{2m}(\Gamma_c)\}$ is 
\begin{equation*}
\dim{H_{c^+}^{2m}(\Gamma)} + \dim{H_{c^-}^{2m}(\Gamma)} = 
\sum_{\phi(q) > c} \lambda_{m-\sigma(q),n} +
\sum_{\phi(p)<c} \lambda_{m-d+\sigma(p),n} \; ,
\end{equation*}
where
$$H_{c^+}(\Gamma) = 
\{ h \in H(\Gamma, \alpha) ; h \mbox{ is supported on } \phi > c \}$$
and 
$$H_{c^-}(\Gamma) = 
\{ h \in H(\Gamma,\alpha) ; h \mbox{ is supported on } \phi < c \}.$$
\end{corollary}

\subsubsection{\bf The surjectivity of the Kirwan map}
\label{ssec:surjectivity}

We will finally prove the graph theoretic analogue of Theorem 
\ref{th:kirw}.

\begin{theorem}
\label{th:surjkirw}
For generic $\xi \in \P$, the Kirwan map $\K_c$ is surjective.
\end{theorem}

\begin{proof}
We will show that 
$\K_c : H(\Gamma, \alpha) \to H(\Gamma_c, \alpha^c)$
is surjective by a dimension count, using induction on the number
of vertices $p\in V_{\Gamma}$ that lie under the level $\phi=c$.

To start, assume there is only one such vertex, $p$. Then $p$ is minimum 
of $\phi$ and the reduced space $\Gamma_c$ is 
a complete one-skeleton with $d$ vertices, $v_1,..,v_d$, 
one for each edge $e_i$ issuing from $p$. 
Let $f \in H^{2m}(\Gamma_c, \alpha^c)$. Then, by (\ref{eq:4.5}),
there exists $f_0 \in \SS(\fg^*) \subset H(\Gamma,\alpha)$ 
such that 
$\rho_{\alpha_{e_i}}(f_0) = f(v_i)$ for all $i$'s. Hence
$\K_c(f_0) =f.$

Let's assume now that the map, $\K_c$, is surjective at the level $a$ and 
choose a regular value $b>a$ such that there is exactly one vertex, $p$,
with $a < \phi(p) < b$. Let $r=\sigma(p)$ and $s = d-r$. 
We have two exact sequences
\begin{equation}
\label{eq:7.2}
0 \to ker (\K_a) \to H^{2m}(\Gamma,\alpha) 
\stackrel{\K_a}{\longrightarrow} H^{2m}(\Gamma_a, \alpha^a) 
\to 0,
\end{equation}
and 
\begin{equation*}
0 \to ker (\K_b) \to H^{2m}(\Gamma,\alpha) 
\stackrel{\K_b}{\longrightarrow} H^{2m}(\Gamma_b, \alpha^b).
\end{equation*}
The last arrow in (\ref{eq:7.2}) is surjective because of our
inductive assumption. 

By Corollary \ref{cor:6.3.1} 
\begin{equation}
\label{eq:7.3}
\dim{ker(\K_b)} - \dim{ker(\K_a)} = \lambda_{m-s,n} - \lambda_{m-r,n}
\end{equation}
and by (\ref{eq:7.3}) and (\ref{eq:6.1new}) 
\begin{multline*}
 \dim{ker(\K_b)}  +  \dim{H^{2m}(\Gamma_b)} = 
\dim{ker(\K_a)}+ \dim{H^{2m}(\Gamma_a)}  \\
 + (\lambda_{m-s,n}  -  \lambda_{m-r,n})+
(\lambda_{m-r,n} - \lambda_{m-s,n}) = \dim H^{2m}(\Gamma,\alpha).
\end{multline*}
This proves that 
\begin{equation*}
\dim({im(\K_b)}) = \dim{H^{2m}(\Gamma)} - \dim{ker(\K_b)} = 
\dim{H^{2m}(\Gamma_b)},
\end{equation*}
hence that $\K_b$ is surjective.
\end{proof}

\section{\bf Applications}
\label{sec:applic}

\subsection{\bf The realization theorem}
\label{ssec:realization}

Recall that an abstract one-skeleton, $(\Gamma,\alpha)$, 
is an abstract GKM one-skeleton if $\alpha$ satisfies \eqref{eq:axiom2gkm} and 
\eqref{eq:axiom3gkm} and the
constants $c_{i,e}$ in \eqref{eq:axiom3gkm} are integers. 
In this section we will 
prove that all such abstract one-skeleta can be realized as the GKM-skeleta of 
GKM spaces.

\begin{theorem}
\label{th:realization}
If $(\Gamma,\alpha)$ is an abstract GKM one-skeleton, 
there exists a complex 
manifold $M$ and a GKM action of $G$ on $M$ for which 
$(\Gamma, \alpha)$ is its GKM one-skeleton.
\end{theorem}

\noindent {\bf Remarks:}
\begin{enumerate}
\item The manifold $M$ which we will construct below is not compact, and 
there does not appear to be a canonical compactification of it. For some 
interesting non-canonical compactifications of it see \cite{GKT}.
\item The manifold $M$ is also not equivariantly formal, but it does have 
the property that the canonical map of $H_G(M)$ into $H(\Gamma,\alpha)$
is surjective.
\end{enumerate}

\begin{proof}
Our construction of $M$ will involve three steps: first we will construct the 
$\CC P^1$'s corresponding to the edges of $\Gamma$; then, for each of these 
$\CC P^1$'s, we will construct a tubular neighborhood of it in $M$. Then we 
will construct $M$ itself by gluing these tubular neighborhoods together.

Let $\rho_{\alpha}$ be the one dimensional representation of $G$ with
weight $\alpha$ and let $V_{\alpha} \simeq \CC$ be the vector space on which 
this representation lives. Let $G$ act on $V_{\alpha} \oplus \CC$ by acting
by $\rho_{\alpha}$ on the first factor and by the trivial representation on 
the second factor. This action induces an action of $G$ on the 
projectivization
\begin{equation}
\label{eq:real1}
X_{\alpha} = \CC P^1 = \PP(V_{\alpha} \oplus \CC).
\end{equation}
The points $q=[1:0]$ and $p=[0:1]$ are the two fixed points of this action
and there are equivariant bijective maps
\begin{equation}
\label{eq:real2}
V_{\alpha} \to X_{\alpha}-\{ q \}, \qquad c \to [c:1],
\end{equation}
and
\begin{equation}
\label{eq:real3}
V_{-\alpha} \to X_{\alpha}-\{ p \}, \qquad c \to [1:c].
\end{equation}
(The equivariance of \eqref{eq:real3} follows from the fact that , 
for $\xi \in \fg$,
$$[e^{i\alpha(\xi)}:c]=[1:e^{-i\alpha(\xi)} c].) $$

We will denote by $\LL_{\alpha}$ the tautological line bundle over 
$X_{\alpha}$. By definition, the fiber of $\LL_{\alpha}$ over $[c_1:c_2]$
is the one dimensional subspace of $V_{\alpha} \oplus \CC$ spanned by 
$(c_1,c_2)$; so from the action of $G$ on $V_{\alpha} \oplus \CC$ one gets
an action of $G$ on $\LL_{\alpha}$ lifting the action of $G$  on 
$X_{\alpha}$. The fiber of $\LL_{\alpha}$ over $q$ is $V_{\alpha}$, so, 
in particular, the following is true:

\begin{lemma}
\label{lem:real1}
The weight of the isotropy action of $G$ on $(\LL_{\alpha})_p$ is zero 
and on $(\LL_{\alpha})_q$ is $\alpha$.
\end{lemma}

The mapping 
\begin{equation}
\label{eq:real4}
[c:1] \to (c,1)
\end{equation}
defines a holomorphic section of $\LL_{\alpha}$ over $X_{\alpha}-\{q\}$, 
and hence a holomorphic trivialization of the restriction of $\LL_{\alpha}$ to 
$X_{\alpha}-\{q\}$.

The vector space $V_{\alpha} \oplus \CC \simeq \CC^2$ can be equipped with
the $G$-invariant Hermitian form
\begin{equation}
\label{eq:real5}
|z|^2 = |z_1|^2 + |z_2|^2
\end{equation}
and since the restriction of this form to each subspace of $\CC^2$ defines 
a Hermitian form on this subspace, we get from this form a $G$-invariant 
Hermitian structure on $\LL_{\alpha}$.

Now let $\alpha_i$ and $\alpha_i'$ be weights of $G$ with 
$\alpha_i' - \alpha_i = m \alpha,$ 
$m$ being an integer, and let $\LL_i$ be the line bundle 
\begin{equation}
\label{eq:real6}
\LL_i = \LL_{\alpha}^m \otimes V_{\alpha_i}.
\end{equation}
From the action of $G$ on $\LL_{\alpha}$ and on $V_{\alpha_i}$ one gets an 
action of $G$ on this line bundle lifting the action of $G$ on $X_{\alpha}$. 
The following is a corollary of Lemma \ref{lem:real1} and of the existence of 
the Hermitian structure \eqref{eq:real5} on $\LL_{\alpha}$ and of the 
trivialization \eqref{eq:real4} of $\LL_{\alpha}$ over $X_{\alpha} - \{q\}$.

\begin{lemma}
\label{lem:real2}
The weight of the isotropy representation of $G$ on $(\LL_i)_p$ is $\alpha_i$ 
and on $(\LL_i)_q$ is $\alpha_i'$. In addition, $\LL_i$ has a $G$-invariant
Hermitian structure and a non-vanishing holomorphic section, 
$s_i : X_{\alpha}-\{q\} \to \LL_i$, which transforms under the action of $G$
according to the weight $\alpha_i$. In particular, the restriction of $\LL_i$
to $X_{\alpha}-\{q\}$ is isomorphic to the trivial bundle over 
$X_{\alpha}-\{q\}$ with fiber $V_{\alpha_i}$.
\end{lemma}

Let's now return to the problem of constructing a manifold $M$ with 
$(\Gamma,\alpha)$ as its GKM one-skeleton. Let $e$ be an oriented edge of 
$\Gamma$ and let $p,p',e_i$ and $e_i'$ be as in \eqref{eq:axiom3gkm} 
with $e_d=e$ and $e_d'=\bar{e}$. Let $X_e=X_{\alpha}$, with
$\alpha=\alpha_e$, and let $\LL_i$ be the line bundle constructed above with
$\alpha_i=\alpha_{e_i}$ and $\alpha_i'=\alpha_{e_i'}$. The $X_e$'s will 
be our candidates for the $G$-invariant $\CC P^1$'s in $M$ and the vector
bundle 
\begin{equation}
\label{eq:real7}
N_e = \bigoplus_{i=1}^{d-1} \LL_i
\end{equation}
will be our candidate for the normal bundle of $X_e$ in $M$. Thus, a 
candidate for a tubular neighborhood of $X_e$ in $M$ will be a convex 
neighborhood of the zero section in $N_e$, for example, the disk bundle
\begin{equation}
\label{eq:real8}
\U_e^{\epsilon} = \{ (x,v_1,..,v_{d-1}) ; x \in X_e, v_i \in (\LL_i)_x, 
|v_1| < \epsilon \}.
\end{equation}

We will construct $M$ by starting with the disjoint union
\begin{equation}
\label{eq:real8'}
\coprod \U_e^{\epsilon}
\end{equation}
over all edges $e$ of $\Gamma$, and making the following obvious 
identifications: Let $N_{e,p}$ be the restriction of the bundle $N_e$ to 
$X_e - \{q \}$. By Lemma \ref{lem:real2} and by \eqref{eq:real2} one gets a
$G$-equivariant bijective map
\begin{equation}
\label{eq:real9}
N_{e,p} \stackrel{\gamma_e}{\longrightarrow} T_pM
\end{equation}
where $T_pM$ is by definition the sum 
\begin{equation}
\label{eq:real10}
\bigoplus_{i=1}^d V_{\alpha_i}.
\end{equation}

If $e$ and $e'$ are edges of $\Gamma$ meeting at $p$ we will identify the 
points $u \in \U_{e,p}^{\epsilon}$ and $u' \in \U_{e',p}^{\epsilon}$
in the set \eqref{eq:real8'} if
$$ \gamma_e(u) = \gamma_{e'}(u').$$
It is easy to check that if we quotient \eqref{eq:real8'} by the equivalence 
relation defined by these identifications (with $\epsilon$ small), we get 
a manifold $M$ with the properties listed in theorem.
\end{proof}

\subsection{\bf A deformation problem}
\label{sec:polytopes}

We return to the example of section \ref{ssec:polytopes}: 
the one-skeleton, $\Gamma$, of 
an edge-reflecting $d$-valent convex polytope $\Delta$ embedded into an 
$n$-dimensional space 
$\fg^*$ by 
$\Phi : \Delta \to \fg^*$, the axial function, $\alpha$, of $\Gamma$ 
being given by
$$\alpha_{pq} = \Phi(q)-\Phi(p).$$

Fix a vector $\xi \in \P$ and let $\phi : V_{\Gamma} \to \RR$ be given by 
\begin{equation}
\label{eq:fixi}
\phi(p) = \langle \Phi(p) , \xi \rangle \; \quad \mbox{ for all } 
p \in V_{\Gamma}.
\end{equation}

Then $\phi$ is $\xi$-compatible and 
the zeroth Betti number of $\Gamma$ is 1;
and since every 2-face of $\Delta$ is convex, $\Gamma$ is non-cyclic.
Moreover, $\Gamma$ is 3-independent if $\dim \fg^* \geq 3$ (see the 
comments at the end of section \ref{ssec:polytopes}).

Assume that there exists a lattice $\ZZ_G^*$ in $\fg^*$ such that the edges
of $\Delta$ are scalar multiples of rational vectors. If we try to deform 
$\Delta$ by changing its vertices such that the above property 
is preserved, we are led to the following definition:

\begin{definition}
A function $f: V_{\Gamma} \to \fg^*$ is called a 
\emph{rational deformation} of $\Phi$ if there exists $\epsilon > 0$
such that for every $t \in [0,\epsilon)$, the map
$\Phi_t : V_{\Gamma} \to \fg^*$ given by
\begin{equation*}
\Phi_t(p) = \Phi(p) + t f(p)\; , \forall p \in V_{\Gamma} \; ,
\end{equation*}
is an embedding of $V_{\Gamma}$ into $\fg^*$ and, for every edge
$e=(p,q)$ of $\Gamma$, $\Phi_t(q)-\Phi_t(p)$ is a positive multiple
of $\alpha_{p,e}$.
\end{definition}

\begin{theorem}
\label{th:accdef}
For a given embedding $\Phi$, the space of rational deformations is
$H^2(\Gamma,\alpha)$.
\end{theorem}

\begin{proof}
Let $f$ be a rational deformation. Then 
\begin{equation}
\label{eq:phit}
\Phi_t(q)-\Phi_t(p) = \Phi(q)-\Phi(p) + t (f(q)-f(p));
\end{equation}
since $\alpha_{pq}$ divides both $\Phi_t(q)-\Phi_t(p)$ and 
$\Phi(q)-\Phi(p)$, it follows that it divides $f(q)-f(p)$ as well,
which means that $f \in H^2(\Gamma,\alpha)$.

Conversely, if $f \in H^2(\Gamma,\alpha)$ then (\ref{eq:phit})
implies that $\Phi_t(q)-\Phi_t(p)$ is a multiple of $\alpha_{pq}$
and, since $\Phi(q)-\Phi(p) = \alpha_{pq}$, for $t$ small enough, it is 
a positive multiple. We can choose $\epsilon$ small enough 
for all edges, which proves that $f$ is a rational deformation.
\end{proof}

By (\ref{eq:dimhg}) 
\begin{equation}
\label{eq:dimh2}
\dim H^2(\Gamma,\alpha) = 
b_2(\Gamma)\lambda_0 + b_0(\Gamma)\lambda_1 = b_2(\Gamma) + n.
\end{equation}
Every translation is a rational deformation and 
the ``$n$'' in (\ref{eq:dimh2}) is the contribution of these ``trivial'' 
deformations to the space of deformations of $\Gamma$.
Hence the non-trivial rational deformations are those corresponding to 
$b_2(\Gamma)$. By Theorem \ref{th:freemodule}, 
these deformations are linear 
combinations of Thom classes, $\tau_p$, for $p$ of index 1.

\begin{example}
\label{eq:octohedron}
Consider an octahedron embedded in $\RR^3$; then $b_2(\Gamma)=1$.
All rational deformations are obtained by composing translations 
with the homothety $p \in V_{\Gamma} \to (1-t)p$. In Figure 
\ref{fig:octahedron} 
the numbers next to 
vertices indicate the index with respect to the height function.
The condition on 2-planes amounts to saying that 04, 22 and 13 
intersect in a point.

\begin{figure}[h]
\begin{center}
\includegraphics{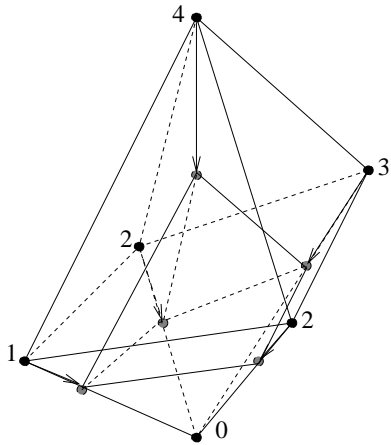}
\caption{Rational deformation of the octahedron}
\label{fig:octahedron}
\end{center}
\end{figure}
\end{example}

\subsection{\bf Schubert polynomials}
\label{ssec:schubs}

For the Grassmannians, the classes $\tau_p$ of Theorem \ref{th:suppthom}
have an alternative description in terms of Schubert polynomials 
(see \cite{BGG}, \cite{LS}, \cite{Dem}, \cite{KK}, \cite{Mac},
\cite{BH}, \cite{Fu1} et al.). This description involves the 
Hecke algebra of divided difference operators, an algebra which
is intrinsically associated to every compact semisimple Lie group $K$: 
Let $G$ be the Cartan subgroup of $K$ and $W=N(G)/G$ the Weyl group. As a 
group of transformations of $\fg$, $W$ is generated by simple 
reflections. Moreover, to each reflection $\sigma \in W$ corresponds a
unique positive root $\alpha=\alpha_{\sigma}\in \fg^*$ with
$\alpha(\sigma\xi)=-\alpha(\xi)$ for all $\xi \in \fg$. In particular,
$\sigma$ leaves fixed the hyperplane $\alpha(\xi)=0$. The \emph{divided 
difference operator}
$$D_{\sigma} : \SS(\fg^*) \to \SS(\fg^*)$$
is the operator defined by
\begin{equation}
\label{eq:schub1}
D_{\sigma}(f)= \frac{f-\sigma f}{\alpha_{\sigma}}.
\end{equation}
(Notice that since $\sigma f (\xi)=f(\sigma(\xi))=f(\xi)$ for $\xi$ on 
the hyperplane $\alpha(\xi)=0$, the left hand side of (\ref{eq:schub1})
is an element of $\SS(\fg^*)$, that is, a polynomial function on $\fg$.)

\emph{The Hecke algebra of divided difference operators} 
$\DD$ is the algebra 
generated by the $D_{\sigma}$'s and the operators ``multiplication by
$f$'' for $f \in \SS(\fg^*)$. We note that if $g \in \SS(\fg^*)^W$ then
$$D_{\sigma}(gf)= \frac{gf-\sigma(gf)}{\alpha_{\sigma}} = 
g \frac{f-\sigma f}{\alpha_{\sigma}} = g D_{\sigma}f,$$
hence, if $D \in \DD$, then 
\begin{equation}
\label{eq:schub2}
D(gf)=gDf,
\end{equation}
so the algebra $\DD$ acts on $\SS(\fg^*)$ as morphisms of 
$S(\fg^*)^W$-modules. More generally, if $\MM_0$ is an 
$S(\fg^*)^W$-module and 
\begin{equation}
\label{eq:schub3}
\MM = \MM_0 \otimes_{S(\fg^*)^W} \SS(\fg^*)
\end{equation}
then one can make $\MM$ into a $\DD$-module by setting
\begin{equation}
\label{eq:schub4}
D( m \otimes f) = m \otimes (Df).
\end{equation}
(In view of (\ref{eq:schub2}) this is a well-defined operator on $\MM$.)

Now let $M$ be a $K$-manifold. From the constant map $M \to pt$ one gets
a map in cohomology
$$H_K(pt) \to H_K(M)$$
and since 
$$H_K(pt) = \SS(\fk^*)^K = \SS(\fg^*)^W,$$
this map makes $H_K(M)$ into a module over $\SS(\fg^*)^W$. For the
following result see, for instance, \cite[Ch. 6]{GS3}:

\begin{theorem}
\label{th:gcohm}
The $G$-equivariant cohomology ring of $M$ is related to the 
$K$-equivariant cohomology ring of $M$ by the following 
ring-theoretic identity:
\begin{equation}
\label{eq:schub5}
H_G(M)= H_K(M) \otimes_{\SS(\fg^*)^W} \SS(\fg^*).
\end{equation}
\end{theorem}

Therefore, by (\ref{eq:schub3}) and (\ref{eq:schub4}) we conclude:

\begin{theorem}
\label{th:hgmd}
The $G$-equivariant cohomology ring $H_G(M)$ is canonically a module
over $\DD$.
\end{theorem}

Suppose now that $M$ is a GKM manifold and is equivariantly formal. 
Then, by Theorem \ref{th:gkm}, 
$$H_G(M) \simeq H(\Gamma,\alpha);$$
so, we can transport this $\DD$-module structure to $H(\Gamma,\alpha)$.
In particular, we get an action of the divided difference operator
$D_{\sigma}$ on $H(\Gamma,\alpha)$.

\begin{theorem}
\label{th:dongamma}
Let $f : V_{\Gamma} \to \SS(\fg^*)$ be a map which satisfies the 
compatibility conditions (\ref{eq:coh6}) (in other words, which belongs to 
$H(\Gamma,\alpha)$). Then, for $p \in V_{\Gamma}$,
\begin{equation}
\label{eq:schub6}
(D_{\sigma}f)(p) = \frac{f(p)-\sigma(f(\sigma^{-1}p))}{\alpha_{\sigma}}.
\end{equation}
\end{theorem}

\noindent {\bf Remarks:}

\begin{enumerate}
\item Since $W$ is by definition $N(G)/G$, it acts on the fixed point set
$M^G$. Since $M^G=V_{\Gamma}$, the expression $\sigma^{-1}p$ on the left 
hand side is unambiguously defined.

\item Since $W$ acts on $\SS(\fg^*)$ by ring automorphisms, the ring 
automorphism $\sigma$, applied to the element 
$f(\sigma^{-1}p) \in \SS(\fg^*)$ on the left hand side is also 
unambiguously defined.

\item For a proof of Theorem \ref{th:dongamma} see \cite{GHZ}.
\end{enumerate}

We now return to section \ref{ssec:grass} and the Bruhat structure of the 
Johnson graph. Let $\phi$ be the Morse function (\ref{eq:grass12})
and let $p \in V_{\Gamma}$ be a vertex of $\Gamma$ of index $r$. Let
$p_0$ be the unique maximum of $\phi$, that is 
$$ p_0 =p_S \quad , \quad S= \{ l+1,..,n \}$$
with $l=n-k$ and let 
\begin{equation}
\label{eq:schub7}
\Delta = \prod_{e \in E_{p_0}} \alpha_e =
\prod_{i \leq k < j}(\alpha_i - \alpha_j).
\end{equation}

By Theorem \ref{th:suppthom}, the Thom class $\tau_{p_0}$ is the map
$$\tau_{p_0} : V_{\Gamma} \to \SS(\fg^*) $$
which takes the value $\Delta$ at $p_0$ and 0 everywhere else. Now let
$\sigma_{i_1}, ..., \sigma_{i_s}$ ($s=kl-r$) be the elementary reflections 
with the properties described in Theorem \ref{th:Bruhatorder}. By 
applying the operator 
$$D_{\sigma_{i_1}} \circ ... \circ D_{\sigma_{i_s}}$$
to $\tau_{p_0}$ one gets a cohomology class $\tau_p'$, which is of degree
$r=kl-s$ and which, by (\ref{eq:grass17}), (\ref{eq:grass18}) and
(\ref{eq:schub6}), is supported on $F_p$. Moreover, it is easy to see 
that 
$$\tau_p(p) = \prod \alpha_e,$$
the product being over the down-ward pointing edges $e \in E_p$. Thus, by
Theorem \ref{th:unictaup} and (\ref{eq:grass12}), $\tau_p= \tau_p'$.

\begin{remark*}
One can regard the $\tau_p$'s as being a doubly indexed family of 
polynomials $f_{p,q} = \tau_p(q)$, indexed by pairs $(p,q)$ of
vertices of $\Gamma$ with $p \prec q$. These polynomials, which are 
called \emph{double Schubert polynomials}, have been studied extensively 
by Billey, Haiman, Stanley, Fomin, Kirillov, Jockusch and others.
(For a succinct and engrossing account of what is known about these 
polynomials we recommend, as collateral reading, the beautiful monograph
`` Young Tableaux'' by William Fulton, which has just been published by
Cambridge University Press.)
\end{remark*}

\end{document}